\documentclass[11pt,reqno]{amsart}

\usepackage[english]{babel}
\usepackage[utf8]{inputenc}   
\usepackage[osf,sc]{mathpazo}
\usepackage[cal=pxtx,scr=dutchcal,bb=pazo]{mathalfa}  %

\usepackage{indentfirst}        
\usepackage{tabularx,booktabs}      
\usepackage{multirow}
\usepackage{caption,subcaption}     
\usepackage[english]{varioref}      
\usepackage[dvipsnames]{xcolor}     
\usepackage{hyperref}         
\usepackage{bookmark}         
\usepackage{textcomp}
\usepackage{enumitem} 
\usepackage{fontawesome5}
\usepackage[margin=1in]{geometry}

\setlist[itemize,1]{leftmargin=1.5em} 
\setlist[enumerate,1]{leftmargin=1.5em} 
\makeatletter
  \def\l@subsection{\@tocline{2}{0pt}{2.5pc}{5pc}{}} 
\makeatother

\usepackage[normalem]{ulem}
\newcommand{\stkout}[1]{\ifmmode\text{\sout{\ensuremath{#1}}}\else\sout{#1}\fi}

\definecolor{highlight}{cmyk}{90,99,0,0}
\hypersetup{
colorlinks=true, linktocpage=true, pdfstartpage=1, pdfstartview=FitV,breaklinks=true, pdfpagemode=UseNone, pageanchor=true, pdfpagemode=UseOutlines,plainpages=false, bookmarksnumbered, bookmarksopen=true, bookmarksopenlevel=1,hypertexnames=true, pdfhighlight=/O,urlcolor=highlight, linkcolor=highlight, citecolor=ForestGreen}

\usepackage{graphicx,float}   
\usepackage{wrapfig}
\usepackage{tikz,tikz-3dplot,tikz-cd,pgfplots}      
\usetikzlibrary{
  arrows,decorations.pathreplacing,decorations.markings,shapes.geometric,positioning,calc,fadings,angles,intersections,3d,perspective,patterns,backgrounds,decorations.pathmorphing
}
\pgfplotsset{compat=1.18}

\makeatletter
\def\@tocline#1#2#3#4#5#6#7{\relax
  \ifnum #1>\c@tocdepth 
  \else
  \par \addpenalty\@secpenalty\addvspace{#2}%
  \begingroup \hyphenpenalty\@M
  \@ifempty{#4}{%
    \@tempdima\csname r@tocindent\number#1\endcsname\relax
  }{%
    \@tempdima#4\relax
  }%
  \parindent\z@ \leftskip#3\relax \advance\leftskip\@tempdima\relax
  \rightskip\@pnumwidth plus4em \parfillskip-\@pnumwidth
  #5\leavevmode\hskip-\@tempdima
  \ifcase #1
  \or\or \hskip 1em \or \hskip 2em \else \hskip 3em \fi%
  #6\nobreak\relax
  \dotfill\hbox to\@pnumwidth{\@tocpagenum{#7}}\par
  \nobreak
  \endgroup
  \fi}
\makeatother

\usepackage[noabbrev]{cleveref}
\expandafter\def\csname ver@etex.sty\endcsname{3000/12/31}

\crefname{lemma}{lemma}{lemmata}
\Crefname{lemma}{Lemma}{Lemmata}
\crefname{exercise}{exercise}{exercises}
\Crefname{exercise}{Exercise}{Exercises}
\crefname{subsection}{subsection}{subsections}
\Crefname{subsection}{Subsection}{Subsections}
\crefname{conjecture}{conjecture}{conjectures}
\Crefname{conjecture}{Conjecture}{Conjectures}
\crefname{question}{question}{questions}
\Crefname{question}{Question}{Questions}
\crefname{warning}{warning}{warnings}
\Crefname{warning}{Warning}{Warnings}

\usepackage{amsmath,amssymb,amsthm,mathtools}
\usepackage{braket,faktor,stmaryrd,bm}
\numberwithin{equation}{section}

\newcommand{\dd}{\mathrm{d}}
\newcommand{\de}{\partial}
\newcommand{\dde}[2]{\frac{\partial #1}{\partial #2}}
\newcommand{\bbraket}[1]{\llbracket #1 \rrbracket}

\newcommand{\Id}{\mathrm{Id}}
\newcommand{\iu}{\mathrm{i}}

\newcommand{\rk}{\mathrm{rk}}

\newcommand{\diag}{\mathrm{diag}}

\newcommand{\cc}{\mathrm{c}}
\newcommand{\ch}{\mathrm{ch}}

\newcommand{\bigO}{\mathrm{O}}

\newcommand{\End}{\mathrm{End}}
\newcommand{\Aut}{\mathrm{Aut}}

\newcommand{\mb}[1]{\mathbb{#1}}
\newcommand{\mc}[1]{\mathcal{#1}}
\newcommand{\mf}[1]{\mathfrak{#1}}
\newcommand{\mr}[1]{\mathrm{#1}}
\newcommand{\ms}[1]{\mathsf{#1}}
\newcommand{\msc}[1]{\mathscr{#1}}

\newcommand{\Z}{\mb{Z}}
\newcommand{\Q}{\mb{Q}}
\newcommand{\R}{\mb{R}}
\newcommand{\C}{\mb{C}}
\renewcommand{\P}{\mb{P}}
\renewcommand{\H}{\mb{H}}

\newcommand{\M}{\mc{M}}
\newcommand{\Mbar}{\overline{\M}}
\newcommand{\Mhyp}{\M^{\textup{hyp}}}

\newcommand{\Vir}{{\ms{Vir}_{\ge -1}}}

\newcommand{\VWP}{V^{\textup{WP}}}
\newcommand{\WP}{\omega_{\textup{WP}}}

\DeclareMathOperator{\tr}{tr}

\theoremstyle{plain}
\newtheorem{theorem}{Theorem}[section]
\newtheorem{proposition}[theorem]{Proposition}

\newtheorem{exercise}{Exercise}[section]

\theoremstyle{definition}
\newtheorem{definition}[theorem]{Definition}

\newtheorem{example}[theorem]{Example}

\newtheorem*{terminology}{Terminology}

\usepackage[
    style = alphabetic,%
    maxalphanames = 5,%
    giveninits = true,%
    sorting = nyt,%
    maxbibnames = 99,%
    backend = bibtex%
  ]{biblatex}
\renewbibmacro{in:}{} 
\usepackage{csquotes}
\title[Moduli spaces of Riemann surfaces]{Les Houches lecture notes on moduli spaces of Riemann surfaces}

\author[A.~Giacchetto]{Alessandro Giacchetto}
\address[A.~Giacchetto]{
  ETH Zürich, Departement Mathematik, Zürich, Switzerland
}
\email{alessandro.giacchetto@math.ethz.ch}

\author[D.~Lewa{\'{n}}ski]{Danilo Lewa{\'{n}}ski}
\address[D.~Lewa{\'{n}}ski]{
  Università degli Studi di Trieste, Dipartimento MIGe \& INFN, Sezione di Trieste, Trieste, Italy
}
\email{danilo.lewanski@units.it}

\begin{document}

\begin{abstract}
  In these lecture notes, we provide an introduction to the moduli space of Riemann surfaces, a fundamental concept in the theories of 2D quantum gravity, topological string theory, and matrix models. We begin by reviewing some basic results concerning the recursive boundary structure of the moduli space and the associated cohomology theory. We then present Witten's celebrated conjecture and its generalisation, framing it as a recursive computation of cohomological field theory correlators via topological recursion. We conclude with a discussion of JT gravity in relation to hyperbolic geometry and topological strings.
\end{abstract}

\maketitle

\vspace{-.5cm}

\setcounter{tocdepth}{2}
\tableofcontents

\newpage
\section{Introduction}
\label{sec:intro}

As a physical theory, 2D gravity is a rather trivial theory, as the Einstein--Hilbert action
\begin{equation}
  S = \frac{1}{2\kappa} \int_{\Sigma} \dd^2x \, R \sqrt{-h}  
\end{equation}
is a topological invariant of the surface $\Sigma$. Consequently, the Einstein equations are automatically satisfied. 
In contrast, 2D \emph{quantum} gravity is a rather rich theory, with deep connections to the theory of integrable systems and algebraic geometry.
In the quantum setting, what is physically realised is not a fixed metric $h$ on the surface $\Sigma$, but rather a dynamically fluctuating metric. 
The quantity of interest, the path integral, is then an integral over the space of all such metrics up to symmetry:
\begin{equation}
  \left.
  \Set{ (\Sigma,h) | \substack{\displaystyle\text{surface $\Sigma$} \\[.5ex] \displaystyle\text{with metric $h$}}}
  \right/
  \substack{\text{diffeomorphism} \\ \text{conformal transf.}}
\end{equation}

In mathematical terms, we are interested in the space parametrising Riemann surfaces, and more precisely in the calculation of integrals over such a moduli space.

A completely different approach to 2D quantum gravity builds upon the idea of discretising the surfaces and counting triangulations, which in turn is related to random matrix theory.
The ``random matrix method'' started with G.~'t~Hooft's discovery in 1974 \cite{Hoo74} from the study of strong nuclear interactions, that matrix integrals are naturally related to graphs drawn on surfaces, weighted by their topology. 
This first example by 't~Hooft was then turned into a general paradigm for enumerating maps, by physicists E.~Brezin, C.~Itzykson, G.~Parisi, and J.-B.~Zuber \cite{BIPZ78}. 
By their method, they recovered some results due to the mathematician W.~T.~Tutte in the '60s, about counting the numbers of triangulations of the sphere \cite{Tut68}. 

In the continuum limit, one would expect the two approaches to coincide. 
The idea that these two models of 2D quantum gravity are equivalent has striking consequences and led E.~Witten to formulate his famous conjecture about the geometry of moduli spaces of Riemann surfaces \cite{Wit91}. 
\begin{wrapfigure}[13]{r}{0.42\textwidth}
  \begin{center}
    \begin{tikzpicture}[x=1pt,y=1pt,xscale=.6,yscale=.67]
        \draw[gray] (328.0749, 701.8622) .. controls (324, 724) and (334.9954, 746.3803) .. (350.1737, 746.3683);
        \draw[gray] (278.1297, 659.3086) .. controls (272, 663.4999) and (270, 677.7499) .. (271.6667, 688.875) .. controls (273.3333, 700) and (278.6667, 708) .. (283.0882, 712.181) .. controls (287.5097, 716.362) and (291.0194, 716.7239) .. (296.9778, 715.1659);
        \draw[gray] (276.6872, 721.7787) .. controls (276, 736) and (283.6175, 748.0911) .. (297.1115, 747.3679);
        \draw[gray] (260.4618, 670.1913) .. controls (250.3736, 675.5806) and (247.1868, 685.7903) .. (246.9267, 695.5618) .. controls (246.6667, 705.3333) and (249.3333, 714.6667) .. (252.849, 721.0041) .. controls (256.3646, 727.3415) and (260.7292, 730.6829) .. (266.9099, 726.6851);
        \draw[gray] (228.1827, 685.2222) .. controls (221.2802, 687.5763) and (217.9267, 695.1557) .. (217.7565, 704.7319) .. controls (217.5863, 714.3081) and (220.5994, 725.881) .. (227.3589, 735.1762) .. controls (234.1184, 744.4715) and (244.6243, 751.4891) .. (254.6582, 750.4118);

        \draw [thick] (294.6667, 694) .. controls (297.3333, 702.6667) and (302.6667, 709.3333) .. (309.3333, 711.3333) .. controls (316, 713.3333) and (324, 710.6667) .. (328.6667, 700.6667) .. controls (333.3333, 690.6667) and (334.6667, 673.3333) .. (333.3333, 661.3333) .. controls (332, 649.3333) and (328, 642.6667) .. (321.3333, 640) .. controls (314.6667, 637.3333) and (305.3333, 638.6667) .. (300, 643.3333) .. controls (294.6667, 648) and (293.3333, 656) .. (292.6667, 665.3333) .. controls (292, 674.6667) and (292, 685.3333) .. cycle;
        \draw [thick] (386, 742) .. controls (390.6667, 736) and (393.3333, 724) .. (394, 712.6667) .. controls (394.6667, 701.3333) and (393.3333, 690.6667) .. (390.6667, 683.3333) .. controls (388, 676) and (384, 672) .. (378.6667, 672.6667) .. controls (373.3333, 673.3333) and (366.6667, 678.6667) .. (362.6667, 688.6667) .. controls (358.6667, 698.6667) and (357.3333, 713.3333) .. (358.6667, 724) .. controls (360, 734.6667) and (364, 741.3333) .. (369.3333, 744.6667) .. controls (374.6667, 748) and (381.3333, 748) .. cycle;
        \draw [thick] (192.0661, 693.9418) .. controls (226.6887, 689.9806) and (263.2759, 672.6531) .. (301.8276, 641.9594);
        \draw [thick] (314.7517, 711.8974) .. controls (290.1371, 715.2197) and (271.1631, 722.2142) .. (257.8298, 732.8808);
        \draw [thick] (377.835, 672.8063) .. controls (357.4611, 678.0091) and (342.6289, 680.6999) .. (333.3384, 680.8788);
        \draw [thick] (377.8473, 746.841) .. controls (309.9491, 744.947) and (252.0648, 748.5166) .. (204.1943, 757.5498);
        \draw [thick] (192.3626, 693.8996) .. controls (177.6654, 696.4724) and (178.7924, 721.9009) .. (182.596, 737.5507) .. controls (186.3996, 753.2004) and (192.8797, 759.0713) .. (204.194, 757.55);
        \draw[dashed] [thick] (192.3626, 693.8996) .. controls (204.5305, 691.9556) and (212.5624, 698.815) .. (217.6451, 707.0424) .. controls (222.7278, 715.2697) and (224.8614, 724.8649) .. (223.9386, 734.7912) .. controls (223.0159, 744.7175) and (219.0367, 754.975) .. (204.194, 757.55);

        \draw (268, 732) -- (268, 612) -- (404, 664) -- (404, 760) -- cycle;
        \draw (132, 760) -- (132, 652) -- (264, 696) -- (264, 788) -- cycle;
        \draw[dotted,->] (132, 642) -- (268, 602);

        \node[rotate=90] at (122, 704) {\small space};
        \node[rotate=-16.75] at (185, 615) {\small time};
      \end{tikzpicture}
  \end{center}
  \caption{A string travelling through spacetime, tracing out a Riemann surface: the worldsheet of the string.}
\end{wrapfigure}
\leavevmode%
The conjecture, later proved by M.~Kontsevich \cite{Kon92}, connects in a beautiful way theoretical physics, algebraic geometry, and mathematical physics. 
Recently, the physics literature has seen a resurgence of such ideas in connection to Jackiw--Teitelboim gravity and its holographic dual, the Sachdev--Ye--Kitaev model \cite{Kit,SSS} (cf. C.~Johnson's and G.~Turiaci's lecture notes \cite{Joh,Tur26}).

Another physical theory, presenting deep connections with the theory of Riemann surfaces, is string theory. As a string travels through spacetime, it traces out a Riemann surface, the worldsheet of the string. These are essentially stringy versions of Feynman diagrams. The path integrals of the theory are mathematically described as integrals over the moduli spaces of Riemann surfaces mapping to the spacetime (cf. C.-C.~M.~Liu's lecture notes \cite{Liu}). The properties satisfied by such integrals are mathematically described by the notion of cohomological field theory.

The goal of these notes is to describe the mathematics related to such ideas, focusing particularly on the moduli space of Riemann surfaces, the concept of cohomological field theory, and its recursive solution.

The main references include:
\begin{itemize}[leftmargin=1.6cm,itemsep=.75em]
  \item[{\cite{Zvo12}}] D.~Zvonkine, ``An introduction to moduli spaces of curves and their intersection theory''
  \\
  {\footnotesize Not too technical notes on the moduli space of curves, its intersection theory, and Witten's conjecture}

  \item[{\cite{Pan19}}] R.~Pandharipande, ``Cohomological field theory calculations''
  \\
  {\footnotesize Not too technical notes on cohomological field theories, focused on examples}

  \item[{\cite{Sch20}}] J.~Schmitt. ``The moduli space of curves''
  \\
  {\footnotesize Algebro-geometric oriented notes on the moduli space of curves and its cohomology}

  \item[{\cite{ACG11}}]  E.~Arbarello , M.~Cornalba , P.~A. Griffiths, ``Geometry of Algebraic Curves, Vol. II''
  \\
  {\footnotesize A comprehensive text on Riemann surfaces and their moduli}
\end{itemize}

\addtocontents{toc}{\protect\setcounter{tocdepth}{0}}
\subsection*{Acknowledgements}
We would like to thank the Les Houches School for providing an extraordinary environment to work and learn. We are also grateful to the participants for their enthusiasm, insightful questions, valuable comments, and the enjoyable atmosphere.

We are deeply indebted to Gaëtan Borot, Bertrand Eynard, Reinier Kramer, Rahul Pandharipande, Adrien Sauvaget, Sergey Shadrin, and Dimitri Zvonkine for sharing their insights on moduli spaces over the years. Special thanks go to Sara Perletti for providing the base manuscript for these notes. We also thank Simon Sieroka for carefully reading an early draft and for his help in spotting several typos, and the referees for their constructive comments, which helped improve the final manuscript.

\subsection*{Funding information}
This work is partly a result of the ERC-SyG project, Recursive and Exact New Quantum Theory (ReNewQuantum) which received funding from the European Research Council (ERC) under the European Union's Horizon 2020 research and innovation programme under grant agreement No 810573. A.~G. was supported by an ETH Fellowship (22-2~FEL-003) and a Hermann-Weyl-Instructorship from the Forschungsinstitut für Mathematik at ETH Zürich. D.~L. is supported by the University of Trieste, by the INdAM group GNSAGA, and by the INFN within the project MMNLP (APINE).
\addtocontents{toc}{\protect\setcounter{tocdepth}{2}}

\newpage
\section{Moduli spaces of Riemann surfaces}
\label{sec:moduli}

In this section, we recall some facts about Riemann surfaces and their moduli space. The latter has been a central object in mathematics since Riemann's work in the mid-19th century, and a compactification was defined more than 50 years ago by Deligne and Mumford \cite{DM69} by including stable curves. For a great one-hour introductory talk to the moduli spaces of Riemann surfaces and their history, see \cite{Pan18}.

\subsection{Definition of the moduli spaces}
\begin{terminology}
  The primary focus of our study is on smooth, connected, compact, complex $1$-dimensional manifolds, simply called \emph{curves} or \emph{Riemann surfaces}, which have $n$ labelled distinct points (see M.~Bertola's lecture notes \cite{Ber}). These will be denoted as
  \begin{equation}
    (\Sigma,p_1,\dots,p_n) \,.
  \end{equation}
  Each compact complex curve has an underlying structure of a real $2$-dimensional orientable compact surface without boundary, uniquely characterised by its genus $g$.
\end{terminology}

\begin{figure}[b]
  \centering
  \begin{tikzpicture}[x=1pt,y=1pt,scale=.5]
    \fill[fill=gray!10] (47.999, 650.3347) arc[start angle=-106.6026, end angle=-15.2577, radius=56] .. controls (101.7321, 667.3896) and (78.3897, 654.4135) .. (47.999, 650.3347) -- cycle;
    \draw[thick] (64, 704) circle[radius=56];
    
    \begin{scope}[xshift=1cm]
      \fill[fill=gray!10] (190.857, 703.97) .. controls (196.8657, 707.8007) and (202.8753, 710.2593) .. (208.886, 711.346) -- (186.286, 708.867) -- cycle;
      \fill[fill=gray!10] (289.07, 687.496) .. controls (287.47, 683.952) and (285.335, 680.536) .. (282.667, 677.333) .. controls (269.333, 661.333) and (242.667, 650.667) .. (216, 650.667) .. controls (210.747, 650.667) and (205.494, 651.081) .. (200.342, 651.868) .. controls (238.572, 656.9007) and (268.148, 668.7767) .. (289.07, 687.496) -- cycle;
      \draw[thick] (149.3333, 730.6667) .. controls (162.6667, 746.6667) and (189.3333, 757.3333) .. (216, 757.3333) .. controls (242.6667, 757.3333) and (269.3333, 746.6667) .. (282.6667, 730.6667) .. controls (296, 714.6667) and (296, 693.3333) .. (282.6667, 677.3333) .. controls (269.3333, 661.3333) and (242.6667, 650.6667) .. (216, 650.6667) .. controls (189.3333, 650.6667) and (162.6667, 661.3333) .. (149.3333, 677.3333) .. controls (136, 693.3333) and (136, 714.6667) .. cycle;
      \draw[thick] (184, 712) .. controls (200, 688) and (232, 688) .. (248, 712);
      \draw[thick] (190.857, 703.97) .. controls (207.619, 714.6567) and (224.3913, 714.666) .. (241.174, 703.998);
    \end{scope}

    \begin{scope}[xshift=1.5cm]
      \node at (320, 704) {$\cdots$};
    \end{scope}
    
    \begin{scope}[xshift=2cm]
      \fill[fill=gray!10] (382.598, 704.208) .. controls (388.604, 707.942) and (394.602, 710.344) .. (400.592, 711.414) -- (378.2347, 708.93) -- cycle;
      \fill[fill=gray!10] (510.598, 704.208) .. controls (516.604, 707.942) and (522.602, 710.344) .. (528.592, 711.414) -- (506.2347, 708.93) -- cycle;
      \fill[fill=gray!10] (654.598, 704.208) .. controls (660.604, 707.942) and (666.602, 710.344) .. (672.592, 711.414) -- (650.2347, 708.93) -- cycle;
      \fill[fill=gray!10] (647.142, 656.215) .. controls (658.094, 653.478) and (669.047, 650.667) .. (680, 650.667) .. controls (701.333, 650.667) and (722.667, 661.333) .. (733.333, 677.333) .. controls (710.955, 661.549) and (682.2247, 654.5097) .. (647.142, 656.215) -- cycle;
      \draw[thick] (376, 712) .. controls (392, 688) and (424, 688) .. (440, 712);
      \draw[thick] (382.598, 704.208) .. controls (399.5327, 714.736) and (416.4027, 714.6763) .. (433.208, 704.029);
      \draw[thick] (504, 712) .. controls (520, 688) and (552, 688) .. (568, 712);
      \draw[thick] (510.598, 704.208) .. controls (527.5327, 714.736) and (544.4027, 714.6763) .. (561.208, 704.029);
      \node at (608, 704) {$\cdots$};
      \draw[thick] (648, 712) .. controls (664, 688) and (696, 688) .. (712, 712);
      \draw[thick] (654.598, 704.208) .. controls (671.5327, 714.736) and (688.4027, 714.6763) .. (705.208, 704.029);
      \draw[thick] (592.3778, 747.0845) .. controls (578.6667, 746.6667) and (557.3333, 757.3333) .. (536, 757.3333) .. controls (514.6667, 757.3333) and (493.3333, 746.6667) .. (472, 746.6667) .. controls (450.6667, 746.6667) and (429.3333, 757.3333) .. (408, 757.3333) .. controls (386.6667, 757.3333) and (365.3333, 746.6667) .. (354.6667, 730.6667) .. controls (344, 714.6667) and (344, 693.3333) .. (354.6667, 677.3333) .. controls (365.3333, 661.3333) and (386.6667, 650.6667) .. (408.0065, 650.6667) .. controls (429.3463, 650.6667) and (450.6927, 661.3333) .. (472.026, 661.3333) .. controls (493.3593, 661.3333) and (514.6797, 650.6667) .. (536.0065, 650.6667) .. controls (557.3333, 650.6667) and (578.6667, 661.3333) .. (592.0232, 660.8775);
      \draw[thick] (624.2546, 660.8468) .. controls (637.3333, 661.3333) and (658.6667, 650.6667) .. (680, 650.6667) .. controls (701.3333, 650.6667) and (722.6667, 661.3333) .. (733.3333, 677.3333) .. controls (744, 693.3333) and (744, 714.6667) .. (733.3333, 730.6667) .. controls (722.6667, 746.6667) and (701.3333, 757.3333) .. (680, 757.3333) .. controls (658.6667, 757.3333) and (637.3333, 746.6667) .. (624.0883, 747.1347);
    \end{scope}
  \end{tikzpicture}
\caption{Some examples of real $2$-dimensional orientable compact surfaces.}
\end{figure}

Our primary examples will be the sphere (genus $0$) and the torus (genus $1$). The sphere has a unique structure as a Riemann surface up to isomorphism, identified as the complex projective line $\P^1$. A complex curve of genus $0$ is called a \emph{rational curve}. The automorphism group of $\P^1$ is
\begin{equation}
  \mr{PSL}(2,\C) = \Set{
    \begin{pmatrix}
      a & b \\
      c & d
    \end{pmatrix}
    \ \bigg| \
    \substack{
      \displaystyle [a:b:c:d] \in \P^3 \\[.5ex]
      \displaystyle ad-bc \neq 0
    }
  }
\end{equation}
acting as
\begin{equation}
  \begin{pmatrix}
    a & b \\
    c & d
  \end{pmatrix}.z
  =
  \frac{az+b}{cz+d} \,.
\end{equation}
As for genus $1$, every Riemann surface structure on the torus is, up to isomorphism, obtained as a quotient $\C/\Lambda$. Here $\Lambda$ is a lattice, that is an additive group of the form
\begin{equation}
  \Lambda = \Set{ n_1\omega_1 + n_2\omega_2 | n_1,n_2 \in \Z}
\end{equation}
for $\omega_1,\omega_2 \in \C$ that are linearly independent over the reals. A complex curve of genus $1$ is referred to as an \emph{elliptic curve}.

As discussed in the introduction, we are interested in the moduli space of Riemann surfaces of a fixed genus $g$ with $n$ marked points (and in particular, we want to make sense of integrals over such a space: the path integrals of 2D quantum gravity).

\begin{definition}
  The \emph{moduli space} $\M_{g,n}$ is the set of isomorphism classes of Riemann surfaces of genus $g$ with $n$ marked points:
  \begin{equation}
    \M_{g,n}
    =
    \left.\Set{
      \substack{
        \displaystyle\text{Riemann surfaces} \\[.5ex]
        \displaystyle\text{of genus $g$ with $n$ marked points}
      }
    }
    \right/ \text{iso.}
  \end{equation}
  For isomorphism between two objects $(\Sigma,p_1,\dots,p_n)$ and $(\Sigma',p'_1,\dots,p'_n)$ we mean a biholomorphism $\phi \colon \Sigma \to \Sigma'$ that preserves the marked points: $\phi(p_i) = p_i'$.
\end{definition}

The above definition is perfectly well-posed, but we want to give it more structure. Recall that our goal is to discuss integrals over the moduli space of Riemann surfaces, so a structure like that of a manifold would be desirable. It turns out that there is a lot of geometry, but it is not as nice as that of a manifold. The main reason is that Riemann surfaces have automorphisms. The simplest example is $\P^1$, whose automorphism group is the infinite group $\mr{PSL}(2,\C)$. Since in the integration we want to quotient out by the group of symmetries, an infinite group of automorphisms is bad news. In other words, $\M_{0,0}$ does not have a nice geometric structure. There is however a way to get rid of automorphism by marking (at least three) points.

\begin{exercise}\leavevmode
  \begin{enumerate}
    \item\label{ex:PSL} Consider a genus $0$ curve with three marked points $(\P^1,p_1,p_2,p_3)$. Find the (unique) $g \in \mr{PSL}(2,\C)$ that maps $(\P^1,p_1,p_2,p_3)$ to $(\P^1,0,1,\infty)$.

    \item Consider a genus $0$ curve with four marked points $(\P^1,p_1,p_2,p_3,p_4)$. The group element $g \in \mr{PSL}(2,\C)$ found in part~(\labelcref{ex:PSL}) maps $(\P^1,p_1,p_2,p_3,p_4)$ to $(\P^1,0,1,\infty,t)$. Find an expression for $t$ as a function\footnote{
      This is known as the \emph{cross-ratio}, defined in deep antiquity (possibly already by Euclid) and considered by Pappus who noted its key invariance property.
    } of $p_1,p_2,p_3,p_4$.
  \end{enumerate}
\end{exercise}

The above exercise shows that
\begin{equation}
\begin{aligned}
  \M_{0,3}
  &=
  \set{(\P^1,0,1,\infty)} = \set{*} \,, \\
  \M_{0,4}
  &=
  \set{ (\P^1,0,1,\infty,t) | t \neq 0,1,\infty } = \P^1 \setminus \set{0,1,\infty} \,.
\end{aligned}
\end{equation}
One can generalise the above analysis to show that, for $n \ge 3$,
\begin{equation}
  \M_{0,n}
  =
  \Set{ (t_1,\dots,t_{n-3}) \in ( \P^1 \setminus \set{0,1,\infty} )^{n-3} | t_i \neq t_j } \,.
\end{equation}
This provides $\M_{0,n}$ with a nice geometric structure.

Another bad example where the automorphism group is infinite is that of an elliptic curve $E$, for which $\Aut(E)$ contains a subgroup isomorphic to $E$ itself acting by translations. Again, we can get rid of automorphisms (in this case, translations) by marking a point. If $E = \C/\Lambda$, a natural choice of marked point is the image of $\Lambda \subset \C$, that is the identity element on the torus. Thus, $\M_{1,1} = \set{\text{lattices}} / \C^*$, where $\C^*$ acts by rescaling. To understand the quotient, let us fix a basis $(\omega_1, \omega_2)$ of $\Lambda$. Multiplying $\Lambda$ by $1/\omega_1$, we obtain an equivalent lattice with basis $(1,\tau)$ for $\tau$ in the upper half-plane $\H$. Choosing another basis of the same lattice, that is acting by the group $\mr{SL}(2,\Z)$ of lattice base changes, we obtain another point $\tau' \in \H$. Thus, we find that
\begin{equation}
  \M_{1,1} = \H / \mr{SL}(2,\Z) \,.
\end{equation}
\begin{figure}[t]
  \centering
  \tikzset{every picture/.style=thick}
  \begin{tikzpicture} 
    \draw [->] (-3,0) -- (3,0);
    \draw [->] (0,0) -- (0,5);

    \node at (0,0) [below] {$0$};
    \node at (2,0) [below] {$1$};
    \node at (-2,0) [below] {$-1$};
    \node at (1,0) [below] {$\frac{1}{2}$};
    \node at (-1,0) [below] {$-\frac{1}{2}$};

    \node at (1,1.73205080757) [above right] {$B$};
    \node at (-1,1.73205080757) [above left] {$B'$};
    \node at (0,2) [below right] {$A$};
    \node at (1,4.5) [right] {$C$};
    \node at (-1,4.5) [left] {$C'$};

    \draw [dotted] (2,0) arc (0:60:2cm);
    \draw [thick] (1,1.73205080757) arc (60:120:2cm);
    \draw [dotted] (-1,1.73205080757) arc (120:180:2cm);

    \draw [dotted] (-1,0) -- (-1,1.73205080757);
    \draw [thick] (-1,1.73205080757) -- (-1,4.7);
    \draw [dotted] (1,0) -- (1,1.73205080757);
    \draw [thick] (1,1.73205080757) -- (1,4.7);

    \fill[opacity=.15,fill=gray] (1,1.73205080757) arc (60:120:2cm) -- (-1,4.7) -- (1,4.7) -- (1,1.73205080757);
  \end{tikzpicture}
  \caption{The moduli space $\M_{1,1}$. The arcs $AB$ and $AB'$ and the half-lines $BC$ and $B'C'$ are identified.}
  \label{fig:M11}
\end{figure}
A fundamental domain for the quotient is shown in \cref{fig:M11}. After identifying the arcs $AB \sim AB'$ and the half-lines $BC \sim B'C'$, we see that $\M_{1,1}$ is topologically $\P^1 \setminus \set{\infty}$. However, lattices have non-trivial automorphisms. Indeed, the matrix $-\Id$ acts trivially on $\H$, so that the automorphism group of each point on $\M_{1,1}$ contains at least $\Z_{2}$ as a subgroup. This involution is called the hyperelliptic involution of a marked elliptic curve. If we write an elliptic curve as (the compactification of) a degree $3$ polynomial equation of the form
\begin{equation}
  E \colon
  \quad
  y^2 = x^3 + ax + b \,,
\end{equation}
then the hyperelliptic involution is simply the map $y \mapsto -y$.

It is actually possible to completely characterise the automorphism group of each point $\tau$ in the fundamental domain (see \cref{fig:aut:lattices}):
\begin{itemize}
  \item for $\tau =  e^{\pi\iu/3} = \frac{1 + \iu\sqrt{3}}{2}$ corresponding to the hexagonal lattice, the automorphism group is $\Z_6$;

  \item for $\tau = e^{\pi\iu/2} = \iu$ corresponding to the square lattice, the automorphism group is $\Z_4$;

  \item for any other $\tau$ in the fundamental domain, the automorphism group is $\Z_2$.
\end{itemize}
\begin{figure}
  \centering
  \begin{tikzpicture}
    \definecolor{blue1}{RGB}{3, 4, 94}
    \definecolor{blue2}{RGB}{2, 62, 138}
    \definecolor{blue3}{RGB}{0, 119, 182}
    \definecolor{blue4}{RGB}{0, 150, 199}
    \definecolor{blue5}{RGB}{0, 180, 216}
    \definecolor{blue6}{RGB}{72, 202, 228}

    \fill [blue1,opacity=.4] (0,0) -- (0:1) -- ($(0:1) + (60:1)$) -- (60:1) -- cycle;
    \fill [blue2,opacity=.4] (0,0) -- (60:1) -- ($(60:1) + (120:1)$) -- (120:1) -- cycle;
    \fill [blue3,opacity=.4] (0,0) -- (120:1) -- ($(120:1) + (180:1)$) -- (180:1) -- cycle;
    \fill [blue4,opacity=.4] (0,0) -- (180:1) -- ($(180:1) + (240:1)$) -- (240:1) -- cycle;
    \fill [blue5,opacity=.4] (0,0) -- (240:1) -- ($(240:1) + (300:1)$) -- (300:1) -- cycle;
    \fill [blue6,opacity=.4] (0,0) -- (300:1) -- ($(300:1) + (0:1)$) -- (0:1) -- cycle;

    \node at (0,0) {\tiny$\bullet$};

    \node at (0:1) {\tiny$\bullet$};
    \node at (60:1) {\tiny$\bullet$};
    \node at (120:1) {\tiny$\bullet$};
    \node at (180:1) {\tiny$\bullet$};
    \node at (240:1) {\tiny$\bullet$};
    \node at (300:1) {\tiny$\bullet$};

    \node at ($(0:1) + (60:1)$) {\tiny$\bullet$};
    \node at ($(60:1) + (120:1)$) {\tiny$\bullet$};
    \node at ($(120:1) + (180:1)$) {\tiny$\bullet$};
    \node at ($(180:1) + (240:1)$) {\tiny$\bullet$};
    \node at ($(240:1) + (300:1)$) {\tiny$\bullet$};
    \node at ($(300:1) + (0:1)$) {\tiny$\bullet$};

    \node at (0,-2.25) {$\tau = e^{\pi\iu/3}$};

    \begin{scope}[xshift=5cm]
      \definecolor{green1}{RGB}{7, 59, 58}
      \definecolor{green2}{RGB}{11, 110, 79}
      \definecolor{green3}{RGB}{8, 160, 69}
      \definecolor{green4}{RGB}{107, 191, 89}

      \fill [green1,opacity=.4] (0,0) -- (0:1) -- ($(0:1) + (90:1)$) -- (90:1) -- cycle;
      \fill [green2,opacity=.4] (0,0) -- (90:1) -- ($(90:1) + (180:1)$) -- (180:1) -- cycle;
      \fill [green3,opacity=.4] (0,0) -- (180:1) -- ($(180:1) + (270:1)$) -- (270:1) -- cycle;
      \fill [green4,opacity=.4] (0,0) -- (270:1) -- ($(270:1) + (0:1)$) -- (0:1) -- cycle;

      \node at (0,0) {\tiny$\bullet$};

      \node at (0:1) {\tiny$\bullet$};
      \node at (90:1) {\tiny$\bullet$};
      \node at (180:1) {\tiny$\bullet$};
      \node at (270:1) {\tiny$\bullet$};

      \node at ($(0:1) + (90:1)$) {\tiny$\bullet$};
      \node at ($(90:1) + (180:1)$) {\tiny$\bullet$};
      \node at ($(180:1) + (270:1)$) {\tiny$\bullet$};
      \node at ($(270:1) + (0:1)$) {\tiny$\bullet$};

      \node at (0,-2.25) {$\tau = e^{\pi\iu/2}$};
    \end{scope}

    \begin{scope}[xshift=10cm]
      \definecolor{red1}{RGB}{255, 190, 11}
      \definecolor{red2}{RGB}{251, 86, 7}

      \fill [red1,opacity=.4] (0,0) -- (0:1) -- ($(0:1) + (75:1)$) -- (75:1) -- cycle;
      \fill [red2,opacity=.4] (0,0) -- (180:1) -- ($(180:1) + (255:1)$) -- (255:1) -- cycle;

      \node at (0,0) {\tiny$\bullet$};

      \node at (0:1) {\tiny$\bullet$};
      \node at (75:1) {\tiny$\bullet$};
      \node at (180:1) {\tiny$\bullet$};
      \node at (255:1) {\tiny$\bullet$};

      \node at ($(0:1) + (75:1)$) {\tiny$\bullet$};
      \node at ($(180:1) + (255:1)$) {\tiny$\bullet$};

      \node at (0,-2.25) {$\tau = e^{5\pi\iu/24}$};
    \end{scope}
  \end{tikzpicture}
  \caption{The automorphism groups of lattices.}
  \label{fig:aut:lattices}
\end{figure}

A theorem by A.~Hurwitz implies that the automorphism group of any Riemann surface satisfying $2g - 2 + n > 0$ is finite. Such a pair $(g,n)$ is called \emph{stable}. Conversely, every Riemann surface with $2g - 2 + n \le 0$ has an infinite group of automorphisms that preserve the marked points. In other words:
\begin{equation}
  \Aut(\Sigma_g,p_1,\dots,p_n) \text{ is finite}
  \qquad\Longleftrightarrow\qquad
  -\chi = 2g - 2 + n > 0 \,.
\end{equation}
This precludes defining the moduli spaces $\M_{0,0}$, $\M_{0,1}$, $\M_{0,2}$, and $\M_{1,0}$ as nice geometric spaces. (While they can still be considered as sets, this is of limited use.)

From now on, we will always assume $2g - 2 + n > 0$. In this case the situation is good, but not as good as it can get: there are still curves with non-trivial automorphism group, as the example of $\M_{1,1}$ showed. Nonetheless, finiteness of the automorphism groups allows us to consider the moduli space of Riemann surfaces as an orbifold.

\begin{theorem}
  For $2g - 2 + n > 0$, the moduli space $\M_{g,n}$ is a connected, smooth, complex orbifold of dimension
  \begin{equation}
    \dim(\M_{g,n}) = 3g - 3 + n \,.
  \end{equation}
\end{theorem}

The definition of a smooth complex orbifold is somewhat technical, but it closely mirrors that of a smooth complex manifold. The key difference is that, locally, an orbifold looks like an open subset of $\C^d/G$, where $G$ is a \emph{finite} group acting holomorphically. A simple example to keep in mind is the global quotient $\C/\Z_{m}$, where $\Z_m$ acts by rotation through an angle $2\pi/m$. In the theorem above, ``smooth'' is meant in this orbifold sense: the moduli space $\M_{g,n}$ is locally modelled on quotients of smooth complex manifolds by finite group actions. In practical terms, this means that $\M_{g,n}$ is a smooth manifold up to finite quotient singularities.

\begin{exercise}
  For the reader familiar with Riemann--Roch and Riemann--Hurwitz, convince yourself that the complex dimension of $\M_{g} = \M_{g,0}$ is $3g-3$. This result was already known to Riemann himself, who also coined the term ``moduli space'' (from the Latin word \emph{modus}, meaning measure):
  \begin{center}
    \includegraphics[width=.75\textwidth]{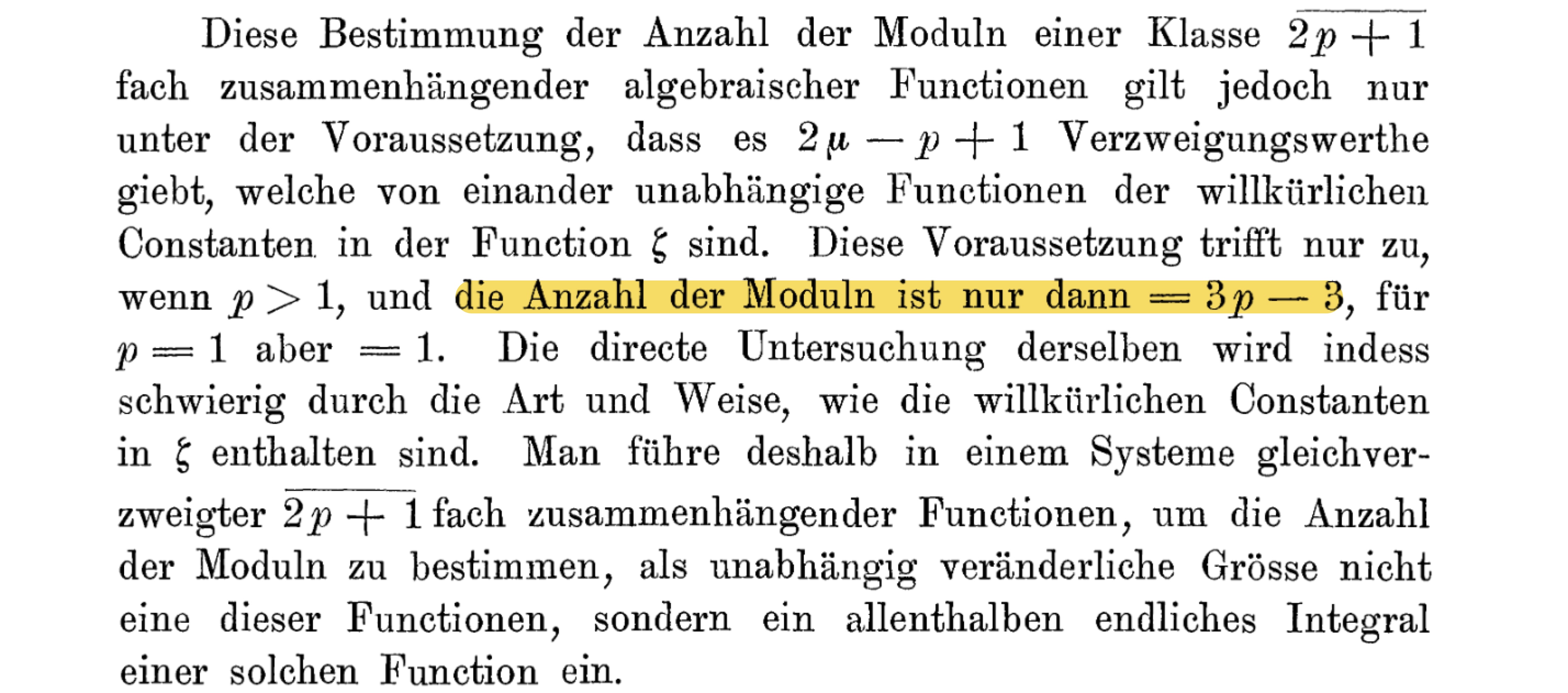}
  \end{center}
  To this end, consider the moduli space of pairs $(\Sigma,f)$, where $\Sigma$ is a genus $g$ Riemann surface and $f$ is a degree $d$ holomorphic map from $\Sigma$ to $\P^1$ (i.e. a meromorphic function on $\Sigma$). Such a space is sometimes called a Hurwitz space, denoted $\mathcal{H}_{g,d}$. Compute its dimension in two different ways.
  \begin{itemize}
    \item The dimension of $\mathcal{H}_{g,d}$ equals the dimension of $\mathcal{M}_g$, counting the ``number of deformation parameters'' of the Riemann surface $\Sigma$, plus the ``number of deformation parameters'' of the function $f$. Compute the latter via Riemann--Roch.

    \item Directly compute the dimension of $\mathcal{H}_{g,d}$ using Riemann--Hurwitz.
  \end{itemize}
  Conclude that $\dim{\mathcal{M}_g} = 3g-3$.

  {\small \emph{\faLightbulb \ Hint.} Consider $d \gg 1$.
  }
\end{exercise}

One important consequence of this local description is that many differential-geometric notions familiar from manifolds, such as differential forms and integration, extend naturally to orbifolds. In particular, it makes sense to talk about integration over complex orbifolds. For the example of $\C/\Z_{m}$, given a function $f \colon \C \to \C$ that is invariant under rotation of $\frac{2\pi}{m}$, we can define
\begin{equation}
  \int_{\C/\Z_{m}} f(z,\bar{z}) \, \dd z \, \dd\bar{z}
  =
  \frac{1}{|\Z_m|} \int_{\C} f(z,\bar{z}) \, \dd z \, \dd\bar{z} \,.
\end{equation}
Most of the results that hold for manifolds extend (with proper modifications) to orbifolds. Here is an example of the Euler characteristic.

\begin{exercise}
  The Euler characteristic of an orbifold $X$ is defined as
  \begin{equation}
    \chi(X) = \sum_{G} \frac{\chi(X_G)}{|G|} \,,
  \end{equation}
  where $X_G$ is the locus of points with automorphism group $G$. Prove that $\chi(\M_{1,1}) = - \frac{1}{12}$. The formula generalises to the celebrated Harer--Zagier formula \cite{HZ86}:
  \begin{equation}
    \chi(\M_{g,n}) = (1-2g)_{n-1} \, \zeta(1-2g) \,,
  \end{equation}
  where $(x)_m$ denotes the Pochhammer symbol (or falling factorial) and $\zeta(x)$ is the Riemann zeta function. 
  Interestingly, the original computation by Harer and Zagier uses matrix model techniques.
\end{exercise}

Although integrals over orbifolds are well-defined, there is another potential issue to deal with: non-compactness. 
The non-compactness problem can be seen already from the examples of $\M_{0,4}$ or $\M_{1,1}$. The latter is topologically $\P^1 \setminus \set{\infty}$, with the missing point at infinity being the source of non-compactness. We can see how this limit point is realised geometrically by considering the family of elliptic curves
\begin{equation}
  E_t \colon
  \quad
  y^2 = x(x-1)(x-t) \,,
  \qquad\qquad
  t \in (0,1) \,.
\end{equation}
In the limit $t \to 0$ or $1$, the Riemann surface $E_t$ becomes degenerate. For instance, as $t \to 0$ we find $y^2 = x^2(x-1)$, which, locally around $x = 0$, looks like the union of the two complex lines $y = \pm x$. This means that at $x = 0$ we have two meeting components, also known as a \emph{nodal singularity}, and the surface $E_0$ will look as in \cref{fig:pinched:torus}. In other words, the limit point of $\M_{1,1}$ is not a torus anymore, but rather a pinched torus.

\begin{figure}
  \tikzset{every picture/.style=thick}
  \centering
  \begin{tikzpicture} [x=1pt,y=1pt,scale=.425]
    \draw(240, 736) .. controls (240, 704) and (264, 704) .. (280, 700) .. controls (296, 696) and (304, 688) .. (304, 672);
    \draw(304, 672) .. controls (304, 656) and (280, 648) .. (257.3333, 649.3333) .. controls (234.6667, 650.6667) and (213.3333, 661.3333) .. (200.908, 673.0217) .. controls (188.4827, 684.71) and (184.9653, 697.42) .. (184.9653, 710.7533) .. controls (184.9653, 724.0867) and (188.4827, 738.0433) .. (195.5747, 750.355) .. controls (202.6667, 762.6667) and (213.3333, 773.3333) .. (226.6667, 781.3333) .. controls (240, 789.3333) and (256, 794.6667) .. (275.3333, 796.6667) .. controls (294.6667, 798.6667) and (317.3333, 797.3333) .. (336.6667, 793.3333) .. controls (356, 789.3333) and (372, 782.6667) .. (382, 770) .. controls (392, 757.3333) and (396, 738.6667) .. (396, 724) .. controls (396, 709.3333) and (392, 698.6667) .. (386.6667, 688.6667) .. controls (381.3333, 678.6667) and (374.6667, 669.3333) .. (356.3333, 659.6667) .. controls (338, 650) and (308, 640) .. (304, 672);
    \draw(304, 672) .. controls (304, 688) and (312, 696) .. (319.3333, 701.3333) .. controls (326.6667, 706.6667) and (333.3333, 709.3333) .. (338.6667, 715.6667) .. controls (344, 722) and (348, 732) .. (344, 744);
    \draw(244, 717.438) .. controls (252, 740) and (272, 746) .. (295, 746) .. controls (318, 746) and (344, 740) .. (343.974, 724.848);
    \node at (304, 780) {\small$\bullet$};
  \end{tikzpicture}
\caption{A pinched torus.}
\label{fig:pinched:torus}
\end{figure}

To make sense of integration over non-compact spaces we have two possibilities. The first one is to consider only functions or differential forms with a certain decay at limit points. The second option is to properly compactify the space of interest, and only consider regular functions or differential forms on such compactification. We will follow the second route. It turns out that for $\M_{g,n}$ the addition of Riemann surfaces with nodes is sufficient to obtain a nice compactification.

\begin{definition}
  A \emph{stable Riemann surface} of genus $g$ with $n$ labelled marked points $p_1,\dots,p_n$ is a possibly singular, compact, connected, complex curve $\Sigma$ such that:
  \begin{itemize}
    \item the genus of the surface obtained from $\Sigma$ by smoothening all its nodes is $g$ (see \cref{fig:sm:norm}),

    \item the only singularities of $\Sigma$ are nodes,

    \item the marked points are distinct and do not coincide with the nodes, and
    
    \item $(\Sigma,p_1,\dots,p_n)$ has a finite number of automorphisms.
  \end{itemize}
  We can then define a moduli space parametrising isomorphism classes of \emph{stable} Riemann surfaces, often called the Deligne--Mumford moduli space \cite{DM69}:
  \begin{equation}
    \Mbar_{g,n}
    =
    \left.\Set{
      \substack{
        \displaystyle\text{stable Riemann surfaces} \\[.5ex]
        \displaystyle\text{of genus $g$ with $n$ marked points}
      }
    }
    \right/ \text{iso.}
  \end{equation}
\end{definition}

\begin{figure}[b]
  \centering
  \begin{tikzpicture}[x=1pt,y=1pt,scale=.325]
    \fill[fill=gray!10] (166.598, 512.208) .. controls (172.604, 515.942) and (178.602, 518.344) .. (184.592, 519.414) -- (162.2347, 516.93) -- cycle;
    \fill[fill=gray!10] (294.598, 512.208) .. controls (300.604, 515.942) and (306.602, 518.344) .. (312.592, 519.414) -- (290.2347, 516.93) -- cycle;
    \fill[fill=gray!10] (422.598, 512.208) .. controls (428.604, 515.942) and (434.602, 518.344) .. (440.592, 519.414) -- (418.2347, 516.93) -- cycle;
    \fill[fill=gray!10] (415.142, 464.215) .. controls (426.094, 461.478) and (437.047, 458.667) .. (448, 458.667) .. controls (469.333, 458.667) and (490.667, 469.333) .. (501.333, 485.333) .. controls (478.955, 469.549) and (450.2247, 462.5097) .. (415.142, 464.215) -- cycle;
    \fill[fill=gray!10] (402.7979, 618.2304) .. controls (408.8066, 622.061) and (414.8162, 624.5197) .. (420.8269, 625.6064) -- (398.2269, 623.1274) -- cycle;

    \draw[thick] (288, 520) .. controls (304, 496) and (336, 496) .. (352, 520);
    \draw[thick] (294.598, 512.208) .. controls (311.5327, 522.736) and (328.4027, 522.6763) .. (345.208, 512.029);
    \draw[thick] (416, 520) .. controls (432, 496) and (464, 496) .. (480, 520);
    \draw[thick] (422.598, 512.208) .. controls (439.5327, 522.736) and (456.4027, 522.6763) .. (473.208, 512.029);
    \draw[thick] (395.9409, 626.2604) .. controls (411.9409, 602.2604) and (443.9409, 602.2604) .. (459.9409, 626.2604);
    \draw[thick] (402.7979, 618.2304) .. controls (419.5599, 628.917) and (436.3322, 628.9264) .. (453.1149, 618.2584);
    \draw[thick] (320, 565.3333) .. controls (298.6667, 565.3333) and (277.3333, 554.6667) .. (256, 554.6667) .. controls (234.6667, 554.6667) and (213.3333, 565.3333) .. (192, 565.3333) .. controls (170.6667, 565.3333) and (149.3333, 554.6667) .. (138.6667, 538.6667) .. controls (128, 522.6667) and (128, 501.3333) .. (138.6667, 485.3333) .. controls (149.3333, 469.3333) and (170.6667, 458.6667) .. (192.0065, 458.6667) .. controls (213.3463, 458.6667) and (234.6927, 469.3333) .. (256.026, 469.3333) .. controls (277.3593, 469.3333) and (298.6797, 458.6667) .. (320.0065, 458.6667) .. controls (341.3333, 458.6667) and (362.6667, 469.3333) .. (384, 469.3333) .. controls (405.3333, 469.3333) and (426.6667, 458.6667) .. (448, 458.6667) .. controls (469.3333, 458.6667) and (490.6667, 469.3333) .. (501.3333, 485.3333) .. controls (512, 501.3333) and (512, 522.6667) .. (501.3333, 538.6667) .. controls (490.6667, 554.6667) and (469.3333, 565.3333) .. (448, 565.3333) .. controls (426.6667, 565.3333) and (405.3333, 554.6667) .. (384, 554.6667) .. controls (362.6667, 554.6667) and (341.3333, 565.3333) .. cycle;
    \draw[thick] (160, 520) .. controls (176, 496) and (208, 496) .. (224, 520);
    \draw[thick] (166.598, 512.208) .. controls (183.5327, 522.736) and (200.4027, 522.6763) .. (217.208, 512.029);
    \draw[thick] (220.2035, 571.2317) .. controls (233.7098, 581.0109) and (249.7098, 599.6775) .. (274.2778, 608.5074) .. controls (298.8458, 617.3373) and (331.9818, 616.3305) .. (339.9818, 625.6638) .. controls (347.9818, 634.9971) and (330.8458, 654.6707) .. (302.2778, 664.5074) .. controls (273.7098, 674.3442) and (233.7098, 674.3442) .. (212.2035, 655.2317) .. controls (190.6973, 636.1192) and (187.6847, 597.8942) .. (191.6847, 579.2275) .. controls (195.6847, 560.5608) and (206.6973, 561.4525) .. cycle;
    \draw[thick] (373.0664, 587.4258) .. controls (401.428, 573.1576) and (449.428, 559.8243) .. (476.0947, 567.8243) .. controls (502.7614, 575.8243) and (508.0947, 605.1576) .. (504.0947, 622.491) .. controls (500.0947, 639.8243) and (486.7614, 645.1576) .. (466.7614, 651.8243) .. controls (446.7614, 658.491) and (420.0947, 666.491) .. (395.7331, 663.4258) .. controls (371.3715, 660.3607) and (349.3149, 646.2304) .. (342.6482, 631.5637) .. controls (335.9816, 616.897) and (344.7048, 601.694) .. cycle;
    \node at (224, 656) {\tiny$\bullet$};
    \node at (480, 632) {\tiny$\bullet$};
    \node at (176, 472) {\tiny$\bullet$};
    \node at (144, 496) {\tiny$\bullet$};

    \begin{scope}[xshift=-15cm]
      \fill[fill=gray!10] (402.7979, 618.2304) .. controls (408.8066, 622.0611) and (414.8162, 624.5197) .. (420.8269, 625.6064) -- (398.2269, 623.1274) -- cycle;
      \fill[fill=gray!10] (166.598, 512.208) .. controls (172.604, 515.942) and (178.602, 518.344) .. (184.592, 519.414) -- (162.2347, 516.93) -- cycle;
      \fill[fill=gray!10] (294.598, 512.208) .. controls (300.604, 515.942) and (306.602, 518.344) .. (312.592, 519.414) -- (290.2347, 516.93) -- cycle;
      \fill[fill=gray!10] (422.598, 512.208) .. controls (428.604, 515.942) and (434.602, 518.344) .. (440.592, 519.414) -- (418.2347, 516.93) -- cycle;
      \fill[fill=gray!10] (415.142, 464.215) .. controls (426.094, 461.478) and (437.047, 458.667) .. (448, 458.667) .. controls (469.333, 458.667) and (490.667, 469.333) .. (501.333, 485.333) .. controls (478.955, 469.549) and (450.2247, 462.5097) .. (415.142, 464.215) -- cycle;
      \draw[thick] (395.9409, 626.2604) .. controls (411.9409, 602.2604) and (443.9409, 602.2604) .. (459.9409, 626.2604);
      \draw[thick] (402.7979, 618.2304) .. controls (419.5599, 628.9171) and (436.3322, 628.9264) .. (453.1149, 618.2584);
      \draw[thick] (160, 520) .. controls (176, 496) and (208, 496) .. (224, 520);
      \draw[thick] (166.598, 512.208) .. controls (183.5327, 522.736) and (200.4027, 522.6763) .. (217.208, 512.029);
      \draw[thick] (288, 520) .. controls (304, 496) and (336, 496) .. (352, 520);
      \draw[thick] (294.598, 512.208) .. controls (311.5327, 522.736) and (328.4027, 522.6763) .. (345.208, 512.029);
      \draw[thick] (416, 520) .. controls (432, 496) and (464, 496) .. (480, 520);
      \draw[thick] (422.598, 512.208) .. controls (439.5327, 522.736) and (456.4027, 522.6763) .. (473.208, 512.029);
      \draw[thick] (373.0663, 587.4258) .. controls (394.2039, 576.792) and (408.6244, 572.2404) .. (421.0735, 570.0424) .. controls (433.5225, 567.8444) and (444, 568) .. (444, 564)(456, 564) .. controls (456, 568) and (469.3023, 565.7866) .. (476.0947, 567.8243) .. controls (502.7613, 575.8243) and (508.0947, 605.1576) .. (504.0947, 622.491) .. controls (500.0947, 639.8243) and (486.7613, 645.1576) .. (466.7613, 651.8243) .. controls (446.7613, 658.491) and (420.0947, 666.491) .. (395.733, 663.4258) .. controls (371.3713, 660.3606) and (361.4458, 654.002) .. (355.603, 646.9118) .. controls (349.7602, 639.8217) and (348, 632) .. (340, 632)(340, 624) .. controls (348, 624) and (344, 604) .. (373.0663, 587.4258)(340, 632) .. controls (336, 632) and (330.8458, 654.6705) .. (302.2778, 664.5074) .. controls (273.7098, 674.3442) and (233.7098, 674.3442) .. (212.2034, 655.2317) .. controls (190.6971, 636.1192) and (189.3414, 618.918) .. (189.6671, 602.8882) .. controls (189.9929, 586.8584) and (192, 572) .. (192, 564) .. controls (192, 556) and (180, 558) .. (168.6106, 556.5825) .. controls (157.2212, 555.1651) and (146.4423, 550.3301) .. (138.6667, 538.6667) .. controls (128, 522.6667) and (128, 501.3333) .. (138.6667, 485.3333) .. controls (149.3333, 469.3333) and (170.6667, 458.6667) .. (192.0065, 458.6667) .. controls (213.3463, 458.6667) and (234.6927, 469.3333) .. (256.026, 469.3333) .. controls (277.3593, 469.3333) and (298.6797, 458.6667) .. (320.0065, 458.6667) .. controls (341.3333, 458.6667) and (362.6667, 469.3333) .. (384, 469.3333) .. controls (405.3333, 469.3333) and (426.6667, 458.6667) .. (448, 458.6667) .. controls (469.3333, 458.6667) and (490.6667, 469.3333) .. (501.3333, 485.3333) .. controls (512, 501.3333) and (512, 522.6667) .. (501.3333, 538.6667) .. controls (490.6667, 554.6667) and (483.3333, 557.3333) .. (474.6667, 558.6667) .. controls (466, 560) and (456, 560) .. (456, 564)(444, 564) .. controls (444, 560) and (432, 560) .. (422.3333, 558.6667) .. controls (412.6667, 557.3333) and (405.3333, 554.6667) .. (384, 554.6667) .. controls (362.6667, 554.6667) and (341.3333, 565.3333) .. (320, 565.3333) -- (320, 565.3333) .. controls (298.6667, 565.3333) and (277.3333, 554.6667) .. (256, 554.6667) .. controls (245.3333, 554.6667) and (208, 556) .. (208, 564) .. controls (208, 572) and (221.6549, 577.4388) .. (232.0823, 584.3582) .. controls (242.5098, 591.2775) and (249.7098, 599.6775) .. (274.2778, 608.5074) .. controls (298.8458, 617.3372) and (332, 624) .. (340, 624);
      \node at (224, 656) {\tiny$\bullet$};
      \node at (480, 632) {\tiny$\bullet$};
      \node at (176, 472) {\tiny$\bullet$};
      \node at (144, 496) {\tiny$\bullet$};

      \node at (330,425) {smoothening};
    \end{scope}

    \begin{scope}[xshift=15cm]
     \fill[fill=gray!10] (166.598, 512.208) .. controls (172.604, 515.942) and (178.602, 518.344) .. (184.592, 519.414) -- (162.2347, 516.93) -- cycle;
     \fill[fill=gray!10] (422.598, 512.208) .. controls (428.604, 515.942) and (434.602, 518.344) .. (440.592, 519.414) -- (418.2347, 516.93) -- cycle;
     \fill[fill=gray!10] (294.598, 512.208) .. controls (300.604, 515.942) and (306.602, 518.344) .. (312.592, 519.414) -- (290.2347, 516.93) -- cycle;
     \fill[fill=gray!10] (415.142, 464.215) .. controls (426.094, 461.478) and (437.047, 458.667) .. (448, 458.667) .. controls (469.333, 458.667) and (490.667, 469.333) .. (501.333, 485.333) .. controls (478.955, 469.549) and (450.2247, 462.5097) .. (415.142, 464.215) -- cycle;
     \fill[fill=gray!10] (418.7979, 634.2304) .. controls (424.8066, 638.0611) and (430.8162, 640.5197) .. (436.8269, 641.6064) -- (414.2269, 639.1274) -- cycle;

      \draw[thick] (411.9409, 642.2604) .. controls (427.9409, 618.2604) and (459.9409, 618.2604) .. (475.9409, 642.2604);
      \draw[thick] (418.7979, 634.2304) .. controls (435.5599, 644.9171) and (452.3322, 644.9264) .. (469.1149, 634.2584);
      \draw[thick] (320, 565.3333) .. controls (298.6667, 565.3333) and (277.3333, 554.6667) .. (256, 554.6667) .. controls (234.6667, 554.6667) and (213.3333, 565.3333) .. (192, 565.3333) .. controls (170.6667, 565.3333) and (149.3333, 554.6667) .. (138.6667, 538.6667) .. controls (128, 522.6667) and (128, 501.3333) .. (138.6667, 485.3333) .. controls (149.3333, 469.3333) and (170.6667, 458.6667) .. (192.0065, 458.6667) .. controls (213.3463, 458.6667) and (234.6927, 469.3333) .. (256.026, 469.3333) .. controls (277.3593, 469.3333) and (298.6797, 458.6667) .. (320.0065, 458.6667) .. controls (341.3333, 458.6667) and (362.6667, 469.3333) .. (384, 469.3333) .. controls (405.3333, 469.3333) and (426.6667, 458.6667) .. (448, 458.6667) .. controls (469.3333, 458.6667) and (490.6667, 469.3333) .. (501.3333, 485.3333) .. controls (512, 501.3333) and (512, 522.6667) .. (501.3333, 538.6667) .. controls (490.6667, 554.6667) and (469.3333, 565.3333) .. (448, 565.3333) .. controls (426.6667, 565.3333) and (405.3333, 554.6667) .. (384, 554.6667) .. controls (362.6667, 554.6667) and (341.3333, 565.3333) .. cycle;
      \draw[thick] (160, 520) .. controls (176, 496) and (208, 496) .. (224, 520);
      \draw[thick] (166.598, 512.208) .. controls (183.5327, 522.736) and (200.4027, 522.6763) .. (217.208, 512.029);
      \draw[thick] (288, 520) .. controls (304, 496) and (336, 496) .. (352, 520);
      \draw[thick] (294.598, 512.208) .. controls (311.5327, 522.736) and (328.4027, 522.6763) .. (345.208, 512.029);
      \draw[thick] (416, 520) .. controls (432, 496) and (464, 496) .. (480, 520);
      \draw[thick] (422.598, 512.208) .. controls (439.5327, 522.736) and (456.4027, 522.6763) .. (473.208, 512.029);
      \draw[thick] (204.2034, 587.2317) .. controls (217.7098, 597.0109) and (233.7098, 615.6775) .. (258.2778, 624.5074) .. controls (282.8458, 633.3372) and (315.9818, 632.3302) .. (323.9818, 641.6635) .. controls (331.9818, 650.9969) and (314.8458, 670.6705) .. (286.2778, 680.5074) .. controls (257.7098, 690.3442) and (217.7098, 690.3442) .. (196.2034, 671.2317) .. controls (174.6971, 652.1192) and (171.6844, 613.8942) .. (175.6844, 595.2275) .. controls (179.6844, 576.5609) and (190.6971, 577.4525) .. cycle;
      \draw[thick] (389.0663, 603.4258) .. controls (417.428, 589.1576) and (465.428, 575.8243) .. (492.0947, 583.8243) .. controls (518.7613, 591.8243) and (524.0947, 621.1576) .. (520.0947, 638.491) .. controls (516.0947, 655.8243) and (502.7613, 661.1576) .. (482.7613, 667.8243) .. controls (462.7613, 674.491) and (436.0947, 682.491) .. (411.733, 679.4258) .. controls (387.3713, 676.3606) and (365.3147, 662.2303) .. (358.648, 647.5636) .. controls (351.9813, 632.897) and (360.7047, 617.694) .. cycle;
      \node at (208, 672) {\tiny$\bullet$};
      \node at (496, 648) {\tiny$\bullet$};
      \node at (176, 472) {\tiny$\bullet$};
      \node at (144, 496) {\tiny$\bullet$};
      \node at (186.3791, 580.5624) {\tiny$\bullet$};
      \node at (186.3791, 580.5624) {\tiny$\bullet$};
      \node at (325.5872, 644.5672) {\tiny$\bullet$};
      \node at (325.5872, 644.5672) {\tiny$\bullet$};
      \node at (357.587, 644.567) {\tiny$\bullet$};
      \node at (202.379, 564.562) {\tiny$\bullet$};
      \node at (202.379, 564.562) {\tiny$\bullet$};
      \node at (467.0156, 581.5871) {\tiny$\bullet$};
      \node at (451.016, 565.587) {\tiny$\bullet$};

      \node at (330,425) {normalisation};
    \end{scope}
  \end{tikzpicture}
  \caption{The smoothening and the normalisation of a singular Riemann surface. From the smoothening, one reads $(g,n)$; from the normalisation, one reads the stability condition.}
  \label{fig:sm:norm}
\end{figure}

The last condition in the above definition can be reformulated as follows. Let $\Sigma_1, \dots , \Sigma_k$ be the connected components of the surface obtained by separating the branches at each node (this process is called normalisation, see \cref{fig:sm:norm}). Let $g(v)$ be the genus of $\Sigma_v$ and $n(v)$ the number of special points, i.e., marked points and preimages of the nodes on $\Sigma_i$. Then, the ``finite automorphisms'' condition is satisfied if and only if $2g(v) - 2 + n(v) > 0$ for all $v$.

The main result about the Deligne--Mumford moduli space is that it provides a compactification of the moduli space of Riemann surfaces.

\begin{theorem}
  For $2g - 2 + n > 0$, the moduli space $\Mbar_{g,n}$
  \begin{itemize}
    \item is a connected, smooth, complex, \emph{compact} orbifold of dimension $\dim(\Mbar_{g,n}) = 3g - 3 + n$;

    \item it contains $\M_{g,n}$ as an open dense subset.
  \end{itemize}
  The set $\de\Mbar_{g,n} = \Mbar_{g,n} \setminus \M_{g,n}$ is called the \emph{boundary} of the moduli space.
\end{theorem}

Now that we have a compact space, we can safely talk about integration. More generally, we have a nice (co)homology algebra
\begin{equation}
  \bigl(
    H_{\bullet}(\Mbar_{g,n},\Q), \frown
  \bigr)
  \qquad\text{and}\qquad
  \bigl(
    H^{\bullet}(\Mbar_{g,n},\Q), \smile
  \bigr) ,
\end{equation}
where the algebra structure is with respect to the cap/cup product (corresponding to intersection of subvarieties/wedge of differential forms respectively). The $\Q$ coefficients are due to the orbifold structure, and one can safely take $\C$ coefficients if they prefer. The two are dual via Poincaré duality:
\begin{equation}
  H^{k}(\Mbar_{g,n},\Q) \cong H_{2(3g-3+n) - k}(\Mbar_{g,n},\Q) \,.
\end{equation}
Most importantly, we have a well-defined fundamental class against which we can integrate cohomology classes to get a number:
\begin{equation}
  \int_{\Mbar_{g,n}} \alpha \in \Q \,,
  \qquad\qquad
  \alpha \in H^{2(3g-3+n)}(\Mbar_{g,n},\Q) \,.
\end{equation}
Since taking cap products in cohomology (i.e. wedges of differential forms) corresponds to taking cup products in homology (i.e. intersection of subvarieties), the theory of integration on compact moduli spaces is often called \emph{intersection theory}.

\subsection{Stratification and tautological maps}
Before discussing the cohomology of $\Mbar_{g,n}$ and its intersection theory further, let us analyse the compactification in more detail. The main picture to keep in mind is the following: most of the points of $\Mbar_{g,n}$ are smooth Riemann surfaces that live on $\M_{g,n}$, but by contracting cycles we produce stable singular Riemann surfaces that live on the boundary $\de\Mbar_{g,n}$. By performing this procedure once, we create a single node. By repeatedly performing this operation, we create Riemann surfaces that are more and more singular. See \cref{fig:Mgn} for an illustration.
\begin{figure}
  \centering
  \begin{tikzpicture}[x=1pt,y=1pt,scale=.75]
    \draw[thick,fill=gray!10]
    (152.806, 544.676).. controls (178.9353, 576.2253) and (202.6667, 618.6667) .. (224, 672)
    .. controls (296, 666.6667) and (360, 666.6667) .. (416, 672)
    .. controls (405.3333, 618.6667) and (400.2337, 576.4947) .. (400.701, 545.484)
    .. controls (325.567, 549.828) and (242.9353, 549.5587) ..(152.806, 544.676);



    
    \node at (224, 576) {\small$\times$};
    \node at (256, 624) {\small$\times$};
    \node at (336, 668) {\small$\times$};
    \node at (403.6515, 595.9072) {\small$\times$};
    \node at (416, 672) {\small$\times$};

    \draw[thick,densely dotted] (212, 580).. controls (188, 594.6667) and (173.3333, 609.3333) .. (168, 628);
    \draw[thick,densely dotted] (256, 676) -- (256, 632);
    \draw[thick,densely dotted] (336, 676).. controls (338.6667, 686.6667) and (342.6667, 694.6667) .. (348, 700);
    \draw[thick,densely dotted] (427, 687).. controls (436, 698.6667) and (446.6667, 709.3333) .. (460, 712);
    \draw[thick,densely dotted] (412, 594).. controls (420, 592) and (426.6667, 593.3333) .. (432, 596);

    \draw[thick] (186.3186, 656.2929).. controls (187.6519, 665.6262) and (179.6519, 676.2929) .. (179.6519, 686.9596).. controls (179.6519, 697.6262) and (187.6519, 708.2929) .. (184.9852, 718.9596).. controls (182.3186, 729.6262) and (168.9852, 740.2929) .. (156.9852, 738.9596).. controls (144.9852, 737.6262) and (134.3186, 724.2929) .. (132.9852, 713.6262).. controls (131.6519, 702.9596) and (139.6519, 694.9596) .. (140.9852, 682.9596).. controls (142.3186, 670.9596) and (136.9852, 654.9596) .. (139.6519, 645.6262).. controls (142.3186, 636.2929) and (152.9852, 633.6262) .. (163.6519, 636.2929).. controls (174.3186, 638.9596) and (184.9852, 646.9596) .. cycle;
    \draw[thick] (160.1755, 654.7269).. controls (170.8422, 665.3935) and (170.8422, 673.3935) .. (160.1755, 678.7269);
    \draw[thick] (164.2625, 659.3569).. controls (156.1755, 662.7269) and (156.1755, 670.7269) .. (165.0445, 675.4519);
    \node at (144, 656) {\tiny${\bullet}$};
    \node at (176, 724) {\tiny${\bullet}$};
    \draw[thick] (270.8848, 706.294).. controls (272.0287, 717.6923) and (268.4533, 730.0557) .. (270.6983, 738.7542).. controls (272.9433, 747.4527) and (281.0087, 752.4863) .. (278.8052, 761.9747).. controls (276.6017, 771.463) and (264.1293, 785.406) .. (252.5598, 785.7108).. controls (240.9903, 786.0157) and (230.3237, 772.6823) .. (231.1022, 762.7093).. controls (231.8807, 752.7363) and (244.1043, 746.1237) .. (245.7065, 736.6892).. controls (247.3087, 727.2547) and (238.2893, 714.9983) .. (238.5705, 704.0847).. controls (238.8517, 693.171) and (248.4333, 683.6) .. (256.1557, 684.0312).. controls (263.878, 684.4623) and (269.741, 694.8957) .. cycle;
    \draw[thick] (248, 756).. controls (264, 756) and (264, 768) .. (252, 776);
    \draw[thick] (253.604, 756.614).. controls (248, 764) and (248, 772) .. (256.037, 772.738);
    \node at (240, 760) {\tiny${\bullet}$};
    \node at (264, 720) {\tiny${\bullet}$};
    \draw[thick] (384, 744) .. controls (384, 728) and (364, 736) .. (352, 737.3333) .. controls (340, 738.6667) and (336, 733.3333) .. (336.6667, 726) .. controls (337.3333, 718.6667) and (342.6667, 709.3333) .. (353.3333, 705.3333) .. controls (364, 701.3333) and (380, 702.6667) .. (392, 709.3333) .. controls (404, 716) and (412, 728) .. (412.6667, 740) .. controls (413.3333, 752) and (406.6667, 764) .. (396.6667, 770.6667) .. controls (386.6667, 777.3333) and (373.3333, 778.6667) .. (362.6667, 777.3333) .. controls (352, 776) and (344, 772) .. (340.9983, 763.1847) .. controls (337.9967, 754.3695) and (339.9933, 740.739) .. (345.3266, 738.0723) .. controls (350.66, 735.4056) and (359.33, 743.7028) .. (366.9681, 745.995) .. controls (374.6062, 748.2872) and (381.2124, 744.5743) .. (381.0369, 736.2876);
    \node at (388, 764) {\tiny${\bullet}$};
    \node at (396, 724) {\tiny${\bullet}$};
    \draw[thick] (507.9455, 731.9235) .. controls (500.0319, 721.5632) and (489.9286, 724.6759) .. (482.2394, 726.9342) .. controls (474.5502, 729.1925) and (469.2751, 730.5962) .. (466.9344, 727.1871) .. controls (464.5937, 723.778) and (465.1873, 715.556) .. (473.1873, 710.2227) .. controls (481.1873, 704.8893) and (496.5937, 702.4447) .. (507.6302, 705.889) .. controls (518.6667, 709.3333) and (525.3333, 718.6667) .. (526, 728) .. controls (526.6667, 737.3333) and (521.3333, 746.6667) .. (511.3333, 752.6667) .. controls (501.3333, 758.6667) and (486.6667, 761.3333) .. (477.4062, 756.7671) .. controls (468.1457, 752.2009) and (464.2914, 740.4018) .. (465.6247, 734.4018) .. controls (466.958, 728.4018) and (473.479, 728.2009) .. (479.9455, 731.3101) .. controls (486.412, 734.4193) and (492.824, 740.8385) .. (497.8877, 739.567) .. controls (502.9514, 738.2956) and (506.6667, 729.3333) .. (505.6562, 729.3443);
    \draw[thick] (540.5888, 714.2895).. controls (552, 709.3333) and (568, 714.6667) .. (573.3333, 725.3333).. controls (578.6667, 736) and (573.3333, 752) .. (561.9222, 756.9562).. controls (550.511, 761.9123) and (533.022, 755.8247) .. (527.6887, 745.158).. controls (522.3553, 734.4913) and (529.1777, 719.2457) .. cycle;
    \node at (556, 752) {\tiny${\bullet}$};
    \node at (560, 724) {\tiny${\bullet}$};
    \draw[thick] (451.3333, 596.6667).. controls (437.3333, 597.3333) and (422.6667, 606.6667) .. (420.6667, 617.3333).. controls (418.6667, 628) and (429.3333, 640) .. (442.6667, 643.3333).. controls (456, 646.6667) and (472, 641.3333) .. (480.6667, 634.6667).. controls (489.3333, 628) and (490.6667, 620) .. (484.6667, 612).. controls (478.6667, 604) and (465.3333, 596) .. cycle;
    \draw[thick] (539.3333, 616.6667).. controls (541.3333, 626.6667) and (530.6667, 641.3333) .. (518, 643.3333).. controls (505.3333, 645.3333) and (490.6667, 634.6667) .. (488.6667, 624.6667).. controls (486.6667, 614.6667) and (497.3333, 605.3333) .. (510, 603.3333).. controls (522.6667, 601.3333) and (537.3333, 606.6667) .. cycle;
    \draw[thick] (484, 608)-- (484, 608);
    \draw[thick] (440, 624).. controls (444, 616) and (460, 612) .. (468, 628);
    \draw[thick] (444.592, 619.204).. controls (448, 628) and (460, 628) .. (465.581, 624.048);
    \node at (520, 636) {\tiny${\bullet}$};
    \node at (524, 616) {\tiny${\bullet}$};
  \end{tikzpicture}
  \caption{An illustration of the compactified moduli space $\Mbar_{g,n}$.}
  \label{fig:Mgn}
\end{figure}

As an example, consider the space $\Mbar_{0,4}$. On the boundary $\de\Mbar_{0,4}$ we find the singular Riemann surface made of two $\P^1$'s  glued together to form a node and each with two marked points. These arise from a smooth rational curve with four marked points by contracting a cycle separating the marked points into two-plus-two. We have three possible configurations, corresponding to the three possible ways of splitting $(p_1,p_2,p_3,p_4)$ into two disjoint sets containing two points each (see \cref{fig:boundary:04}). Notice that all such stable Riemann surfaces have no moduli: each rational component of their normalisations has three special points (the two marked points and a branch of the node), which can always be brought to $(0,1,\infty)$. Another way of saying it is that we can realise each of the above stable Riemann surfaces as the point $\mc{M}_{0,3} \times \mc{M}_{0,3}$. Recalling that $\M_{0,4} = \P^1 \setminus \set{0,1,\infty}$, we obtain that
\begin{equation}
  \Mbar_{0,4} = \M_{0,4} \sqcup (\M_{0,3} \times \M_{0,3})^{\sqcup 3} = \P^1 \,,
\end{equation}
which is indeed compact.

\begin{figure}[b]
  \centering
  \begin{tikzpicture}[scale=.85]
    \draw[thick] (-1,0) circle[radius=1cm];
    \draw[thick] (1,0) circle[radius=1cm];

    \node at (-1.5,.5) {\tiny$\bullet$};
    \node at ($(-1.5,.5) + (135:.65)$) {$p_1$};
    \node at (-1.5,-.5) {\tiny$\bullet$};
    \node at ($(-1.5,-.5) + (-135:.65)$) {$p_2$};

    \node at (1.5,.5) {\tiny$\bullet$};
    \node at ($(1.5,.5) + (45:.65)$) {$p_3$};
    \node at (1.5,-.5) {\tiny$\bullet$};
    \node at ($(1.5,-.5) + (-45:.65)$) {$p_4$};

    \begin{scope}[xshift=6cm]
      \draw[thick] (-1,0) circle[radius=1cm];
      \draw[thick] (1,0) circle[radius=1cm];

      \node at (-1.5,.5) {\tiny$\bullet$};
      \node at ($(-1.5,.5) + (135:.65)$) {$p_1$};
      \node at (-1.5,-.5) {\tiny$\bullet$};
      \node at ($(-1.5,-.5) + (-135:.65)$) {$p_3$};

      \node at (1.5,.5) {\tiny$\bullet$};
      \node at ($(1.5,.5) + (45:.65)$) {$p_2$};
      \node at (1.5,-.5) {\tiny$\bullet$};
      \node at ($(1.5,-.5) + (-45:.65)$) {$p_4$};
    \end{scope}

    \begin{scope}[xshift=12cm]
      \draw[thick] (-1,0) circle[radius=1cm];
      \draw[thick] (1,0) circle[radius=1cm];

      \node at (-1.5,.5) {\tiny$\bullet$};
      \node at ($(-1.5,.5) + (135:.65)$) {$p_1$};
      \node at (-1.5,-.5) {\tiny$\bullet$};
      \node at ($(-1.5,-.5) + (-135:.65)$) {$p_4$};

      \node at (1.5,.5) {\tiny$\bullet$};
      \node at ($(1.5,.5) + (45:.65)$) {$p_2$};
      \node at (1.5,-.5) {\tiny$\bullet$};
      \node at ($(1.5,-.5) + (-45:.65)$) {$p_3$};
    \end{scope}
  \end{tikzpicture}
  \caption{The three points on the boundary $\de\Mbar_{0,4}$.}
  \label{fig:boundary:04}
\end{figure}

As for $\Mbar_{1,1}$, the only element in the boundary $\de\Mbar_{1,1}$ is the pinched torus with a marked point encountered before. Again, the pinched torus has no moduli, as its normalisation is a rational curve with three marked points. However, the pinched torus has $\Z_2$ as an automorphism group. Another way of saying it is to realise it as $\mc{M}_{0,3}/\Z_2$. This gives
\begin{equation}
  \Mbar_{1,1}
  =
  \M_{1,1} \sqcup (\mc{M}_{0,3}/\Z_2) \,,
\end{equation}
which is topologically a $\P^1$ but with orbifold structure given by a point of automorphism $\Z_6$, a point of automorphism $\Z_4$, and all other points of automorphism $\Z_2$.

It should be clear from the above examples that the compactification of $\Mbar_{g,n}$ has a kind of recursive structure, obtained by pinching cycles and reducing the topology of the Riemann surface by breaking it up into pieces. We can keep track of this via certain graphs. Consider \cref{fig:stable:graph} for an illustration.

\begin{figure}[t]
  \centering
  \begin{tikzpicture}[x=1pt,y=1pt,xscale=.45,,yscale=.375]
    \draw[thick] (221.3465, 685.3958) .. controls (220, 672) and (207.6712, 660.4984) .. (196.9497, 658.5815) .. controls (186.2282, 656.6645) and (177.1141, 664.3323) .. (165.8904, 668.8328) .. controls (154.6667, 673.3333) and (141.3333, 674.6667) .. (134.3665, 662.9682) .. controls (127.3997, 651.2697) and (126.7993, 626.5393) .. (142.7993, 609.8727) .. controls (158.7993, 593.206) and (191.3997, 584.603) .. (223.0332, 586.3015) .. controls (254.6667, 588) and (285.3333, 600) .. (305.3333, 620.6667) .. controls (325.3333, 641.3333) and (334.6667, 670.6667) .. (333.3333, 700.6667) .. controls (332, 730.6667) and (320, 761.3333) .. (295.5366, 779.3405) .. controls (271.0732, 797.3476) and (234.1464, 802.6952) .. (204.1464, 796.0285) .. controls (174.1464, 789.3618) and (151.0732, 770.6809) .. (139.5463, 748.0344) .. controls (128.0195, 725.3878) and (128.0389, 698.7756) .. (136.5474, 684.7451) .. controls (145.0559, 670.7146) and (162.0534, 669.266) .. (174.0282, 674.5308) .. controls (186.003, 679.7957) and (192.955, 691.7741) .. (216.6159, 672.6254);
    \draw[thick] (350.8209, 613.5285) .. controls (384, 608) and (416, 592) .. (432, 568) .. controls (448, 544) and (448, 512) .. (424, 496) .. controls (400, 480) and (352, 480) .. (318.8209, 501.5285) .. controls (285.6419, 523.0571) and (267.2837, 566.1142) .. (275.2837, 590.1142) .. controls (283.2837, 614.1142) and (317.6419, 619.0571) .. cycle;
    \draw[thick] (444.2799, 615.5693) .. controls (464, 629.3333) and (480, 650.6667) .. (485.3333, 677.3333) .. controls (490.6667, 704) and (485.3333, 736) .. (472, 762.6667) .. controls (458.6667, 789.3333) and (437.3333, 810.6667) .. (413.377, 818.4799) .. controls (389.4206, 826.2932) and (362.8411, 820.5863) .. (343.7106, 804.3723) .. controls (324.58, 788.1584) and (312.8983, 761.4373) .. (320.8546, 748.2907) .. controls (328.811, 735.1441) and (356.4055, 735.5721) .. (361.7682, 727.1387) .. controls (367.131, 718.7052) and (350.2619, 701.4105) .. (340.6764, 688.5937) .. controls (331.0909, 675.7769) and (328.7889, 667.4381) .. (333.5236, 662.4668) .. controls (338.2584, 657.4956) and (350.0298, 655.8919) .. (358.5618, 644.278) .. controls (367.0939, 632.6641) and (372.3865, 611.0399) .. (386.753, 603.3251) .. controls (401.1194, 595.6104) and (424.5597, 601.8052) .. cycle;
    \draw[thick] (184, 740) .. controls (204, 720) and (260, 724) .. (256, 760);
    \draw[thick] (204.25, 729.625) .. controls (208, 748) and (236, 752) .. (254.1784, 746.171);
    \draw[thick] (332, 588) .. controls (324, 556) and (360, 524) .. (404, 536);
    \draw[thick] (335.1198, 561.5158) .. controls (360, 568) and (380, 552) .. (387.7078, 533.4673);
    \draw[thick] (368, 784) .. controls (392, 748) and (424, 756) .. (424, 784);
    \draw[thick] (396, 696) .. controls (388, 668) and (416, 640) .. (448, 680);
    \draw[thick] (395.2528, 679.7994) .. controls (412, 688) and (432, 684) .. (435.9985, 667.8575);
    \draw[thick] (382.841, 767.6393) .. controls (384, 788) and (408, 792) .. (421.1756, 771.0678);
    \node at (160, 616) {\tiny$\bullet$};
    \node at (316, 528) {\tiny$\bullet$};
    \node at (388, 500) {\tiny$\bullet$};
    \node at (145, 580) {$p_1$};
    \node at (284, 500) {$p_2$};
    \node at (396, 464) {$p_3$};
  \end{tikzpicture}
  \hspace{2cm}
  \begin{tikzpicture}
    \draw[thick] (150:1.5) -- ($(150:1.5) + (240:1.25)$);
    \node at ($(150:1.5) + (240:1.5)$) {$1$};
    \draw[thick] (-90:1.5) -- ($(-90:1.5) + (-120:1.25)$);
    \node at ($(-90:1.5) + (-120:1.5)$) {$2$};
    \draw[thick] (-90:1.5) -- ($(-90:1.5) + (-60:1.25)$);
    \node at ($(-90:1.5) + (-60:1.5)$) {$3$};

    \draw[thick] (30:1.5) to [bend right=30] (150:1.5) -- (-90:1.5) -- cycle;
    \draw[thick] (30:1.5) to [bend left=30] (150:1.5);


    \draw[thick,shift={(150:1.5)},rotate=60] (0,0) to[out=45,in=-90] (.45,.75) to[out=90,in=0] (0,1.2) to[out=180,in=90] (-.45,.75) to[out=-90,in=135] (0,0);

    \draw[thick,fill=white] (30:1.5) circle[radius=.3cm];
    \node at (30:1.5) {$2$};
    \draw[thick,fill=white] (150:1.5) circle[radius=.3cm];
    \node at (150:1.5) {$1$};
    \draw[thick,fill=white] (-90:1.5) circle[radius=.3cm];
    \node at (-90:1.5) {$1$};

  \end{tikzpicture}
  \caption{A stable Riemann surface and the associated stable graph.}
  \label{fig:stable:graph}
\end{figure}

\begin{definition}
  The \emph{stable graph} associated with a stable Riemann surface $(\Sigma,p_1,\dots,p_n) \in \Mbar_{g,n}$ is the graph $\Gamma$ obtained by associating:
  \begin{itemize}
    \item a vertex $v$ to each component of the normalisation, decorated by the genus $g(v)$ of the component;
    \item a leaf to each marked point $p_i$, labelled by $i$ accordingly;
    \item an edge to each node.
  \end{itemize}
  The genus of a stable graph $\Gamma$ is
  \begin{equation}
    g(\Gamma) = \sum_{v\in V(\Gamma)} g(v) + h^1(\Gamma) \,,
  \end{equation}
  where $V(\Gamma)$ is the set of the vertices of the graph and $h^1(\Gamma)$ denotes the first Betti number (i.e. the number of faces) of $\Gamma$. It coincides with the genus of $\Sigma$. We also denote by $E(\Gamma)$ the set of edges and by $n(v)$ the valency of the vertex $v$ (that is, the number of leaves and half-edges incident to $v$). The latter corresponds to the number of special points (that is, marked points and branches of nodes) on the component corresponding to the vertex $v$.
\end{definition}

We remark that the stability condition implies that $2g(v) - 2 + n(v) > 0$ for all $v \in V(\Gamma)$. This guarantees that for each $(g,n)$, called the type, there are only finitely many stable graphs of genus $g$ with $n$ leaves. Such stable graphs provide a \emph{stratification} of $\Mbar_{g,n}$: for a given $\Gamma$ of type $(g,n)$, set
\begin{equation}
  \M_{\Gamma}
  =
  \Set{
    (\Sigma,p_1,\dots,p_n) \in \Mbar_{g,n}
    |
    \substack{
      \displaystyle \Gamma \text{ is the stable graph} \\[.5ex]
      \displaystyle \text{associated with }(\Sigma,p_1,\dots,p_n) 
    }
  } .
\end{equation}
Then we get the stratification
\begin{equation}
  \Mbar_{g,n}
  =
  \bigsqcup_{\substack{\Gamma \text{ of type }(g,n)}} \M_{\Gamma} \,.
\end{equation}
We have already analysed thoroughly the cases of $\Mbar_{0,4}$ and $\Mbar_{1,1}$, whose stable graphs are given in \cref{fig:stble:04:11}. Another example is that of $\Mbar_{2}$, see \cref{fig:stble:20}. There we drew the strata corresponding to the type of stable Riemann surface dual to the graph and on different levels according to the number of edges. Note that contraction of cycles is dual to contraction of edges.

\begin{figure}[t]
  \begin{tikzpicture}[scale=.625]
    \node at (-3,0) {$(0,4)$:};

    \draw[thick] (0,0) -- (135:1);
    \node at (135:1.3) {$1$};
    \draw[thick] (0,0) -- (-135:1);
    \node at (-135:1.3) {$2$};
    \draw[thick] (0,0) -- (-45:1);
    \node at (-45:1.3) {$3$};
    \draw[thick] (0,0) -- (45:1);
    \node at (45:1.3) {$4$};

    \draw[thick,fill=white] (0,0) circle[radius=.3cm];
    \node at (0,0) {\small$0$};

    \begin{scope}[xshift = 4.75cm]
      \draw[thick] (-1,0) -- (1,0);

      \draw[thick] (-1,0) -- ($(-1,0) + (135:1)$);
      \node at ($(-1,0) + (135:1.3)$) {$1$};
      \draw[thick] (-1,0) -- ($(-1,0) + (-135:1)$);
      \node at ($(-1,0) + (-135:1.3)$) {$2$};
      
      \draw[thick] (1,0) -- ($(1,0) + (45:1)$);
      \node at ($(1,0) + (45:1.3)$) {$3$};
      \draw[thick] (1,0) -- ($(1,0) + (-45:1)$);
      \node at ($(1,0) + (-45:1.3)$) {$4$};

      \draw[thick,fill=white] (-1,0) circle[radius=.3cm];
      \node at (-1,0) {\small$0$};
      \draw[thick,fill=white] (1,0) circle[radius=.3cm];
      \node at (1,0) {\small$0$};

      \begin{scope}[xshift = 5.75cm]
        \draw[thick] (-1,0) -- (1,0);

        \draw[thick] (-1,0) -- ($(-1,0) + (135:1)$);
        \node at ($(-1,0) + (135:1.3)$) {$1$};
        \draw[thick] (-1,0) -- ($(-1,0) + (-135:1)$);
        \node at ($(-1,0) + (-135:1.3)$) {$3$};
        
        \draw[thick] (1,0) -- ($(1,0) + (45:1)$);
        \node at ($(1,0) + (45:1.3)$) {$2$};
        \draw[thick] (1,0) -- ($(1,0) + (-45:1)$);
        \node at ($(1,0) + (-45:1.3)$) {$4$};

        \draw[thick,fill=white] (-1,0) circle[radius=.3cm];
        \node at (-1,0) {\small$0$};
        \draw[thick,fill=white] (1,0) circle[radius=.3cm];
        \node at (1,0) {\small$0$};

        \begin{scope}[xshift = 5.75cm]
          \draw[thick] (-1,0) -- (1,0);

          \draw[thick] (-1,0) -- ($(-1,0) + (135:1)$);
          \node at ($(-1,0) + (135:1.3)$) {$1$};
          \draw[thick] (-1,0) -- ($(-1,0) + (-135:1)$);
          \node at ($(-1,0) + (-135:1.3)$) {$4$};
          
          \draw[thick] (1,0) -- ($(1,0) + (45:1)$);
          \node at ($(1,0) + (45:1.3)$) {$2$};
          \draw[thick] (1,0) -- ($(1,0) + (-45:1)$);
          \node at ($(1,0) + (-45:1.3)$) {$3$};

          \draw[thick,fill=white] (-1,0) circle[radius=.3cm];
          \node at (-1,0) {\small$0$};
          \draw[thick,fill=white] (1,0) circle[radius=.3cm];
          \node at (1,0) {\small$0$};
        \end{scope}
      \end{scope}
    \end{scope}

    \begin{scope}[xshift=.8cm, yshift=-2.5cm]
      \node at (-3.8,0) {$(1,1)$:};
      \draw[thick] (0,0) -- (180:1.5);
      \node at (180:1.7) {$1$};

      \draw[thick,fill=white] (0,0) circle[radius=.3cm];
      \node at (0,0) {\small$1$};

      \begin{scope}[xshift = 4.25cm]
        \draw[thick] (0,0) -- (180:1.5);
        \node at (180:1.7) {$1$};
      
        \draw[thick,rotate=-90] (0,0) to[out=45,in=-90] (.45,.75) to[out=90,in=0] (0,1.2) to[out=180,in=90] (-.45,.75) to[out=-90,in=135] (0,0);

        \draw[thick,fill=white] (0,0) circle[radius=.3cm];
        \node at (0,0) {\small$0$};
      \end{scope}
      \end{scope}
  \end{tikzpicture}
\caption{All stable graphs of type $(0,4)$ and $(1,1)$.}
\label{fig:stble:04:11}
\end{figure}

\begin{figure}[b]
  \begin{tikzpicture}[x=1pt,y=1pt,scale=.265]
    \node at (400, 880) {$\vphantom{\dim(\mc{M}_{\Gamma})}$};
    \draw[thick] (349.3333, 690.6667) .. controls (368, 693.3333) and (384, 714.6667) .. (384, 736) .. controls (384, 757.3333) and (368, 778.6667) .. (349.3333, 781.3333) .. controls (330.6667, 784) and (309.3333, 768) .. (288, 768) .. controls (266.6667, 768) and (245.3333, 784) .. (226.6667, 781.3333) .. controls (208, 778.6667) and (192, 757.3333) .. (192, 736) .. controls (192, 714.6667) and (208, 693.3333) .. (226.6667, 690.6667) .. controls (245.3333, 688) and (266.6667, 704) .. (288, 704) .. controls (309.3333, 704) and (330.6667, 688) .. cycle;
    \draw[thick] (220, 744) .. controls (228, 724) and (252, 724) .. (260, 744);
    \draw[thick] (223.8918, 737.1747) .. controls (234.6306, 747.0582) and (245.3106, 746.9846) .. (255.9319, 736.9537);
    \draw[thick] (316, 744) .. controls (324, 724) and (348, 724) .. (356, 744);
    \draw[thick] (319.892, 737.175) .. controls (330.6307, 747.0583) and (341.3107, 746.9847) .. (351.932, 736.954);
    \draw[thick] (221.3333, 546.6667) .. controls (240, 549.3333) and (256, 570.6667) .. (256, 592) .. controls (256, 613.3333) and (240, 634.6667) .. (221.3333, 637.3333) .. controls (202.6667, 640) and (181.3333, 624) .. (160, 624) .. controls (138.6667, 624) and (117.3333, 640) .. (98.6667, 637.3333) .. controls (80, 634.6667) and (64, 613.3333) .. (64, 592) .. controls (64, 570.6667) and (80, 549.3333) .. (98.6667, 546.6667) .. controls (117.3333, 544) and (138.6667, 560) .. (160, 560) .. controls (181.3333, 560) and (202.6667, 544) .. cycle;
    \draw[thick] (188, 600) .. controls (196, 580) and (220, 580) .. (228, 600);
    \draw[thick] (191.892, 593.175) .. controls (202.6307, 603.0583) and (213.3107, 602.9847) .. (223.932, 592.954);
    \draw[thick] (112, 596) .. controls (88, 576) and (76, 578) .. (70, 581) .. controls (64, 584) and (64, 588) .. (64, 591.3333) .. controls (64, 594.6667) and (64, 597.3333) .. (68, 599.6667) .. controls (72, 602) and (80, 604) .. (107.0627, 591.9835);
    \draw[thick] (416, 592) .. controls (424, 560) and (456, 552) .. (478.6667, 556) .. controls (501.3333, 560) and (514.6667, 576) .. (514.6667, 592) .. controls (514.6667, 608) and (501.3333, 624) .. (478.6667, 628) .. controls (456, 632) and (424, 624) .. (416, 592) .. controls (408, 624) and (376, 632) .. (353.3333, 628) .. controls (330.6667, 624) and (317.3333, 608) .. (317.3333, 592) .. controls (317.3333, 576) and (330.6667, 560) .. (353.3333, 556) .. controls (376, 552) and (408, 560) .. (416, 592);
    \draw[thick] (344, 600) .. controls (352, 580) and (376, 580) .. (384, 600);
    \draw[thick] (347.892, 593.175) .. controls (358.6307, 603.0583) and (369.3107, 602.9847) .. (379.932, 592.954);
    \draw[thick] (448, 600) .. controls (456, 580) and (480, 580) .. (488, 600);
    \draw[thick] (451.892, 593.175) .. controls (462.6307, 603.0583) and (473.3107, 602.9847) .. (483.932, 592.954);
    \draw[thick] (93.3333, 402.6667) .. controls (112, 405.3333) and (128, 426.6667) .. (128, 448) .. controls (128, 469.3333) and (112, 490.6667) .. (93.3333, 493.3333) .. controls (74.6667, 496) and (53.3333, 480) .. (32, 480) .. controls (10.6667, 480) and (-10.6667, 496) .. (-29.3333, 493.3333) .. controls (-48, 490.6667) and (-64, 469.3333) .. (-64, 448) .. controls (-64, 426.6667) and (-48, 405.3333) .. (-29.3333, 402.6667) .. controls (-10.6667, 400) and (10.6667, 416) .. (32, 416) .. controls (53.3333, 416) and (74.6667, 400) .. cycle;
    \draw[thick] (-16, 452) .. controls (-40, 432) and (-52, 434) .. (-58, 437) .. controls (-64, 440) and (-64, 444) .. (-64, 447.3333) .. controls (-64, 450.6667) and (-64, 453.3333) .. (-60, 455.6667) .. controls (-56, 458) and (-48, 460) .. (-20.937, 447.984);
    \draw[thick] (80, 452) .. controls (104, 432) and (116, 432) .. (122, 434.6667) .. controls (128, 437.3333) and (128, 442.6667) .. (128, 446.6667) .. controls (128, 450.6667) and (128, 453.3333) .. (124, 455.6667) .. controls (120, 458) and (112, 460) .. (85.0078, 447.9572);
    \draw[thick] (316, 456) .. controls (324, 436) and (348, 436) .. (356, 456);
    \draw[thick] (319.892, 449.175) .. controls (330.6307, 459.0583) and (341.3107, 458.9847) .. (351.932, 448.954);
    \draw[thick] (288, 448) .. controls (296, 416) and (328, 408) .. (350.6667, 412) .. controls (373.3333, 416) and (386.6667, 432) .. (386.6667, 448) .. controls (386.6667, 464) and (373.3333, 480) .. (350.6667, 484) .. controls (328, 488) and (296, 480) .. (288, 448) .. controls (280, 480) and (248, 488) .. (225.3333, 484) .. controls (202.6667, 480) and (189.3333, 464) .. (189.3333, 448) .. controls (189.3333, 432) and (202.6667, 416) .. (225.3333, 412) .. controls (248, 408) and (280, 416) .. (288, 448);
    \draw[thick] (237.3334, 452.0459) .. controls (213.3334, 432.0459) and (201.3334, 434.0459) .. (195.3334, 437.0459) .. controls (189.3334, 440.0459) and (189.3334, 444.0459) .. (189.3334, 447.3793) .. controls (189.3334, 450.7126) and (189.3334, 453.3793) .. (193.3334, 455.7126) .. controls (197.3334, 458.0459) and (205.3334, 460.0459) .. (232.3964, 448.0299);
    \draw[thick] (-34.6667, 258.6667) .. controls (-16, 261.3333) and (0, 282.6667) .. (0, 304) .. controls (0, 325.3333) and (-16, 346.6667) .. (-34.6667, 349.3333) .. controls (-53.3333, 352) and (-74.6667, 336) .. (-96, 336) .. controls (-117.3333, 336) and (-138.6667, 352) .. (-157.3333, 349.3333) .. controls (-176, 346.6667) and (-192, 325.3333) .. (-192, 304) .. controls (-192, 282.6667) and (-176, 261.3333) .. (-157.3333, 258.6667) .. controls (-138.6667, 256) and (-117.3333, 272) .. (-96, 272) .. controls (-74.6667, 272) and (-53.3333, 256) .. cycle;
    \draw[thick] (-96, 304) .. controls (-104, 320) and (-124, 324) .. (-144, 324) .. controls (-164, 324) and (-184, 320) .. (-192, 304);
    \draw[thick] (-192, 304) .. controls (-184, 288) and (-168, 284) .. (-148, 284) .. controls (-128, 284) and (-104, 288) .. (-96, 304);
    \draw[thick] (0, 304) .. controls (-8, 320) and (-28, 324) .. (-48, 324) .. controls (-68, 324) and (-88, 320) .. (-96, 304);
    \draw[thick] (-96, 304) .. controls (-88, 288) and (-72, 284) .. (-52, 284) .. controls (-32, 284) and (-8, 288) .. (0, 304);
    \draw[thick] (160, 304) .. controls (168, 272) and (200, 264) .. (222.6667, 268) .. controls (245.3333, 272) and (258.6667, 288) .. (258.6667, 304) .. controls (258.6667, 320) and (245.3333, 336) .. (222.6667, 340) .. controls (200, 344) and (168, 336) .. (160, 304) .. controls (152, 336) and (120, 344) .. (97.3333, 340) .. controls (74.6667, 336) and (61.3333, 320) .. (61.3333, 304) .. controls (61.3333, 288) and (74.6667, 272) .. (97.3333, 268) .. controls (120, 264) and (152, 272) .. (160, 304);
    \draw[thick] (109.333, 308.046) .. controls (85.333, 288.046) and (73.333, 290.046) .. (67.333, 293.046) .. controls (61.333, 296.046) and (61.333, 300.046) .. (61.333, 303.3793) .. controls (61.333, 306.7127) and (61.333, 309.3793) .. (65.333, 311.7127) .. controls (69.333, 314.046) and (77.333, 316.046) .. (104.396, 304.03);
    \draw[thick] (210.6667, 308) .. controls (234.6667, 288) and (246.6667, 288) .. (252.6667, 290.6667) .. controls (258.6667, 293.3333) and (258.6667, 298.6667) .. (258.6667, 302.6667) .. controls (258.6667, 306.6667) and (258.6667, 309.3333) .. (254.6667, 311.6667) .. controls (250.6667, 314) and (242.6667, 316) .. (215.6745, 303.957);
    \draw[-{>[scale=.75]}] (336, 672) -- (352, 656);
    \draw[-{>[scale=.75]}] (240, 672) -- (224, 656);
    \draw[-{>[scale=.75]}] (208, 528) -- (224, 512);
    \draw[-{>[scale=.75]}] (112, 528) -- (96, 512);
    \draw[-{>[scale=.75]}] (80, 384) -- (96, 368);
    \draw[-{>[scale=.75]}] (-16, 384) -- (-32, 368);
    \draw[-{>[scale=.75]}] (368, 528) -- (352, 512);
    \draw[-{>[scale=.75]}] (240, 384) -- (224, 368);
  \end{tikzpicture}
  \begin{tikzpicture}[x=1pt,y=1pt,scale=.265]
    \node at (600, 840) {$\dim(\mc{M}_{\Gamma})$}; 
    \node at (600, 736) {$3$}; 
    \filldraw[thick,fill=white] (288, 736) circle[radius=24];
    \node at (288, 736) {\small$2$};
    \draw[thick] (160, 592) .. controls (112, 640) and (88, 616) .. (88, 592) .. controls (88, 568) and (112, 544) .. (160, 592);
    \filldraw[thick,fill=white] (160, 592) circle[radius=24];
    \node at (160, 592) {\small$1$};
    \node at (600, 592) {$2$}; 
    \draw[thick] (352, 592) -- (480, 592);
    \filldraw[thick,fill=white] (352, 592) circle[radius=24];
    \node at (352, 592) {\small$1$};
    \filldraw[thick,fill=white] (480, 592) circle[radius=24];
    \node at (480, 592) {\small$1$};
    \draw[thick] (32, 448) .. controls (80, 496) and (104, 472) .. (104, 448) .. controls (104, 424) and (80, 400) .. (32, 448);
    \draw[thick] (32, 448) .. controls (-16, 496) and (-40, 472) .. (-40, 448) .. controls (-40, 424) and (-16, 400) .. (32, 448);
    \filldraw[thick,fill=white] (32, 448) circle[radius=24];
    \node at (32, 448) {\small$0$};
    \node at (600, 448) {$1$}; 
    \draw[thick] (240, 448) .. controls (192, 496) and (168, 472) .. (168, 448) .. controls (168, 424) and (192, 400) .. (240, 448);
    \draw[thick] (240, 448) -- (336, 448);
    \filldraw[thick,fill=white] (240, 448) circle[radius=24];
    \node at (240, 448) {\small$0$};
    \filldraw[thick,fill=white] (336, 448) circle[radius=24];
    \node at (336, 448) {\small$1$};
    \draw[thick] (-96, 248) .. controls (-53.3333, 285.3333) and (-53.3333, 322.6667) .. (-96, 360);
    \draw[thick] (-96, 248) .. controls (-138.6667, 285.3333) and (-138.6667, 322.6667) .. (-96, 360);
    \draw[thick] (-96, 248) -- (-96, 360);
    \filldraw[thick,fill=white] (-96, 360) circle[radius=24];
    \node at (-96, 360) {\small$0$};
    \filldraw[thick,fill=white] (-96, 248) circle[radius=24];
    \node at (-96, 248) {\small$0$};
    \draw[thick] (208, 304) .. controls (256, 352) and (280, 328) .. (280, 304) .. controls (280, 280) and (256, 256) .. (208, 304);
    \draw[thick] (112, 304) .. controls (64, 352) and (40, 328) .. (40, 304) .. controls (40, 280) and (64, 256) .. (112, 304);
    \draw[thick] (112, 304) -- (208, 304);
    \filldraw[thick,fill=white] (112, 304) circle[radius=24];
    \node at (112, 304) {\small$0$};
    \filldraw[thick,fill=white] (208, 304) circle[radius=24];
    \node at (208, 304) {\small$0$};
    \node at (600, 304) {$0$}; 
    \draw[{<[scale=.75]}-] (336, 672) -- (352, 656);
    \draw[{<[scale=.75]}-] (240, 672) -- (224, 656);
    \draw[{<[scale=.75]}-] (208, 528) -- (224, 512);
    \draw[{<[scale=.75]}-] (112, 528) -- (96, 512);
    \draw[{<[scale=.75]}-] (80, 384) -- (96, 368);
    \draw[{<[scale=.75]}-] (-16, 384) -- (-32, 368);
    \draw[{<[scale=.75]}-] (368, 528) -- (352, 512);
    \draw[{<[scale=.75]}-] (240, 384) -- (224, 368);
  \end{tikzpicture}
\caption{All possible stable Riemann surfaces of genus $2$, together with the associated stable graphs and the dimensions of the corresponding strata.}
\label{fig:stble:20}
\end{figure}

\noindent 

\begin{exercise}\leavevmode
  \begin{enumerate}
    \item
    List all strata of $\Mbar_{2,1}$.

    \item
    Consider a stable graph $\Gamma$ of type $(g,n)$. Show that the dimension of the stratum is given by $\dim(\M_{\Gamma}) = \dim(\Mbar_{g,n}) - |E_{\Gamma}|$.

  \end{enumerate}
\end{exercise}

The fact that the strata of $\Mbar_{g,n}$ are parametrised by smaller-dimensional spaces \smash{$\M_{\Gamma}$} is sometimes called the \emph{recursive boundary structure} of $\Mbar_{g,n}$. It is one of the most important features of the moduli space of Riemann surfaces and the proofs of many results about $\Mbar_{g,n}$ (including the computation of integrals) use it in a very essential way.

One way of taking advantage of it is by defining \emph{gluing maps}. More precisely, for each stable graph $\Gamma$ of type $(g,n)$ we define
\begin{equation}
  \xi_{\Gamma} \colon
  \Mbar_{\Gamma} = \prod_{v \in V(\Gamma)} \Mbar_{g(v),n(v)}
  \longrightarrow
  \Mbar_{g,n} \,,
\end{equation}
which sends the stable Riemann surface $((\Sigma_v)_{v \in V(\Gamma)}, (q_{h},q_{h'})_{e=(h,h') \in E(\Gamma)}, p_1,\dots,p_n)$ to the stable Riemann surface $(\Sigma,p_1,\dots,p_n)$ obtained by gluing all pairs $(q_h,q_{h'})$ of points corresponding to pairs $e = (h, h')$ of half-edges forming an edge $e$ of $\Gamma$. The image of $\Mbar_{\Gamma}$ under $\xi_{\Gamma}$ coincides with the closure of $\M_{\Gamma}$.

The easiest case is that of a stable graph $\Gamma$ with a single edge $e$. We have two possible cases: the edge is non-separating (i.e. a loop) or it is.

\paragraph{Non-separating edge.}
It corresponds to the following stable graph:
\begin{equation}
  \begin{tikzpicture}[baseline]
    \draw[thick] (0,0) -- (145:1.5);
    \node at (145:1.7) {$1$};
    \node[rotate=90] at (180:1.35) {$\cdots$};
    \draw[thick] (0,0) -- (-145:1.5);
    \node at (-145:1.7) {$n$};
    
    \draw[thick,rotate=-90] (0,0) to[out=45,in=-90] (.45,.75) to[out=90,in=0] (0,1.2) to[out=180,in=90] (-.45,.75) to[out=-90,in=135] (0,0);

    \draw[thick,fill=white] (0,0) circle[radius=.4cm];
    \node at (0,0) {\tiny$g-1$};
  \end{tikzpicture}
  \;\;.
\end{equation}
Thus, the gluing map, called the gluing map of \emph{non-separating kind}, is given by
\begin{equation}
  \rho \colon
  \Mbar_{g-1,n+2}
  \longrightarrow
  \Mbar_{g,n} \,,
  \qquad\text{e.g.}\qquad
  \vcenter{\hbox{
    \begin{tikzpicture}[x=1pt,y=1pt,scale=.3]
      \fill[fill=gray!10] (130.7979, 642.2304) .. controls (136.8066, 646.0611) and (142.8162, 648.5197) .. (148.8269, 649.6064) -- (126.2269, 647.1274) -- cycle;
      \draw[thick] (123.9409, 650.2604) .. controls (139.9409, 626.2604) and (171.9409, 626.2604) .. (187.9409, 650.2604);
      \draw[thick] (130.7979, 642.2304) .. controls (147.5599, 652.9171) and (164.3322, 652.9264) .. (181.1149, 642.2584);
      \draw[thick] (108, 606.6667) .. controls (130.6667, 594.6667) and (173.3333, 589.3333) .. (198.9047, 597.7485) .. controls (224.476, 606.1637) and (232.952, 628.3273) .. (228.952, 645.6607) .. controls (224.952, 662.994) and (208.476, 675.497) .. (186.9047, 681.7485) .. controls (165.3333, 688) and (138.6667, 688) .. (120, 682.6667) .. controls (101.3333, 677.3333) and (90.6667, 666.6667) .. (86.6667, 652) .. controls (82.6667, 637.3333) and (85.3333, 618.6667) .. cycle;
      \node at (216, 648) {\tiny$\bullet$};
      \node at (104, 656) {\tiny$\bullet$};
      \node at (210, 616) {\tiny$\bullet$};
      \node at (80, 685) {\small$p_1$};
      \node at (232, 685) {\small$q$};
      \node at (224, 580) {\small$q'$};
    \end{tikzpicture}
  }}
  \longmapsto
  \vcenter{\hbox{
    \begin{tikzpicture}[x=1pt,y=1pt,scale=.3]
      \fill[fill=gray!10] (441.7167, 646.4241) .. controls (439.4751, 643.8509) and (437.2825, 640.1825) .. (435.139, 635.419) -- (432.5828, 641.9544) -- cycle;
      \fill[fill=gray!10] (346.798, 642.2304) .. controls (352.8067, 646.0611) and (358.8163, 648.5197) .. (364.827, 649.6064) -- (342.227, 647.1274) -- cycle;
      \draw[thick] (339.941, 650.2604) .. controls (355.941, 626.2604) and (387.941, 626.2604) .. (403.941, 650.2604);
      \draw[thick] (346.798, 642.2304) .. controls (363.56, 652.9171) and (380.3323, 652.9264) .. (397.115, 642.2584);
      \draw[thick] (332, 606.6667) .. controls (354.6667, 594.6667) and (397.3333, 589.3333) .. (425.3333, 598.6667) .. controls (453.3333, 608) and (466.6667, 632) .. (462.6667, 649.3333) .. controls (458.6667, 666.6667) and (437.3333, 677.3333) .. (413.3333, 682.6667) .. controls (389.3333, 688) and (362.6667, 688) .. (344, 682.6667) .. controls (325.3333, 677.3333) and (314.6667, 666.6667) .. (310.6667, 652) .. controls (306.6667, 637.3333) and (309.3333, 618.6667) .. cycle;
      \node at (328, 656) {\tiny$\bullet$};
      \draw[thick] (432, 648) .. controls (432, 624) and (456, 624) .. (463.2164, 640.0147);
      \draw[thick] (463.2164, 640.0147) .. controls (453.0721, 656.0049) and (443.713, 654.473) .. (435.1389, 635.4191);
      \draw (435.1389, 635.4191) -- (435.1389, 635.4191);
      \node at (304, 685) {\small$p_1$};
    \end{tikzpicture}
  }} \;\;.
\end{equation}
To be pedantic, $\rho$ should depend on $(g,n)$. We omit the dependence for lighter notation.

\paragraph{Separating edge.}
It corresponds to the following stable graph:
\begin{equation}
  \begin{tikzpicture}[baseline]
    \draw[thick] (0,0) -- (155:1.2);
    \node[rotate=90] at (180:.8) {$\cdots$};
    \draw[thick] (0,0) -- (-155:1.2);
    \node at (180:1.5) {$I_1$};

    \draw[thick] (1.5,0) -- ($(1.5,0) + (25:1.2)$);
    \node[rotate=90] at ($(1.5,0) + (0:.8)$) {$\cdots$};
    \draw[thick] (1.5,0) -- ($(1.5,0) + (-25:1.2)$);
    \node at ($(1.5,0) + (0:1.5)$) {$I_2$};

    \draw[thick] (0,0) -- (1.5,0);

    \draw[thick,fill=white] (0,0) circle[radius=.3cm];
    \node at (0,0) {\footnotesize$g_1$};

    \draw[thick,fill=white] (1.5,0) circle[radius=.3cm];
    \node at (1.5,0) {\footnotesize$g_2$};
  \end{tikzpicture}
\end{equation}
where $g = g_1 + g_2$ is a splitting of the genus and $I_1 \sqcup I_2 = \set{p_1,\dots,p_n}$ is a splitting of the marked points.
Thus, the corresponding gluing map, called the gluing map of \emph{separating kind}, is given by
\begin{equation}
  \sigma \colon
  \Mbar_{g_1,1+|I_1|} \times \Mbar_{g_2,1+|I_2|}
  \longrightarrow
  \Mbar_{g,n} \,,
  \quad\text{e.g.}\quad
  \vcenter{\hbox{
    \begin{tikzpicture}[x=1pt,y=1pt,scale=.3]
      \fill[fill=gray!10] (66.798, 642.2304) .. controls (72.8067, 646.0611) and (78.8163, 648.5197) .. (84.827, 649.6064) -- (62.227, 647.1274) -- cycle;
      \draw[thick] (59.941, 650.2604) .. controls (75.941, 626.2604) and (107.941, 626.2604) .. (123.941, 650.2604);
      \draw[thick] (66.798, 642.2304) .. controls (83.56, 652.9171) and (100.3323, 652.9264) .. (117.115, 642.2584);
      \draw[thick] (44, 606.6667) .. controls (66.6667, 594.6667) and (109.3333, 589.3333) .. (134.9047, 597.7485) .. controls (160.476, 606.1637) and (168.952, 628.3273) .. (164.952, 645.6607) .. controls (160.952, 662.994) and (144.476, 675.497) .. (122.9047, 681.7485) .. controls (101.3333, 688) and (74.6667, 688) .. (56, 682.6667) .. controls (37.3333, 677.3333) and (26.6667, 666.6667) .. (22.6667, 652) .. controls (18.6667, 637.3333) and (21.3333, 618.6667) .. cycle;
      \node at (152, 640) {\tiny$\bullet$};
      \node at (40, 656) {\tiny$\bullet$};
      \node at (16, 685) {\small$p_1$};
      \node at (160, 690) {\small$q$};
      \node at (180, 640) {$\vphantom{q}$,};
      \fill[fill=gray!10] (220.994, 598.93)  .. controls (244.3313, 604.9767) and (264.7397, 613.2347) .. (282.219, 623.704)  arc[start angle=21.106, end angle=114.8333, x radius=45.2548, y radius=-45.2548]  -- cycle;
      \draw[thick] (240, 640) circle[radius=45.2548];
      \node at (208, 640) {\tiny$\bullet$};
      \node at (256, 672) {\tiny$\bullet$};
      \node at (256, 608) {\tiny$\bullet$};
      \node at (192, 590) {\small$q'$};
      \node at (272, 700) {\small$p_2$};
      \node at (272, 575) {\small$p_3$};
    \end{tikzpicture}
  }}
  \longmapsto\!\!
  \vcenter{\hbox{
    \begin{tikzpicture}[x=1pt,y=1pt,scale=.3]
      \fill[fill=gray!10] (402.798, 642.2304) .. controls (408.8067, 646.0611) and (414.8163, 648.5197) .. (420.827, 649.6064) -- (398.227, 647.1274) -- cycle;
      \fill[fill=gray!10] (528.046, 598.9717) .. controls (551.3833, 605.0184) and (571.7917, 613.2764) .. (589.271, 623.7457) arc[start angle=21.106, end angle=114.8333, x radius=45.2548, y radius=-45.2548] -- cycle;
      \draw[thick] (395.941, 650.2604) .. controls (411.941, 626.2604) and (443.941, 626.2604) .. (459.941, 650.2604);
      \draw[thick] (402.798, 642.2304) .. controls (419.56, 652.9171) and (436.3323, 652.9264) .. (453.115, 642.2584);
      \draw[thick] (380, 606.6667) .. controls (402.6667, 594.6667) and (445.3333, 589.3333) .. (470.9047, 597.7485) .. controls (496.476, 606.1637) and (504.952, 628.3273) .. (500.952, 645.6607) .. controls (496.952, 662.994) and (480.476, 675.497) .. (458.9047, 681.7485) .. controls (437.3333, 688) and (410.6667, 688) .. (392, 682.6667) .. controls (373.3333, 677.3333) and (362.6667, 666.6667) .. (358.6667, 652) .. controls (354.6667, 637.3333) and (357.3333, 618.6667) .. cycle;
      \node at (376, 656) {\tiny$\bullet$};
      \node at (352, 685) {\small$p_1$};
      \draw[thick] (547.052, 640.0417) circle[radius=45.2548];
      \node at (563.052, 672.0417) {\tiny$\bullet$};
      \node at (563.052, 608.0417) {\tiny$\bullet$};
      \node at (579.052, 700) {\small$p_2$};
      \node at (579.052, 575) {\small$p_3$};
    \end{tikzpicture}
  }} \,.
\end{equation}
To be pedantic, $\sigma$ should depend on $(g,n)$ and on the splitting of genus and marked points.

Notice how the above terms corresponds to the terms appearing in the topological recursion formula (see V.~Bouchard's lecture notes \cite{Bou26}). This is not a coincidence, as we will see in \cref{sec:CohFT}.

We conclude this section with one more natural map between moduli spaces: the \emph{forgetful map}. This is the map that forgets the last marked point:
\begin{equation}
  \pi \colon \Mbar_{g,n+1} \longrightarrow \Mbar_{g,n} \,,
  \qquad
  (\Sigma,p_1,\dots,p_n,p_{n+1}) \longmapsto (\Sigma,p_1,\dots,p_n)^{\textup{stab}} \,.
\end{equation}
Again, to be pedantic, $\pi$ should depend on $(g,n)$. We omit the dependence for a lighter notation.
The suffix `stab' stands for `stabilisation'. Indeed, it may happen that, when forgetting a marked point, the resulting Riemann surface is not stable. This is the case of a marked point $p_{n+1}$ on a rational component with only three special points. The stabilisation process simply contracts this component to a point. If the resulting Riemann surface is still not stable, we keep contracting unstable components until we find a stable result. For example:
\begin{equation}
  \vcenter{\hbox{
    \begin{tikzpicture}[x=1pt,y=1pt,scale=.3]
      \fill[fill=gray!10] (402.798, 642.2304) .. controls (408.8067, 646.0611) and (414.8163, 648.5197) .. (420.827, 649.6064) -- (398.227, 647.1274) -- cycle;
      \fill[fill=gray!10] (528.046, 598.9717) .. controls (551.3833, 605.0184) and (571.7917, 613.2764) .. (589.271, 623.7457) arc[start angle=21.106, end angle=114.8333, x radius=45.2548, y radius=-45.2548] -- cycle;
      \draw[thick] (395.941, 650.2604) .. controls (411.941, 626.2604) and (443.941, 626.2604) .. (459.941, 650.2604);
      \draw[thick] (402.798, 642.2304) .. controls (419.56, 652.9171) and (436.3323, 652.9264) .. (453.115, 642.2584);
      \draw[thick] (380, 606.6667) .. controls (402.6667, 594.6667) and (445.3333, 589.3333) .. (470.9047, 597.7485) .. controls (496.476, 606.1637) and (504.952, 628.3273) .. (500.952, 645.6607) .. controls (496.952, 662.994) and (480.476, 675.497) .. (458.9047, 681.7485) .. controls (437.3333, 688) and (410.6667, 688) .. (392, 682.6667) .. controls (373.3333, 677.3333) and (362.6667, 666.6667) .. (358.6667, 652) .. controls (354.6667, 637.3333) and (357.3333, 618.6667) .. cycle;
      \node at (376, 656) {\tiny$\bullet$};
      \node at (352, 685) {\small$p_1$};
      \draw[thick] (547.052, 640.0417) circle[radius=45.2548];
      \node at (563.052, 672.0417) {\tiny$\bullet$};
      \node at (563.052, 608.0417) {\tiny$\bullet$};
      \node at (579.052, 700) {\small$p_2$};
      \node at (579.052, 575) {\small$p_3$};
    \end{tikzpicture}
  }}
  \longmapsto
  \left(\!\!
    \vcenter{\hbox{
      \begin{tikzpicture}[x=1pt,y=1pt,scale=.3]
        \fill[fill=gray!10] (402.798, 642.2304) .. controls (408.8067, 646.0611) and (414.8163, 648.5197) .. (420.827, 649.6064) -- (398.227, 647.1274) -- cycle;
        \fill[fill=gray!10] (528.046, 598.9717) .. controls (551.3833, 605.0184) and (571.7917, 613.2764) .. (589.271, 623.7457) arc[start angle=21.106, end angle=114.8333, x radius=45.2548, y radius=-45.2548] -- cycle;
        \draw[thick] (395.941, 650.2604) .. controls (411.941, 626.2604) and (443.941, 626.2604) .. (459.941, 650.2604);
        \draw[thick] (402.798, 642.2304) .. controls (419.56, 652.9171) and (436.3323, 652.9264) .. (453.115, 642.2584);
        \draw[thick] (380, 606.6667) .. controls (402.6667, 594.6667) and (445.3333, 589.3333) .. (470.9047, 597.7485) .. controls (496.476, 606.1637) and (504.952, 628.3273) .. (500.952, 645.6607) .. controls (496.952, 662.994) and (480.476, 675.497) .. (458.9047, 681.7485) .. controls (437.3333, 688) and (410.6667, 688) .. (392, 682.6667) .. controls (373.3333, 677.3333) and (362.6667, 666.6667) .. (358.6667, 652) .. controls (354.6667, 637.3333) and (357.3333, 618.6667) .. cycle;
        \node at (376, 656) {\tiny$\bullet$};
        \node at (352, 685) {\small$p_1$};
        \draw[thick] (547.052, 640.0417) circle[radius=45.2548];
        \node at (563.052, 672.0417) {\tiny$\bullet$};
        \node at (579.052, 700) {\small$p_2$};
      \end{tikzpicture}
    }}
  \right)^{\textup{stab}}
  =
  \vcenter{\hbox{
    \begin{tikzpicture}[x=1pt,y=1pt,scale=.3]
      \fill[fill=gray!10] (402.798, 642.2304) .. controls (408.8067, 646.0611) and (414.8163, 648.5197) .. (420.827, 649.6064) -- (398.227, 647.1274) -- cycle;
      \draw[thick] (395.941, 650.2604) .. controls (411.941, 626.2604) and (443.941, 626.2604) .. (459.941, 650.2604);
      \draw[thick] (402.798, 642.2304) .. controls (419.56, 652.9171) and (436.3323, 652.9264) .. (453.115, 642.2584);
      \draw[thick] (380, 606.6667) .. controls (402.6667, 594.6667) and (445.3333, 589.3333) .. (470.9047, 597.7485) .. controls (496.476, 606.1637) and (504.952, 628.3273) .. (500.952, 645.6607) .. controls (496.952, 662.994) and (480.476, 675.497) .. (458.9047, 681.7485) .. controls (437.3333, 688) and (410.6667, 688) .. (392, 682.6667) .. controls (373.3333, 677.3333) and (362.6667, 666.6667) .. (358.6667, 652) .. controls (354.6667, 637.3333) and (357.3333, 618.6667) .. cycle;
      \node at (376, 656) {\tiny$\bullet$};
      \node at (352, 685) {\small$p_1$};
      \node at (490, 640) {\tiny$\bullet$};
      \node at (530, 640) {\small$p_2$};
    \end{tikzpicture}
  }} \,.
\end{equation}

The gluing maps $\rho, \sigma$ and the forgetful map $\pi$ are sometimes referred to as the \emph{tautological maps}. We will see shortly that they play a crucial role in the intersection theory of the moduli space of Riemann surfaces. The main takeaway is that, thanks to the compactification and the introduction of the tautological maps, we can think about the moduli spaces as a collection of spaces connected by maps (rather than ``isolated'' spaces). In particular, we can talk about pullback and pushforward in cohomology.

Let us recall these operations for an arbitrary smooth map
\begin{equation}
  \phi \colon M \longrightarrow N
\end{equation}
between smooth real orbifolds of real dimensions $\dim_{\R}(M) = m$ and $\dim_{\R}(N) = n$.

\paragraph{Pullback.} The pullback is always a well-defined contravariant operation in cohomology corresponding to pre-composition. More precisely, it is a degree-preserving map
\begin{equation}
  \phi^* \colon H^k(N) \longrightarrow H^k(M) \,.
\end{equation}
In terms of differential forms, write locally $\phi$ as $(x_1,\ldots,x_m) \mapsto (y_1(x),\ldots,y_n(x))$ and let $\eta$ be a $k$-form on $N$ locally expressed as $\eta = \eta^{\mu_1,\dots,\mu_k}(y) \, \dd y_{\mu_1} \wedge \cdots \wedge \dd y_{\mu_k}$ (we use Einstein's notation for the summation over repeated indices). Then
\begin{equation}
  \phi^* \eta
  =
  \eta^{\mu_1,\dots,\mu_k}(y(x)) \, \dd y_{\mu_1}(x) \wedge \cdots \wedge \dd y_{\mu_k}(x) \,.
\end{equation}
The pullback is compatible with both addition and cup product.

\paragraph{Pushforward.} The pushforward is well-defined only for maps $\phi$ with compact fibres. In this case, the pushforward defines a covariant operation in cohomology, which corresponds to the geometric idea of ``integration along fibres''. More precisely, if we denote by $r$ the real dimension of the fibres of $\phi$, then
\begin{equation}
  \phi_* \colon H^k(M) \longrightarrow H^{k-r}(N) \,.
\end{equation}
In terms of differential forms, write locally $\phi$ as $(x_1,\ldots,x_r,y_1,\ldots,y_n) \mapsto (y_1,\ldots,y_n)$ and let $\omega$ be a $k$-form on $M$ locally expressed as $\omega = \omega^{\nu_1,\dots,\nu_{k-r}}(x,y) \, \dd x_{1} \wedge \cdots \wedge \dd x_{r} \wedge \dd y_{\nu_1} \wedge \cdots \wedge \dd y_{\nu_{k-r}} + \cdots$. The dots stand for terms with a lower number of $\dd x$'s. Then
\begin{equation}
  (\phi_* \omega)_{q}
  =
  \left( \int_{\phi^{-1}(q)} \omega^{\nu_1,\dots,\nu_{k-r}}(x,q) \, \dd x_{1} \wedge \cdots \wedge \dd x_{r} \right)
  \dd_q y_{\nu_1} \wedge \cdots \wedge \dd_q y_{\nu_{k-r}}
\end{equation}
for all $q \in N$. The pushforward is compatible with the addition, but it does not respect the cup product.

The definition generalises via Poincaré duality whenever both $M$ and $N$ are compact. In this case, the pushforward is simply the pre-composition and post-composition of the pushforward in homology by Poincaré duality:
\begin{equation}
  \phi_* \colon
  H^k(M)
  \overset{\textup{PD}}{\cong}
  H_{m-k}(M)
  \longrightarrow
  H_{m-k}(N)
  \overset{\textup{PD}}{\cong}
  H^{k-(m-n)}(N) \,.
\end{equation}
It coincides with the ``integration along fibres'' whenever $\phi$ has compact fibres of dimension $r = m-n$. 

\paragraph{Projection formula.}
In the case of compact fibres, there is a useful formula, known as projection formula, which expresses integrals over $M$ as integrals over $N$. More precisely: if $\omega \in H^k(M)$ and $\eta \in H^{m-k}(N)$, then
\begin{equation}
  \int_M \omega \wedge \phi^*\eta
  =
  \int_N \phi_* \omega \wedge \eta \,.
\end{equation}

\subsection{Intersection theory and Witten's conjecture}
Recall our main goal: to define and compute integrals over the moduli space of Riemann surfaces. Since $\Mbar_{g,n}$ is a compact orbifold, we can finally discuss integrals of top cohomology classes. However, we do not yet have any natural classes to integrate. There are two natural sources of cohomology classes.
\begin{itemize}
  \item The Poincaré dual of natural (complex) subspaces.

  \item Chern classes of natural complex vector bundles.
\end{itemize}
In both cases, cohomology classes of even degree are produced. For this reason, when multiplying classes in cohomology, we shall always omit the cap product since the cap product of even-degree cohomology classes is commutative.

We have already encountered several subspaces of $\Mbar_{g,n}$, namely the boundary strata. Recall that for a fixed stable graph $\Gamma$ of type $(g,n)$, the associated subspace $\Mbar_{\Gamma}$ has complex dimension $\dim(\Mbar_{g,n}) - |E(\Gamma)|$. We deduce that the Poincaré dual, denoted by brackets $[\,\cdot\,]$, lives in
\begin{equation}
  [\Gamma] \in H^{2|E(\Gamma)|}(\Mbar_{g,n},\Q) \,.
\end{equation}
It can be expressed as a pushforward along the gluing maps:
\begin{equation}\label{eq:Gamma:PD:xi}
  [\Gamma] = \frac{1}{|\Aut(\Gamma)|} \xi_{\Gamma,*} \bm{1} \,.
\end{equation}
The element $\bm{1}$ on the right-hand side is the unit in $H^{\bullet}(\Mbar_{\Gamma},\Q)$. In particular, the Poincaré dual of the entire space, corresponding to the stable graph with a single vertex of genus $g$, no edges, and $n$ leaves, is the unit in cohomology:
\begin{equation}
  \left[
  \begin{tikzpicture}[baseline]
    \draw[thick] (0,.1) -- ($(0,.1) + (150:.75)$);
    \node at ($(0,.1) + (150:.75)$) [left] {\footnotesize$1$};
    \node[rotate=90] at ($(0,.1) + (180:.55)$) {\tiny$\cdots$};
    \draw[thick] (0,.1) -- ($(0,.1) + (-150:.75)$);
    \node at ($(0,.1) + (-150:.75)$) [left] {\footnotesize$n$};
    \draw[thick,fill=white] (0,.1) circle[radius=.25cm];
    \node at (0,.1) {\footnotesize$g$};
  \end{tikzpicture}
  \right]
  =
  \bm{1} \in H^{0}(\Mbar_{g,n},\Q) \,.
\end{equation}

Let us discuss now Chern classes of complex vector bundles. A complex vector bundle over $\Mbar_{g,n}$ is the assignment of a complex vector space to each isomorphism class of stable Riemann surfaces so that, as the stable Riemann surface varies within the moduli space, the assigned vector spaces vary smoothly and are coherently glued together. Once a complex vector bundle $\mc{V} \to \Mbar_{g,n}$ is given, we can consider its Chern classes:
\begin{equation}
  \cc_k(\mc{V}) \in H^{2k}(\Mbar_{g,n},\Q) \,,
  \qquad
  k = 0,1,\dots,\rk(\mc{V}) \,,
\end{equation}
where $\rk(\mc{V})$ denotes the complex rank of $\mc{V}$, that is, the complex dimension of the fibres. The zeroth Chern class is always the unit in cohomology: $\cc_0(\mc{V}) = \bm{1}$. Chern classes are topological invariants associated with complex vector bundles and offer a simple test to determine whether two vector bundles are not isomorphic: if the Chern classes of a pair of vector bundles differ, then the vector bundles are distinct (the converse, however, is not necessarily true). Geometrically, they provide information about the number of linearly independent sections a vector bundle has and may be expressed as coefficients of the characteristic polynomial of the curvature form of any Hermitian connection $\nabla$ on $\mc{V}$ (the cohomology class does not depend on the choice of connection):
\begin{equation}
  \cc(\mc{V};t) = \sum_{k=0}^{\rk(\mc{V})} \cc_k(\mc{V}) \, t^k
  =
  \det\left( \Id - t \, \frac{F_{\nabla}}{2\pi\iu} \right) .
\end{equation} 

The first example of such a holomorphic vector bundle is the so-called \emph{$i$-th cotangent line bundle}: for each $i \in \set{1,\dots,n}$, set
\begin{equation}
  \mc{L}_i \longrightarrow \Mbar_{g,n} \,,
  \qquad\qquad
  \mc{L}_i|_{(\Sigma,p_1,\dots,p_n)}
  =
  T_{p_i}^*\Sigma \,.
\end{equation}
In other words, the fibre over $(\Sigma,p_1,\dots,p_n)$ is the holomorphic cotangent space at the $i$-th marked point. Since $T_{p_i}^*\Sigma$ is a complex vector space of dimension $1$, the associated bundle $\mc{L}_i$ has complex rank $1$: it is a line bundle. We then consider its first Chern class:
\begin{equation}
  \psi_i
  =
  \cc_1(\mc{L}_i) \in H^2(\Mbar_{g,n},\Q) \,.
\end{equation}
These are called cotangent line classes, or simply \emph{$\psi$-classes}. As usual, strictly speaking, $\psi$-classes should depend on $(g,n)$. We omit this dependence, that is hopefully clear from the context. As we shall see shortly, $\psi$-classes appear in the seminal work of Witten on topological 2D gravity \cite{Wit91} and represent a cornerstone of all physical theories connected to the moduli space of Riemann surfaces, such as JT gravity and topological string theory.

\begin{figure}
  \centering
  \begin{tikzpicture}[x=1pt,y=1pt,scale=.75]
    \filldraw[thick,fill=gray!10] (152.806, 544.676) .. controls (178.9353, 576.2253) and (202.6667, 618.6667) .. (224, 672) .. controls (293.3333, 666.6667) and (357.3333, 666.6667) .. (416, 672) .. controls (405.3333, 618.6667) and (400.2337, 576.4947) .. (400.701, 545.484) .. controls (325.567, 549.828) and (242.9353, 549.5587) .. (152.806, 544.676) -- cycle;

    \filldraw[thick,fill=white] (250.3186, 648.2929) .. controls (251.6519, 657.6262) and (243.6519, 668.2929) .. (243.6519, 678.9596) .. controls (243.6519, 689.6262) and (251.6519, 700.2929) .. (248.9852, 710.9596) .. controls (246.3186, 721.6262) and (232.9852, 732.2929) .. (220.9852, 730.9596) .. controls (208.9852, 729.6262) and (198.3186, 716.2929) .. (196.9852, 705.6262) .. controls (195.6519, 694.9596) and (203.6519, 686.9596) .. (204.9852, 674.9596) .. controls (206.3186, 662.9596) and (200.9852, 646.9596) .. (203.6519, 637.6262) .. controls (206.3186, 628.2929) and (216.9852, 625.6262) .. (227.6519, 628.2929) .. controls (238.3186, 630.9596) and (248.9852, 638.9596) .. (250.3186, 648.2929);
    \fill[gray!10] (228.2625, 651.3569) .. controls (233.1104, 657.8444) and (233.371, 663.2094) .. (229.0445, 667.4519) .. controls (220.1755, 662.7269) and (220.1755, 654.7269) .. (228.2625, 651.3569) -- cycle;
    \draw[thick] (224.1755, 646.7269) .. controls (234.8421, 657.3935) and (234.8421, 665.3935) .. (224.1755, 670.7269);
    \draw[thick] (219.6519, 694.9596) .. controls (231.6519, 702.9596) and (231.6519, 714.9596) .. (219.6519, 714.9596);
    \draw[thick] (225.6449, 700.4066) .. controls (219.6519, 702.9596) and (219.6519, 706.9596) .. (224.4539, 714.1186);
    \draw[thick] (228.2625, 651.3569) .. controls (233.1104, 657.8444) and (233.371, 663.2094) .. (229.0445, 667.4519);
    \draw[thick] (228.2625, 651.3569) .. controls (220.1755, 654.7269) and (220.1755, 662.7269) .. (229.0445, 667.4519);
    \filldraw[BurntOrange,fill opacity=.5] (224.607, 717.088) .. controls (229.8037, 711.8573) and (233.201, 705.7487) .. (234.799, 698.762) .. controls (239.6737, 703.7427) and (248.9117, 708.3963) .. (262.513, 712.723) .. controls (253.2297, 721.219) and (246.2657, 730.1567) .. (241.621, 739.536) .. controls (236.3637, 731.0453) and (230.6923, 723.5627) .. (224.607, 717.088) -- cycle;
    \node at (224, 563) {$(\Sigma,p_1,\dots,p_n)$};
    \node at (224, 576) {\small$\bullet$};
    \draw[thick,densely dotted] (224, 620) -- (224, 588);
    \node at (208, 648) {\small$\bullet$};
    \node at (240, 716) {\small$\bullet$};
    \node[BurntOrange] at (235, 753) {$T^*_{p_i}\Sigma$};

    \filldraw[thick,fill=white] (306.0782, 683.6978) .. controls (307.3183, 694.0123) and (302.8857, 704.5067) .. (303.7557, 717.7955) .. controls (304.6257, 731.0843) and (310.7983, 747.1677) .. (308.1317, 757.8343) .. controls (305.465, 768.501) and (293.959, 773.751) .. (282.8727, 769.7093) .. controls (271.7863, 765.6677) and (261.1197, 752.3343) .. (259.7863, 741.6677) .. controls (258.453, 731.001) and (266.453, 723.001) .. (269.037, 713.1213) .. controls (271.621, 703.2417) and (268.789, 691.4823) .. (271.5488, 681.1678) .. controls (274.3087, 670.8533) and (282.6603, 661.9837) .. (290.2927, 662.6162) .. controls (297.925, 663.2487) and (304.838, 673.3833) .. (306.0782, 683.6978);
    \draw[thick] (284, 684) .. controls (292, 688) and (292, 700) .. (280, 700);
    \draw[thick] (276, 736) .. controls (292, 740) and (296, 752) .. (284, 760);
    \draw[thick] (283.105, 738.677) .. controls (276, 744) and (280, 756) .. (288.95, 755.359);
    \draw[thick] (288.91, 689.167) .. controls (284, 688) and (280, 696) .. (284.979, 699.109);
    \filldraw[BurntOrange,fill opacity=.5] (276.31, 760.116) .. controls (285.1593, 758.5193) and (292.9373, 755.339) .. (299.644, 750.575) .. controls (302.0453, 759.2883) and (307.3927, 766.9377) .. (315.686, 773.523) .. controls (304.006, 776.835) and (292.954, 781.3503) .. (282.53, 787.069) .. controls (282.5707, 778.0417) and (280.4973, 769.0573) .. (276.31, 760.116) -- cycle;
    \node at (288, 608) {\small$\bullet$};
    \draw[thick,densely dotted] (288, 656) -- (288, 624);
    \node at (284, 676) {\small$\bullet$};
    \node at (292, 764) {\small$\bullet$};

    \filldraw[thick,fill=white] (366.8848, 698.294) .. controls (368.0287, 709.6923) and (364.4533, 722.0557) .. (366.6983, 730.7542) .. controls (368.9433, 739.4527) and (377.0087, 744.4863) .. (374.8052, 753.9747) .. controls (372.6017, 763.463) and (360.1293, 777.406) .. (348.5598, 777.7108) .. controls (336.9903, 778.0157) and (326.3237, 764.6823) .. (327.1022, 754.7093) .. controls (327.8807, 744.7363) and (340.1043, 738.1237) .. (341.7065, 728.6892) .. controls (343.3087, 719.2547) and (334.2893, 706.9983) .. (334.5705, 696.0847) .. controls (334.8517, 685.171) and (344.4333, 675.6) .. (352.1557, 676.0312) .. controls (359.878, 676.4623) and (365.741, 686.8957) .. (366.8848, 698.294);
    \draw[thick] (344, 748) .. controls (360, 748) and (360, 760) .. (348, 768);
    \draw[thick] (349.4481, 691.4675) .. controls (357.4481, 699.4675) and (353.4481, 707.4675) .. (349.4481, 707.4675);
    \draw[thick] (352.3401, 695.1435) .. controls (345.4481, 695.4675) and (345.4481, 703.4675) .. (352.7391, 705.5015);
    \draw[thick] (349.604, 748.614) .. controls (344, 756) and (344, 764) .. (352.037, 764.738);
    \filldraw[BurntOrange,fill opacity=.5] (314.698, 750.268) .. controls (322.6707, 744.754) and (328.9343, 738.756) .. (333.489, 732.274) .. controls (338.2563, 739.504) and (345.1143, 744.3027) .. (354.063, 746.67) .. controls (345.9437, 754.0913) and (339.912, 761.7963) .. (335.968, 769.785) .. controls (330.6067, 762.6117) and (323.5167, 756.106) .. (314.698, 750.268) -- cycle;
    \node at (352, 616) {\small$\bullet$};
    \draw[thick,densely dotted] (352, 664) -- (352, 632);
    \node at (360, 712) {\small$\bullet$};
    \node at (336, 752) {\small$\bullet$};
  \end{tikzpicture}
  \caption{An illustration of the cotangent line bundle $\mc{L}_i$.}
\end{figure}

From the $\psi$-classes, we can derive new cohomology classes that are projections of forgotten points: the Morita--Miller--Mumford classes, or simply \emph{$\kappa$-classes}, defined as
\begin{equation}
  \kappa_m = \pi_* \bigl( \psi_{n+1}^{m+1} \bigr) \in H^{2m}(\Mbar_{g,n},\Q) \,,
  \qquad
  m = 0, \dots, 3g-3+n \,,
\end{equation}
where \smash{$\pi \colon \Mbar_{g,n+1} \to \Mbar_{g,n}$} is the forgetful map. Since the fibres of $\pi$ are compact and one-dimensional, the pushforward is well-defined in cohomology and decreases the complex cohomological degree by $1$. As we shall see shortly, the class $2\pi^2 \kappa_1$, called the \emph{Weil--Petersson class}, plays a fundamental role in JT theory and hyperbolic geometry.

A third collection of natural cohomology classes consists of those arising from the most natural vector space associated with a Riemann surface: the space of holomorphic differentials. More precisely, define the Hodge bundle
\begin{equation}
  \mc{H} \longrightarrow \Mbar_{g,n} \,,
  \qquad\qquad
  \mc{H}|_{(\Sigma,p_1,\dots,p_n)}
  =
  \Omega(\Sigma) \,.
\end{equation}
Here, $\Omega(\Sigma)$ denotes the space of holomorphic forms on $\Sigma$, which is a complex vector space of dimension $g$. One should be cautious, however, regarding the definition of holomorphic forms on Riemann surfaces with nodes (that is, the definition of $\mc{H}$ on the boundary of the moduli space). In order to understand how holomorphic forms should be defined on nodal Riemann surfaces, consider the example
\begin{equation}
  E_t \colon \quad
  y^2 = x (x-1) (x-t) \,.
\end{equation}
For $t \neq 0$, the space of holomorphic forms on $E_t$ is one dimensional and generated by
\begin{equation}
  \omega_t = \frac{\dd x}{y} = \frac{\dd x}{\sqrt{x(x-1)(x-t)}} \,.
\end{equation}
As $t \to 0$, the torus degenerates into a pinched torus, and the holomorphic form $\omega_t$ limits to
\begin{equation}
  \omega_0 = \frac{\dd x}{x \sqrt{x-1}} \,.
\end{equation}
One can verify in local coordinates that $\omega_0$ is no longer holomorphic, but is instead meromorphic with a simple pole at the node and opposite residues at the two branches of the node. The presence of this simple pole is crucial. Indeed, the pinched torus is a $\P^1$ with two points identified; on $\P^1$, there are no non-trivial holomorphic forms; however, there exists a one-dimensional complex vector space of meromorphic forms with simple poles at the two special points and opposite residues. In other words, the dimension of $\Omega(E_t)$ is preserved even in the limit $t \to 0$.

The definition of $\Omega(\Sigma)$ is thus
\begin{equation}
  \Omega(\Sigma)
  =
  \Set{\substack{
    \displaystyle\text{meromorphic form on $\Sigma$} \\
    \displaystyle\text{with at most simple poles at the nodes, opposite residues} \\
    \displaystyle\text{and holomorphic everywhere else}
  }} \,,
\end{equation}
which has constant dimension $g$ as $\Sigma$ varies in $\Mbar_{g,n}$ (independently of the marked points). We then define the Hodge classes, or simply \emph{$\lambda$-classes}, as the Chern classes of the Hodge bundle:
\begin{equation}
  \lambda_k = \cc_k(\mc{H}) \in H^{2k}(\Mbar_{g,n},\Q) \,,
  \qquad
  k = 0, \dots, g \,,
\end{equation}
As we will briefly mention in \cref{sec:next}, the Hodge class plays a fundamental role in topological string theory.

We conclude this section with a brief overview of Witten's conjecture. We begin with two facts regarding $\psi$-class intersection numbers, also known as \emph{Witten's correlators}: the string and dilaton equations. These equations relate integrals of $\psi$-classes over different moduli spaces. Such integrals are conveniently written following Witten's notation as 
\begin{equation}
  \braket{\tau_{d_1} \cdots \tau_{d_n}}_{g}
  =
  \int_{\Mbar_{g,n}} \psi_1^{d_1} \cdots \psi_n^{d_n} \,,
  \qquad\qquad
  d_i \ge 0 \,.
\end{equation}
The integral is set to be zero unless $d_1 + \cdots + d_n = 3g-3+n$, in which case the integrand is a top-dimensional cohomology class.
\begin{itemize}
  \item
  \textbf{Geometric string equation.}
  The pullback of $\psi$-classes along the forgetful map is
  \begin{equation}\label{eq:geom:string}
    \pi^{*} \psi_i
    =
    \psi_i - D_i \,,
    \qquad\quad
    D_i =
    \left[
      \begin{tikzpicture}[baseline]
        \draw[thick] (0,.1) -- ($(0,.1) + (150:.75)$);
        \node at ($(0,.1) + (150:1)$) {\footnotesize$1$};
        \node[rotate=90] at ($(0,.1) + (180:.55)$) {\tiny$\cdots$};
        \draw[thick] (0,.1) -- ($(0,.1) + (-150:.75)$);
        \node at ($(0,.1) + (-150:1)$) {\footnotesize$n$};
        \node at ($(0,.1) + (180:1.2)$) {\footnotesize$\widehat{i}$};
        \draw[thick] (0,.1) -- (1,.1);
        \draw[thick] (1,.1) -- ($(1,.1) + (30:.75)$);
        \node at ($(1,.1) + (30:1)$) {\footnotesize$i$};
        \draw[thick] (1,.1) -- ($(1,.1) + (-30:.75)$);
        \node at ($(1,.1) + (-30:1)$) {\footnotesize$n+1$};
        \draw[thick,fill=white] (0,.1) circle[radius=.25cm];
        \node at (0,.1) {\footnotesize$g$};
        \draw[thick,fill=white] (1,.1) circle[radius=.25cm];
        \node at (1,.1) {\footnotesize$0$};
      \end{tikzpicture}
    \right],
  \end{equation}
  where as usual, a caret as in $\widehat{i}$ denotes omission. The $\psi$-class on the left-hand side lies in $\Mbar_{g,n}$, while that on the right lies in $\Mbar_{g,n+1}$.

  \item
  \textbf{Geometric dilaton equation.}
  The $0$-th $\kappa$-class on $\Mbar_{g,n}$ is equal to (minus) the Euler characteristic:
  \begin{equation}
    \kappa_0 = (2g-2+n) \, \bm{1} \in H^0(\Mbar_{g,n},\Q) \,.
  \end{equation}
\end{itemize}

\begin{exercise}
  Employ the geometric string and dilaton equations, together with the projection formula and the expression \labelcref{eq:Gamma:PD:xi} for the Poincaré dual of boundary strata, to prove the following equations satisfied by Witten's correlators.
  \begin{itemize}
    \item
    \textbf{String equation.}
    Integrals over $\Mbar_{g,n+1}$ with no $\psi_{n+1}$ are reduced to integrals over $\Mbar_{g,n}$:
    \begin{equation}
      \int_{\Mbar_{g,n+1}} \psi_1^{d_1} \cdots \psi_{n}^{d_n}
      =
      \sum_{i=1}^n 
        \int_{\Mbar_{g,n}} \left( \prod_{j \neq i} \psi_j^{d_j} \right) \psi_i^{d_i - 1} \,.
    \end{equation}
    In Witten's notation, the string equation amounts to the removal of a $\tau_0$:
    \begin{equation}
      \braket{\tau_{d_1} \cdots \tau_{d_n} \tau_0}_g
      =
      \sum_{i=1}^n \braket{\tau_{d_1} \cdots \tau_{d_i - 1} \cdots \tau_{d_n}}_g \,.
    \end{equation}
    {\small \emph{\faLightbulb \ Hints.} Consider the following facts.
    \begin{itemize}
      \item By looking at cohomological degrees, what can you say about the integral \smash{$\int_{\Mbar_{g,n+1}} \pi^*\alpha$} for a cohomology class $\alpha \in H^{2(3g-3+n)}(\Mbar_{g,n},\Q)$?

      \item Let $D_i$ as in \cref{eq:geom:string}. Interpreting it as a Poincaré dual, one can see that $D_i \cdot D_j = 0$ for all $i \neq j$.
    \end{itemize}}

    \item
    \textbf{Dilaton equation.}
    Integrals over $\Mbar_{g,n+1}$ with a single power of $\psi_{n+1}$ are reduced to integrals over $\Mbar_{g,n}$:
    \begin{equation}
      \int_{\Mbar_{g,n+1}} \psi_1^{d_1} \cdots \psi_{n}^{d_n} \psi_{n+1}
      =
      (2g-2+n)
      \int_{\Mbar_{g,n}} \psi_1^{d_1} \cdots \psi_{n}^{d_n} \,.
    \end{equation}
    In Witten's notation, the dilaton equation amounts to the removal of a $\tau_1$:
    \begin{equation}
      \braket{\tau_{d_1} \cdots \tau_{d_n} \tau_1}_g
      =
      (2g-2+n) \braket{\tau_{d_1} \cdots \tau_{d_n}}_g \,.
    \end{equation}    
  \end{itemize}
\end{exercise}

The string and dilaton equations allow for the computation of all Witten's correlators in genus $0$ and $1$.

\begin{exercise}
  Knowing the string equation and the integral $\braket{\tau_0^3}_0 = 1$, show that all genus $0$, $\psi$-class intersection numbers are determined. Can you prove the following closed formula:
  \begin{equation}
    \braket{\tau_{d_1} \cdots \tau_{d_n}}_0
    =
    \binom{n-3}{d_1, \ldots, d_n} \,,
  \end{equation}
  where $\binom{D}{d_1, \ldots, d_n} = \frac{D!}{d_1! \cdots d_n!}$ is the multinomial coefficient?
\end{exercise}

\begin{exercise}
  Knowing the string equation, the dilaton equation, and the integral $\braket{\tau_1}_1 = \frac{1}{24}$, show that all genus $1$, $\psi$-class intersection numbers are determined. Can you prove the following closed formula:
  \begin{equation}
    \braket{\tau_{d_1} \cdots \tau_{d_n}}_1
    =
    \frac{1}{24}
    \left(
      \binom{n}{d_1, \ldots, d_n}
      -
      \sum_{\epsilon_1,\ldots,\epsilon_n \in \set{0,1}}
        \binom{n - |\epsilon|}{d_1-\epsilon_1, \ldots, d_n-\epsilon_n} (|\epsilon|-2)!
    \right) ,
  \end{equation}
  where $|\epsilon| = \epsilon_1 + \cdots + \epsilon_n$?
\end{exercise}

While the genus $0$ initial value $\braket{\tau_0^3}_0 = 1$ is trivially satisfied, the genus $1$ case $\braket{\tau_1}_1 = \frac{1}{24}$ is rather non-trivial. This can be computed using the geometry of the moduli space $\Mbar_{1,1}$ and its connection to modular forms.

\begin{exercise}
  Prove that $\braket{\tau_1}_1 = \frac{1}{24}$ using the following facts.
  \begin{enumerate}
    \item
    The following identity holds for arbitrary line bundle $\mc{L}$: $\cc_1(\mc{L}) = \frac{1}{k} \cc_1(\mc{L}^{\otimes k})$.

    \item
    For an arbitrary line bundle $\mc{L}$, we have $\cc_1(\mc{L}) = [Z - P]$, where $Z$ and $P$ are the divisors of zeros and poles of a generic meromorphic section of $\mc{L}$ and $[\;\cdot\;]$ denotes the Poincaré dual\footnote{
      Poincaré duality for orbifolds involves the automorphism group. More precisely, if $Z$ is a sub-orbifold of $X$ with underlying topological space $\hat{Z}$, then $[Z] = \frac{1}{|G|}[\hat{Z}]$, where $G$ is the automorphism group of a generic point in $\hat{Z}$.
    }.

    \item
    Consider the cotangent line bundle $\mc{L}_1^{\otimes k} \to \Mbar_{1,1}$. There is a canonical identification of the vector space of holomorphic sections of $\mc{L}_1^{\otimes k}$ and the vector space of weight $k$ modular forms.

    \item
    The following (combination of) Eisenstein series
    \begin{equation}
    \begin{aligned}
      G_4(\tau) &= \sum_{\lambda \in (\Z + \tau \Z)\setminus \set{0}} \frac{1}{\lambda^4} \,, \\
      G_6(\tau) &= \sum_{\lambda \in (\Z + \tau \Z)\setminus \set{0}} \frac{1}{\lambda^6} \,, \\
      \tilde{G}_{12}(\tau)
      &=\left( \frac{G_4(\tau)}{2\zeta(4)} \right)^3
      -
      \left( \frac{G_6(\tau)}{2\zeta(6)} \right)^2,
    \end{aligned}
    \end{equation}
    are modular forms of weight $4$, $6$, and $12$ respectively. Furthermore, they have a unique simple zero at $\tau = \frac{1+\iu\sqrt{3}}{2}$, $\tau = \iu$, and $\tau = +\iu \infty$ respectively.
  \end{enumerate}
\end{exercise}

We can now state Witten's conjecture. To start with, let us package Witten's correlators in a single generating series: let $t_d$ (for $d \ge 0$) be a set of formal variables and set
\begin{equation}
  Z(t_0,t_1,t_2,\dots;\hbar)
  =
  \exp\left(
    \sum_{\substack{g \ge 0, \, n \ge 1 \\ 2g-2+n > 0}}
    \frac{
      \hbar^{2g-2+n}
    }{n!}
    \sum_{d_1,\dots,d_n \ge 0}
      \braket{\tau_{d_1}\cdots\tau_{d_n}}_g \,
      \prod_{i=1}^n t_{d_i}
  \right) .
\end{equation}
The generating series $Z$ arises as a partition function in topological 2D quantum gravity. The string and dilaton equations may be written as differential operators annihilating $Z$ in the following way.

\begin{exercise}
  Define the differential operators
  \begin{align}
    \label{eq:L-1}
    & L_{-1}
    =
    \hbar
    \dde{}{t_0}
    - \hbar^2
    \left(
      \sum_{k \ge 1} t_{k} \, \dde{}{t_{k-1}}
      +
      \frac{t_0^2}{2}
    \right)
    \,, \\
    \label{eq:L0}
    & L_{0}
    =
    \hbar
    \dde{}{t_1}
    - \hbar^2
    \left(
      \sum_{k \ge 0} \frac{2k+1}{3} \, t_{k} \, \dde{}{t_{k}}
      +
      \frac{1}{24}
    \right)
    \,.
  \end{align}
  Prove the following:
  \begin{itemize}
    \item
    The string equation and $\braket{\tau_0^3}_0 = 1$ are equivalent to the equation $L_{-1}~Z = 0$.

    \item
    The dilaton equation and $\braket{\tau_1}_1 = \frac{1}{24}$ are equivalent to the equation $L_{0}~Z = 0$.
  \end{itemize}
\end{exercise}

The operators $L_{-1}$ and $L_0$ may be viewed as the beginning of (a representation of a subalgebra of) the Virasoro algebra. More precisely, consider the Lie algebra $\Vir$ of holomorphic differential operators spanned by
\begin{equation}
  \msc{L}_n =  - \, z^{n+1} \dde{}{z} \,,
  \qquad\qquad
  n \ge -1 \,.
\end{equation}
The bracket is given by $[\msc{L}_m, \msc{L}_n] = (m - n)\msc{L}_{m+n}$.

The collection $(L_{-1}, L_0)$ of differential operators can be uniquely extended (under a certain homogeneity restriction) to a complete representation of (an $\hbar$-deformation of) $\Vir$. For $n \ge 1$, these are given by
\begin{multline}\label{eq:Ln}
  L_n
  =
  \hbar
  \dde{}{t_{n+1}}
  -
  \hbar^2
  \left(
    \vphantom{\sum_{\substack{a,b \ge 0 \\ a+b = n-1}}}
    \sum_{k \ge 0}
      \frac{(2n+2k+1)!!}{(2n+3)!! (2k-1)!!} \,
      t_{k} \, \dde{}{t_{k+n}}
  \right.\\
  \left.
    +
    \frac{1}{2}
    \sum_{\substack{a,b \ge 0 \\ a+b = n-1}}
      \frac{(2a+1)!! (2b+1)!!}{(2n+3)!!} \,
      \frac{\de^2}{\de t_{a} \de t_{b}}
  \right) .
\end{multline}
Here $m!!$ denotes the double factorial, defined recursively as $m!! = m \cdot (m-2)!!$ with initial conditions $0!! = 1!! = 1$. These are precisely the differential constraints appearing in Bouchard's course \cite{Bou26}!!

\begin{exercise}
  Prove that the collection $(\ms{L}_{n} \coloneqq - \frac{(2n+3)!!}{2} L_n )_{n \ge -1}$ of differential operators defined by \cref{eq:L-1,eq:L0,eq:Ln} is indeed a representation of $\Vir$:
  \begin{equation}
    [\ms{L}_m, \ms{L}_n] = \hbar^2 (m - n) \ms{L}_{m+n} \,.
  \end{equation}
  This, together with the form \eqref{eq:Ln} of the operators, proves that $(\ms{L}_{n})_{n \ge -1}$ form an Airy ideal \cite{Bou26}.
\end{exercise}

\begin{theorem}[{Witten's conjecture/Kontsevich's theorem}]
  The differential operators $(L_{n})_{n \ge -1}$ annihilate the partition function $Z$:
  \begin{equation}\label{eq:Vir}
    L_{n}~Z = 0
    \qquad\qquad
    \forall n \ge -1 \,.
   \end{equation}
   Moreover, the above system of equations (known as \emph{Virasoro constraints}) uniquely determine all intersection numbers.
\end{theorem}

We remark that Witten's original formulation of his conjecture states that $Z$ is the unique tau-function of the Korteweg--de~Vries (KdV) hierarchy satisfying the string equation $L_{-1}~Z = 0$. The KdV hierarchy is an infinite sequence of partial differential equations which extends in a certain sense the KdV equation. The equivalent statement in terms of Virasoro constraints was proved by R.~Dijkgraaf, H.~Verlinde, E.~Verlinde \cite{DVV91}.

\begin{exercise}
  Show that the Virasoro constraints are equivalent to the following \emph{topological recursion} for Witten's correlators:
  \begin{multline}\label{eq:TR:psi:classes}
    \braket{\tau_{d_1} \cdots \tau_{d_n}}_g
    =
    \sum_{m=2}^n
      \frac{(2d_1 + 2d_m - 1)!!}{(2d_1 + 1)!! \, (2d_m - 1)!!}
      \braket{\tau_{d_1 + d_m - 1} \tau_{d_2} \cdots \widehat{\tau_{d_m}} \cdots \tau_{d_n}}_{g}
    \\
    +
    \frac{1}{2} \sum_{a + b = d_1 - 2}
      \frac{(2a+1)!! \, (2b+1)!!}{(2d_1 + 1)!!}
      \biggl(
        \braket{\tau_a \tau_b \tau_{d_2} \cdots \tau_{d_n}}_{g-1} \\
        +
        \sum_{ \substack{g_1 + g_2 = g \\ I_1 \sqcup I_2 = \set{d_2,\ldots,d_n} }}
          \braket{\tau_a \tau_{I_1}}_{g_1} \braket{\tau_b \tau_{I_2}}_{g_2}
      \biggr) \,.
  \end{multline}
  Prove that the above recursion is equivalent to the Eynard--Orantin topological recursion formula \cite{EO07} (see \cite{Bou26}) on the Airy spectral curve $(\P^1,\; x(z) = \frac{z^2}{2},\; y(z) = z,\; \omega_{0,2}(z_1,z_2) = \frac{\dd z_1 \dd z_2}{(z_1 - z_2)^2})$:
  \begin{equation}\label{eq:omegagn:Airy}
    \omega_{g,n}(z_1,\dots,z_n)
    =
    (-1)^n
    \sum_{\substack{d_1,\dots,d_n \ge 0 \\ d_1 + \cdots + d_n = 3g-3+n}}
      \braket{\tau_{d_1} \cdots \tau_{d_n}}_g
      \prod_{i=1}^n \frac{(2d_i+1)!!}{z_i^{2d_i+2}} \dd z_i \,.
  \end{equation}
\end{exercise}

\begin{table}
  \centering
  {\small
  \renewcommand{\arraystretch}{1.3}
  \hfill
  \begin{tabular}[t]{c|c|c}
    \toprule
    $(g,n)$ & $\braket{\tau_{d_1} \!\cdots \tau_{d_n}}_g$ & $\ast$ \\
    \midrule
    $(0,3)$ & $\braket{\tau_0^3}_0$ & $1$ \\
    \midrule
    $(0,4)$ & $\braket{\tau_0^3 \tau_1}_0$ & $1$ \\
    \midrule
    \multirow{2}{*}{$(0,5)$}
    & $\braket{\tau_0^4 \tau_2}_0$ & $1$ \\
    & $\braket{\tau_0^3 \tau_1^2}_0$ & $2$ \\
    \midrule
    \multirow{3}{*}{$(0,6)$}
    & $\braket{\tau_0^5 \tau_3}_0$ & $1$ \\
    & $\braket{\tau_0^4 \tau_1 \tau_2}_0$ & $3$ \\
    & $\braket{\tau_0^3 \tau_1^3}_0$ & $6$ \\
    \midrule
    \multirow{4}{*}{$(0,7)$}
    & $\braket{\tau_0^6 \tau_4}_0$ & $1$ \\
    & $\braket{\tau_0^5 \tau_1 \tau_3}_0$ & $4$ \\
    & $\braket{\tau_0^5 \tau_2^2}_0$ & $6$ \\
    & $\braket{\tau_0^4 \tau_1^2 \tau_2}_0$ & $12$ \\
    & $\braket{\tau_0^3 \tau_1^4}_0$ & $24$ \\
    \bottomrule
  \end{tabular}
  \hfill
  \begin{tabular}[t]{c|c|c}
    \toprule
    $(g,n)$ & $\braket{\tau_{d_1} \!\cdots \tau_{d_n}}_g$ & $\ast$ \\
    \midrule
    $(1,1)$ & $\braket{\tau_1}_1$ & $\tfrac{1}{24}$ \\
    \midrule
    \multirow{2}{*}{$(1,2)$}
    & $\braket{\tau_0 \tau_2}_1$ & $\tfrac{1}{24}$ \\
    & $\braket{\tau_1^2}_1$ & $\tfrac{1}{24}$ \\
    \midrule
    \multirow{3}{*}{$(1,3)$}
    & $\braket{\tau_0^2 \tau_3}_1$ & $\tfrac{1}{24}$ \\
    & $\braket{\tau_0 \tau_1 \tau_2}_1$ & $\tfrac{1}{12}$ \\
    & $\braket{\tau_1^3}_1$ & $\tfrac{1}{12}$ \\
    \midrule
    \multirow{4}{*}{$(1,4)$}
    & $\braket{\tau_0^3 \tau_4}_1$ & $\tfrac{1}{24}$ \\
    & $\braket{\tau_0^2 \tau_1 \tau_3}_1$ & $\tfrac{1}{8}$ \\
    & $\braket{\tau_0^2 \tau_2^2}_1$ & $\tfrac{1}{6}$ \\
    & $\braket{\tau_0 \tau_1^2 \tau_2}_1$ & $\tfrac{1}{4}$ \\
    & $\braket{\tau_1^4}_1$ & $\tfrac{1}{4}$ \\
    \bottomrule
  \end{tabular}
  \hfill
  \begin{tabular}[t]{c|c|c}
    \toprule
    $(g,n)$ & $\braket{\tau_{d_1} \!\cdots \tau_{d_n}}_g$ & $\ast$ \\
    \midrule
    $(2,1)$ & $\braket{\tau_4}_2$ & $\tfrac{1}{1152}$ \\
    \midrule
    \multirow{3}{*}{$(2,2)$}
    & $\braket{\tau_0 \tau_5}_2$ & $\tfrac{1}{1152}$ \\
    & $\braket{\tau_1 \tau_4}_2$ & $\tfrac{1}{384}$ \\
    & $\braket{\tau_2 \tau_3}_2$ & $\tfrac{29}{5760}$ \\
    \midrule
    $(3,1)$ & $\braket{\tau_7}_3$ & $\tfrac{1}{82944}$ \\
    \midrule
    \multirow{5}{*}{$(3,2)$}
    & $\braket{\tau_0 \tau_8}_3$ & $\tfrac{1}{82944}$ \\
    & $\braket{\tau_1 \tau_7}_3$ & $\tfrac{5}{82944}$ \\
    & $\braket{\tau_2 \tau_6}_3$ & $\tfrac{77}{414720}$ \\
    & $\braket{\tau_3 \tau_5}_3$ & $\tfrac{503}{1451520}$ \\
    & $\braket{\tau_4^2}_3$ & $\tfrac{607}{1451520}$ \\
    \midrule
    $(4,1)$ & $\braket{\tau_{10}}_4$ & $\tfrac{1}{7962624}$ \\
    \bottomrule
  \end{tabular}
  \hfill
  }
  \caption{Some $\psi$-classes intersection numbers, computed using the topological recursion relation \labelcref{eq:TR:psi:classes}.}
  \label{table:psi:classes}
\end{table}

As mentioned in the introduction, Witten's motivation for the above conjecture finds its roots in 2D quantum gravity. In the classical setting, the spacetime is a surface while the gravitational field is a Riemannian metric on the surface itself. In an attempt to quantise such a theory, one should compute a certain path integral over the space of all possible Riemannian metrics on all possible surfaces. The space of Riemannian metrics over a fixed topological surface is infinite-dimensional, and there are two possible ways to give meaning to such an ill-defined quantity.
\begin{itemize}
  \item
  The first way is to approximate the Riemann surface by small triangles. Thus, the integral over all metrics is replaced by a sum over triangulations. This combinatorial problem can be solved, and the Virasoro constraints appeared in works devoted to the enumeration of triangulations on surfaces, which admit a matrix model formulation.

  \item
  Alternatively, one can compute the partition function by integrating first over all conformally equivalent metrics. Afterward, the remaining integral is performed over the moduli space of Riemann surfaces, and, more precisely, one has to compute integrals of the form $\braket{\tau_{d_1} \cdots \tau_{d_n}}_{g}$. 
\end{itemize}
Witten's conjecture states that the partition functions resulting from the two approaches coincide, based on the physical expectation that there is a unique theory of gravity.

Kontsevich's proof follows the matrix model/discretisation idea (see \cite{Zvo} for a rigorous proof). He started by considering the moduli space of metric ribbon graphs of genus $g$ with $n$ faces of fixed length $L_1,\ldots,L_n$, which comes with a natural (symplectic) volume form. Interpreting metric ribbon graphs as a discretisation of Riemannian metrics, Kontsevich expressed these volumes precisely as the $\psi$-class intersection numbers
\begin{equation}
\begin{split}
  V_{g,n}(L_1,\ldots,L_n)
  &=
  \int_{\Mbar_{g,n}} \exp\left( \frac{1}{2} \sum_{i=1}^n L_i^2 \psi_i \right) \\
  &=
  \sum_{\substack{d_1,\dots,d_n \ge 0 \\ d_1 + \cdots + d_n = 3g-3+n}}
    \braket{\tau_{d_1} \cdots \tau_{d_n}}_g
    \prod_{i=1}^n \frac{L_i^{2d_i}}{2^{d_i} d_i!} \,.
\end{split}
\end{equation}
Note that $V_{g,n}(L_1,\ldots,L_n)$ is a symmetric polynomial in the boundary lengths squared. The Laplace transform of such a volume is computed as the rational function
\begin{equation}\label{eq:vol:tau}
\begin{split}
  \widehat{V}_{g,n}(\lambda_1,\ldots,\lambda_n)
  &=
  \left( \prod_{i=1}^n \int_{0}^{\infty} \dd L_i \, e^{-\lambda_i L_i} \right)
    V_{g,n}(L_1,\ldots,L_n) \\
  &=
  \sum_{\substack{d_1,\dots,d_n \ge 0 \\ d_1 + \cdots + d_n = 3g-3+n}}
    \braket{\tau_{d_1} \cdots \tau_{d_n}}_g
    \prod_{i=1}^n \frac{(2d_i-1)!!}{\lambda_i^{2d_i+1}} \,.
\end{split}
\end{equation}
Notice that $(\dd_{\lambda_1} \cdots \dd_{\lambda_n}) \widehat{V}_{g,n}(\lambda) = \omega_{g,n}(\lambda)$ is precisely the topological recursion correlator from \labelcref{eq:omegagn:Airy} computed from the Airy spectral curve.

As $V_{g,n}(L_1,\ldots,L_n)$ is the volume of the moduli space of metric ribbon graphs of genus $g$ with $n$ faces of fixed length $L_1,\ldots,L_n$, he obtained an expression for the Laplace transform as a sum over ribbon graphs:
\begin{equation}\label{eq:vol:ribbon}
  \widehat{V}_{g,n}(\lambda_1,\ldots,\lambda_n)
  =
  2^{2g-2+n}
  \sum_{\mb{G}}
    \frac{1}{|\Aut(\mb{G})|}
    \prod_{e = (i,j) \in E(\mb{G})} \frac{1}{\lambda_i + \lambda_j} \,,
\end{equation}
where the sum is over all trivalent ribbon graphs of genus $g$ with $n$ faces labelled by $1,\dots,n$ and $e = (i,j)$ denotes the two (possibly equal) faces adjacent to the edge $e$ in the graph $\mb{G}$.

\begin{figure}[b]
  \centering
  \begin{tikzpicture}[baseline,scale=.6]
      \draw[line width=10.5pt] (0,1.7) to[out=90,in=180] (.5,2.3) to[out=0,in=90] (1,1.7) to[out=-90,in=90] (0,0) -- (0,-.1);
      \draw[line width=9pt,white] (0,1.7) to[out=90,in=180] (.5,2.3) to[out=0,in=90] (1,1.7) to[out=-90,in=90] (0,0) -- (0,-.2);
      \draw[BrickRed] (0,1.7) to[out=90,in=180] (.5,2.3) to[out=0,in=90] (1,1.7) to[out=-90,in=90] (0,0) -- (0,-.2);

      \draw[line width=10.5pt] (0,-1.7) -- (0,-.1);
      \draw[line width=10.5pt] (0,0) ellipse (2cm and 1.7cm);
      
      \draw[line width=9pt,white] (0,-1.7) -- (0,0);
      \draw[line width=9pt,white] (0,1.7) to[out=90,in=180] (.5,2.3);
      \draw[line width=9pt,white] (0,0) ellipse (2cm and 1.7cm);
      
      \draw[BrickRed] (0,1.7) to[out=90,in=180] (.5,2.3);
      \draw[BrickRed] (0,0) -- (0,-1.7);

      \draw[BrickRed] (0,0) ellipse (2cm and 1.7cm);
      \node[BrickRed] at (0,1.7) {\tiny$\bullet$};
      \node[BrickRed] at (0,-1.7) {\tiny$\bullet$};
  \end{tikzpicture}
  \caption{The unique trivalent ribbon graph of type $(1,1)$.}
  \label{fig:ribbon}
\end{figure}

For example, take $g=n=1$. In this case there is a single trivalent ribbon graph given in \cref{fig:ribbon}, which has automorphism group $\Z_{6}$ (the cyclic permutations of the three edges and an additional $\Z_2$ symmetry swapping the vertices). Then Kontsevich's formula \labelcref{eq:vol:ribbon} gives
\begin{equation}
  \widehat{V}_{1,1}(\lambda_1)
  =
  2 \cdot \frac{1}{6} \cdot \left(\frac{1}{2\lambda_1}\right)^3
  =
  \frac{1}{24} \frac{1}{\lambda_1^3} \,,
\end{equation}
which indeed gives $\braket{\tau_1}_1 = \frac{1}{24}$, following \labelcref{eq:vol:tau}.

\begin{figure}[b]
  \centering
  \begin{tikzpicture}
  \begin{axis}[
      width=.95\textwidth,
      height=.45\textwidth,
      axis y line = left, 
      axis x line = bottom,
      xlabel = {$g$},
      xlabel style={at={(1,0)},above right},
      xmin = 0,
      xmax = 50,
      yticklabel style = {
          /pgf/number format/sci
      },
      ymode=log,
      yminorticks=true,
  ]
  \addplot[thick, color=BrickRed, smooth, domain = 0:50] {1};
  \addplot[NavyBlue, mark = o, mark size=1, only marks] table[col sep=&,row sep=\\,y expr={\thisrow{y}}] {
      x & y \\
      3 & 0.8583222213 \\
      3 & 0.8693203094 \\
      3 & 0.8649996319 \\
      3 & 0.8425321090 \\
      4 & 0.8955247040 \\
      4 & 0.8931843370 \\
      4 & 0.8970849486 \\
      4 & 0.8939273106 \\
      5 & 0.9137995337 \\
      5 & 0.9143160036 \\
      5 & 0.9142034397 \\
      5 & 0.9134257251 \\
      5 & 0.9152403924 \\
      5 & 0.9130044631 \\
      5 & 0.9055513654 \\
      6 & 0.9273763888 \\
      6 & 0.9272593716 \\
      6 & 0.9274231957 \\
      6 & 0.9273507350 \\
      6 & 0.9270213683 \\
      6 & 0.9280094684 \\
      6 & 0.9263720454 \\
      7 & 0.9370201110 \\
      7 & 0.9370470316 \\
      7 & 0.9370429581 \\
      7 & 0.9370070153 \\
      7 & 0.9370717124 \\
      7 & 0.9370284151 \\
      7 & 0.9368657225 \\
      7 & 0.9374622622 \\
      7 & 0.9362194710 \\
      7 & 0.9325408091 \\
      8 & 0.9444475219 \\
      8 & 0.9444412513 \\
      8 & 0.9444493135 \\
      8 & 0.9444468524 \\
      8 & 0.9444333888 \\
      8 & 0.9444630087 \\
      8 & 0.9444364267 \\
      8 & 0.9443470148 \\
      8 & 0.9447344665 \\
      8 & 0.9437618837 \\
      9 & 0.9502930523 \\
      9 & 0.9502945251 \\
      9 & 0.9502943543 \\
      9 & 0.9502925121 \\
      9 & 0.9502954070 \\
      9 & 0.9502940738 \\
      9 & 0.9502882705 \\
      9 & 0.9503033591 \\
      9 & 0.9502863637 \\
      9 & 0.9502332142 \\
      9 & 0.9504989616 \\
      9 & 0.9497183284 \\
      9 & 0.9475325555 \\
      10 & 0.9550281319 \\
      10 & 0.9550277835 \\
      10 & 0.9550282094 \\
      10 & 0.9550281100 \\
      10 & 0.9550274705 \\
      10 & 0.9550286581 \\
      10 & 0.9550279259 \\
      10 & 0.9550251480 \\
      10 & 0.9550334817 \\
      10 & 0.9550221831 \\
      10 & 0.9549886296 \\
      10 & 0.9551787663 \\
      10 & 0.9545389412 \\
      11 & 0.9589386214 \\
      11 & 0.9589387043 \\
      11 & 0.9589386965 \\
      11 & 0.9589385970 \\
      11 & 0.9589387407 \\
      11 & 0.9589386902 \\
      11 & 0.9589384384 \\
      11 & 0.9589389779 \\
      11 & 0.9589385609 \\
      11 & 0.9589371184 \\
      11 & 0.9589420227 \\
      11 & 0.9589342471 \\
      11 & 0.9589120309 \\
      11 & 0.9590527332 \\
      11 & 0.9585190696 \\
      11 & 0.9570722484 \\
      12 & 0.9622235863 \\
      12 & 0.9622235665 \\
      12 & 0.9622235899 \\
      12 & 0.9622235855 \\
      12 & 0.9622235527 \\
      12 & 0.9622236074 \\
      12 & 0.9622235817 \\
      12 & 0.9622234723 \\
      12 & 0.9622237380 \\
      12 & 0.9622234909 \\
      12 & 0.9622226916 \\
      12 & 0.9622257291 \\
      12 & 0.9622202148 \\
      12 & 0.9622049254 \\
      12 & 0.9623119517 \\
      12 & 0.9618602172 \\
      13 & 0.9650218143 \\
      13 & 0.9650218190 \\
      13 & 0.9650218186 \\
      13 & 0.9650218131 \\
      13 & 0.9650218207 \\
      13 & 0.9650218185 \\
      13 & 0.9650218064 \\
      13 & 0.9650218293 \\
      13 & 0.9650218158 \\
      13 & 0.9650217643 \\
      13 & 0.9650219041 \\
      13 & 0.9650217522 \\
      13 & 0.9650212850 \\
      13 & 0.9650232472 \\
      13 & 0.9650192337 \\
      13 & 0.9650083689 \\
      13 & 0.9650916656 \\
      13 & 0.9647044349 \\
      13 & 0.9636765133 \\
      14 & 0.9674340828 \\
      14 & 0.9674340816 \\
      14 & 0.9674340830 \\
      14 & 0.9674340810 \\
      14 & 0.9674340837 \\
      14 & 0.9674340827 \\
      14 & 0.9674340778 \\
      14 & 0.9674340881 \\
      14 & 0.9674340808 \\
      14 & 0.9674340549 \\
      14 & 0.9674341325 \\
      14 & 0.9674340360 \\
      14 & 0.9674337505 \\
      14 & 0.9674350635 \\
      14 & 0.9674320756 \\
      14 & 0.9674241442 \\
      14 & 0.9674902399 \\
      14 & 0.9671546767 \\
      15 & 0.9695350833 \\
      15 & 0.9695350836 \\
      15 & 0.9695350836 \\
      15 & 0.9695350829 \\
      15 & 0.9695350840 \\
      15 & 0.9695350813 \\
      15 & 0.9695350863 \\
      15 & 0.9695350822 \\
      15 & 0.9695350685 \\
      15 & 0.9695351136 \\
      15 & 0.9695350504 \\
      15 & 0.9695348693 \\
      15 & 0.9695357746 \\
      15 & 0.9695335059 \\
      15 & 0.9695275811 \\
      15 & 0.9695809048 \\
      15 & 0.9692873524 \\
      15 & 0.9685196000 \\
      16 & 0.9713814158 \\
      16 & 0.9713814158 \\
      16 & 0.9713814156 \\
      16 & 0.9713814161 \\
      16 & 0.9713814148 \\
      16 & 0.9713814173 \\
      16 & 0.9713814149 \\
      16 & 0.9713814073 \\
      16 & 0.9713814346 \\
      16 & 0.9713813920 \\
      16 & 0.9713812734 \\
      16 & 0.9713819140 \\
      16 & 0.9713801612 \\
      16 & 0.9713756466 \\
      16 & 0.9714192884 \\
      16 & 0.9711603468 \\
      17 & 0.9730167381 \\
      17 & 0.9730167380 \\
      17 & 0.9730167382 \\
      17 & 0.9730167376 \\
      17 & 0.9730167389 \\
      17 & 0.9730167375 \\
      17 & 0.9730167330 \\
      17 & 0.9730167501 \\
      17 & 0.9730167208 \\
      17 & 0.9730166409 \\
      17 & 0.9730171044 \\
      17 & 0.9730157293 \\
      17 & 0.9730122291 \\
      17 & 0.9730483983 \\
      17 & 0.9728183038 \\
      17 & 0.9722231258 \\
      18 & 0.9744752687 \\
      18 & 0.9744752687 \\
      18 & 0.9744752685 \\
      18 & 0.9744752692 \\
      18 & 0.9744752683 \\
      18 & 0.9744752657 \\
      18 & 0.9744752767 \\
      18 & 0.9744752560 \\
      18 & 0.9744752009 \\
      18 & 0.9744755429 \\
      18 & 0.9744744493 \\
      18 & 0.9744716938 \\
      18 & 0.9745020038 \\
      18 & 0.9742962016 \\
      19 & 0.9757842069 \\
      19 & 0.9757842067 \\
      19 & 0.9757842072 \\
      19 & 0.9757842066 \\
      19 & 0.9757842049 \\
      19 & 0.9757842122 \\
      19 & 0.9757841974 \\
      19 & 0.9757841585 \\
      19 & 0.9757844154 \\
      19 & 0.9757835349 \\
      19 & 0.9757813362 \\
      19 & 0.9758069873 \\
      19 & 0.9756218313 \\
      19 & 0.9751469606 \\
      20 & 0.9769654443 \\
      20 & 0.9769654445 \\
      20 & 0.9769654441 \\
      20 & 0.9769654431 \\
      20 & 0.9769654480 \\
      20 & 0.9769654372 \\
      20 & 0.9769654092 \\
      20 & 0.9769656052 \\
      20 & 0.9769648884 \\
      20 & 0.9769631127 \\
      20 & 0.9769850129 \\
      20 & 0.9768175505 \\
      21 & 0.9780367999 \\
      21 & 0.9780367999 \\
      21 & 0.9780367997 \\
      21 & 0.9780367991 \\
      21 & 0.9780368025 \\
      21 & 0.9780367945 \\
      21 & 0.9780367740 \\
      21 & 0.9780369255 \\
      21 & 0.9780363362 \\
      21 & 0.9780348865 \\
      21 & 0.9780537329 \\
      21 & 0.9779015478 \\
      21 & 0.9775138762 \\
      22 & 0.9790129245 \\
      22 & 0.9790129245 \\
      22 & 0.9790129240 \\
      22 & 0.9790129264 \\
      22 & 0.9790129204 \\
      22 & 0.9790129052 \\
      22 & 0.9790130239 \\
      22 & 0.9790125349 \\
      22 & 0.9790113397 \\
      22 & 0.9790276748 \\
      22 & 0.9788887707 \\
      23 & 0.9799059755 \\
      23 & 0.9799059754 \\
      23 & 0.9799059751 \\
      23 & 0.9799059768 \\
      23 & 0.9799059722 \\
      23 & 0.9799059608 \\
      23 & 0.9799060548 \\
      23 & 0.9799056458 \\
      23 & 0.9799046515 \\
      23 & 0.9799189024 \\
      23 & 0.9797916160 \\
      23 & 0.9794691570 \\
      24 & 0.9807261250 \\
      24 & 0.9807261247 \\
      24 & 0.9807261259 \\
      24 & 0.9807261224 \\
      24 & 0.9807261137 \\
      24 & 0.9807261889 \\
      24 & 0.9807258442 \\
      24 & 0.9807250104 \\
      24 & 0.9807375171 \\
      24 & 0.9806204511 \\
      25 & 0.9814819496 \\
      25 & 0.9814819496 \\
      25 & 0.9814819494 \\
      25 & 0.9814819504 \\
      25 & 0.9814819476 \\
      25 & 0.9814819409 \\
      25 & 0.9814820017 \\
      25 & 0.9814817090 \\
      25 & 0.9814810046 \\
      25 & 0.9814920406 \\
      25 & 0.9813840130 \\
      25 & 0.9811115955 \\
      26 & 0.9821807315 \\
      26 & 0.9821807316 \\
      26 & 0.9821807314 \\
      26 & 0.9821807321 \\
      26 & 0.9821807299 \\
      26 & 0.9821807247 \\
      26 & 0.9821807743 \\
      26 & 0.9821805242 \\
      26 & 0.9821799250 \\
      26 & 0.9821897121 \\
      26 & 0.9820897161 \\
      27 & 0.9828286935 \\
      27 & 0.9828286934 \\
      27 & 0.9828286939 \\
      27 & 0.9828286922 \\
      27 & 0.9828286880 \\
      27 & 0.9828287288 \\
      27 & 0.9828285139 \\
      27 & 0.9828280010 \\
      27 & 0.9828367208 \\
      27 & 0.9827438937 \\
      27 & 0.9825107121 \\
      28 & 0.9834311849 \\
      28 & 0.9834311853 \\
      28 & 0.9834311839 \\
      28 & 0.9834311806 \\
      28 & 0.9834312143 \\
      28 & 0.9834310287 \\
      28 & 0.9834305870 \\
      28 & 0.9834383892 \\
      28 & 0.9833519875 \\
      29 & 0.9839928300 \\
      29 & 0.9839928300 \\
      29 & 0.9839928302 \\
      29 & 0.9839928291 \\
      29 & 0.9839928265 \\
      29 & 0.9839928546 \\
      29 & 0.9839926934 \\
      29 & 0.9839923111 \\
      29 & 0.9839993200 \\
      29 & 0.9839186994 \\
      29 & 0.9837168492 \\
      30 & 0.9845176463 \\
      30 & 0.9845176465 \\
      30 & 0.9845176456 \\
      30 & 0.9845176434 \\
      30 & 0.9845176670 \\
      30 & 0.9845175263 \\
      30 & 0.9845171937 \\
      30 & 0.9845235134 \\
      30 & 0.9844481129 \\
      31 & 0.9850091413 \\
      31 & 0.9850091414 \\
      31 & 0.9850091407 \\
      31 & 0.9850091390 \\
      31 & 0.9850091589 \\
      31 & 0.9850090355 \\
      31 & 0.9850087448 \\
      31 & 0.9850144627 \\
      31 & 0.9849437916 \\
      31 & 0.9847673574 \\
      32 & 0.9854703908 \\
      32 & 0.9854703909 \\
      32 & 0.9854703903 \\
      32 & 0.9854703889 \\
      32 & 0.9854704058 \\
      32 & 0.9854702972 \\
      32 & 0.9854700420 \\
      32 & 0.9854752322 \\
      32 & 0.9854088593 \\
      33 & 0.9859041034 \\
      33 & 0.9859041035 \\
      33 & 0.9859041030 \\
      33 & 0.9859041018 \\
      33 & 0.9859041162 \\
      33 & 0.9859040203 \\
      33 & 0.9859037952 \\
      33 & 0.9859085208 \\
      33 & 0.9858460658 \\
      33 & 0.9856905328 \\
      34 & 0.9863126735 \\
      34 & 0.9863126732 \\
      34 & 0.9863126722 \\
      34 & 0.9863126845 \\
      34 & 0.9863125994 \\
      34 & 0.9863124003 \\
      34 & 0.9863167150 \\
      34 & 0.9862578412 \\
      35 & 0.9866982259 \\
      35 & 0.9866982259 \\
      35 & 0.9866982256 \\
      35 & 0.9866982248 \\
      35 & 0.9866982354 \\
      35 & 0.9866981596 \\
      35 & 0.9866979828 \\
      35 & 0.9867019329 \\
      35 & 0.9866463410 \\
      35 & 0.9865082036 \\
      36 & 0.9870626523 \\
      36 & 0.9870626521 \\
      36 & 0.9870626514 \\
      36 & 0.9870626606 \\
      36 & 0.9870625929 \\
      36 & 0.9870624353 \\
      36 & 0.9870660608 \\
      36 & 0.9870134840 \\
      37 & 0.9874076429 \\
      37 & 0.9874076430 \\
      37 & 0.9874076427 \\
      37 & 0.9874076421 \\
      37 & 0.9874076501 \\
      37 & 0.9874075895 \\
      37 & 0.9874074486 \\
      37 & 0.9874107841 \\
      37 & 0.9873609836 \\
      37 & 0.9872374786 \\
      38 & 0.9877347122 \\
      38 & 0.9877347120 \\
      38 & 0.9877347115 \\
      38 & 0.9877347185 \\
      38 & 0.9877346640 \\
      38 & 0.9877345376 \\
      38 & 0.9877376132 \\
      38 & 0.9876903750 \\
      39 & 0.9880452211 \\
      39 & 0.9880452212 \\
      39 & 0.9880452210 \\
      39 & 0.9880452206 \\
      39 & 0.9880452267 \\
      39 & 0.9880451776 \\
      39 & 0.9880450639 \\
      39 & 0.9880479059 \\
      39 & 0.9880030373 \\
      39 & 0.9878919571 \\
      40 & 0.9883403965 \\
      40 & 0.9883403965 \\
      40 & 0.9883403964 \\
      40 & 0.9883403960 \\
      40 & 0.9883404014 \\
      40 & 0.9883403570 \\
      40 & 0.9883402545 \\
      40 & 0.9883428860 \\
      40 & 0.9883002131 \\
      41 & 0.9886213467 \\
      41 & 0.9886213467 \\
      41 & 0.9886213463 \\
      41 & 0.9886213510 \\
      41 & 0.9886213109 \\
      41 & 0.9886212182 \\
      41 & 0.9886236595 \\
      41 & 0.9885830249 \\
      41 & 0.9884825847 \\
      42 & 0.9888890759 \\
      42 & 0.9888890758 \\
      42 & 0.9888890755 \\
      42 & 0.9888890797 \\
      42 & 0.9888890433 \\
      42 & 0.9888889592 \\
      42 & 0.9888912282 \\
      42 & 0.9888524895 \\
      43 & 0.9891444958 \\
      43 & 0.9891444955 \\
      43 & 0.9891444992 \\
      43 & 0.9891444661 \\
      43 & 0.9891443897 \\
      43 & 0.9891465022 \\
      43 & 0.9891095297 \\
      43 & 0.9890182708 \\
      44 & 0.9893884364 \\
      44 & 0.9893884364 \\
      44 & 0.9893884361 \\
      44 & 0.9893884394 \\
      44 & 0.9893884091 \\
      44 & 0.9893883396 \\
      44 & 0.9893903097 \\
      44 & 0.9893549854 \\
      45 & 0.9896216543 \\
      45 & 0.9896216543 \\
      45 & 0.9896216541 \\
      45 & 0.9896216570 \\
      45 & 0.9896216294 \\
      45 & 0.9896215660 \\
      45 & 0.9896234062 \\
      45 & 0.9895896223 \\
      45 & 0.9895063410 \\
      46 & 0.9898448416 \\
      46 & 0.9898448414 \\
      46 & 0.9898448440 \\
      46 & 0.9898448187 \\
      46 & 0.9898447607 \\
      46 & 0.9898464821 \\
      46 & 0.9898141401 \\
      47 & 0.9900586315 \\
      47 & 0.9900586314 \\
      47 & 0.9900586337 \\
      47 & 0.9900586105 \\
      47 & 0.9900585573 \\
      47 & 0.9900601701 \\
      47 & 0.9900291795 \\
      47 & 0.9899528737 \\
      48 & 0.9902636055 \\
      48 & 0.9902636054 \\
      48 & 0.9902636075 \\
      48 & 0.9902635861 \\
      48 & 0.9902635373 \\
      48 & 0.9902650503 \\
      48 & 0.9902353283 \\
      49 & 0.9904602978 \\
      49 & 0.9904602977 \\
      49 & 0.9904602996 \\
      49 & 0.9904602799 \\
      49 & 0.9904602351 \\
      49 & 0.9904616564 \\
      49 & 0.9904331265 \\
      49 & 0.9903629552 \\
      50 & 0.9906492004 \\
      50 & 0.9906492003 \\
      50 & 0.9906492020 \\
      50 & 0.9906491839 \\
      50 & 0.9906491425 \\
      50 & 0.9906504794 \\
      50 & 0.9906230715 \\
  };
  \addplot[ForestGreen, mark = o, mark size=1, only marks] table[col sep=&,row sep=\\,y expr={\thisrow{y}}] {
      x & y \\
      3 & 0.9548697235 \\
      4 & 0.9659214565 \\
      5 & 0.9726292443 \\
      6 & 0.9771321575 \\
      7 & 0.9803634147 \\
      8 & 0.9827948715 \\
      9 & 0.9846906950 \\
      10 & 0.9862102794 \\
      11 & 0.9874554944 \\
      12 & 0.9884944837 \\
      13 & 0.9893745536 \\
      14 & 0.9901295708 \\
      15 & 0.9907844184 \\
      16 & 0.9913577890 \\
      17 & 0.9918639971 \\
      18 & 0.9923141895 \\
      19 & 0.9927171761 \\
      20 & 0.9930800113 \\
      21 & 0.9934084108 \\
      22 & 0.9937070545 \\
      23 & 0.9939798112 \\
      24 & 0.9942299067 \\
      25 & 0.9944600525 \\
      26 & 0.9946725440 \\
      27 & 0.9948693374 \\
      28 & 0.9950521103 \\
      29 & 0.9952223094 \\
      30 & 0.9953811891 \\
      31 & 0.9955298422 \\
      32 & 0.9956692253 \\
      33 & 0.9958001793 \\
      34 & 0.9959234463 \\
      35 & 0.9960396838 \\
      36 & 0.9961494766 \\
      37 & 0.9962533460 \\
      38 & 0.9963517589 \\
      39 & 0.9964451342 \\
      40 & 0.9965338490 \\
      41 & 0.9966182438 \\
      42 & 0.9966986266 \\
      43 & 0.9967752769 \\
      44 & 0.9968484487 \\
      45 & 0.9969183735 \\
      46 & 0.9969852628 \\
      47 & 0.9970493101 \\
      48 & 0.9971106927 \\
      49 & 0.9971695734 \\
      50 & 0.9972261023 \\
  };
  \end{axis}
  \end{tikzpicture}
  \caption{
    The $2$-point correlators, normalised by their leading asymptotics: note the convergence to $1 + \bigO(g^{-1})$. Also observe the differing convergence behaviour of the correlators with a $\tau_0$ insertion (in green) versus without (in blue); this suggests that the subleading terms do depend on the partition $(d_1,\dots,d_n)$. This is indeed the case, and it can be proved via resurgence.
  }
  \label{fig:asymt}
\end{figure}

On the one hand, Kontsevich's theorem gives a sum of graphs, where each graph is weighted by its symmetry factor and by a product of edge weights. This is the typical kind of graph obtained from Wick's theorem, and therefore it can be obtained with a perturbation of a Gaussian Hermitian matrix integral. Specifically, trivalent ribbon graphs are generated by a cubic formal matrix integral, the so-called \emph{Airy matrix integral}:
\begin{equation}
  Z(\Lambda)
  =
  \frac{1}{Z_0(\Lambda)}
  \int
    \dd X \,
    \exp\left(
      N \tr
      \left[
        \frac{X^3}{3} - \Lambda X^2 
      \right]
      \right) ,
  \qquad
  \Lambda = \diag(\lambda_1, \ldots, \lambda_N) \,.
\end{equation}
Here $Z_0(\Lambda) = (\pi/N)^{N^2/2} \prod_{i,j} (\lambda_i + \lambda_j)^{-1/2}$ is a normalisation constant. By Wick's theorem, one can write the large $N$ expansion of $\log{ Z(\Lambda) }$ as a sum over trivalent ribbon graphs:
\begin{equation}
  \log{ Z(\Lambda) }
  =
  \sum_{\substack{g \ge 0, \, n \ge 1 \\ 2g-2+n > 0}}
  \frac{
    N^{-(2g-2+n)}
  }{n!}
  \sum_{\mb{G}}
    \frac{1}{|\Aut(\mb{G})|}
    \prod_{e = (i,j) \in E(\mb{G})} \frac{1}{\lambda_{i} + \lambda_{j}} \,,
\end{equation}
where the sum is over all trivalent ribbon graphs of genus $g$ with $n$ labelled faces.

To conclude, integration by parts (also known as \emph{Schwinger--Dyson equations} in this context) shows that $Z(\Lambda)$ satisfies the Virasoro constraints \labelcref{eq:Vir}, upon identification $\hbar = 2/N$ and the times with the normalised traces of $\Lambda$ which appear naturally in the expansion of the matrix model:
\begin{equation}
  t_d =
  \frac{\tr(\Lambda^{-2d-1})}{(2d-1)!!} \,,
  \qquad\qquad d \ge 0 \,.
\end{equation}

It is worth mentioning that, through resurgence techniques (see \cite{ABS19} or I.~Aniceto's and M.~Mariño's lecture notes \cite{Ani,Mar26}), one can compute the large genus asymptotic Witten's correlators \cite{EGGGL}, see \cref{fig:asymt}:
\begin{equation}\label{eq:asymt}
  \braket{\tau_{d_1} \cdots \tau_{d_n}}_g \, \prod_{i=1}^n (2d_i + 1)!!
  =
  \frac{2^{n-1}}{2\pi} \, \frac{\Gamma(2g-2+n)}{(\frac{2}{3})^{2g-2+n}}
  \left(
    1 + \bigO(g^{-1})
  \right) .
\end{equation}
The subleading asymptotics are also accessible using resurgence. The first proof of this result, using combinatorial and probabilistic arguments, is due to A.~Aggarwal \cite{Agg21}. Notice the Stokes constant $S = \iu$ and the instanton action $A = 2/3$, corresponding to those of the Airy function. This is, naturally, not a coincidence!

\newpage
\section{Cohomological field theories}
\label{sec:CohFT}

The Virasoro constraints satisfied by Witten's correlators provide a recursive method to compute all $\psi$-class intersection numbers. The main geometric property underpinning the constraints is the recursive nature of $\Mbar_{g,n}$. By looking at Witten's correlators as the intersections of the unit with $\psi$-classes, we can rephrase the recursive structure purely in cohomological terms. The unit $\bm{1}_{g,n} \in H^0(\Mbar_{g,n},\Q)$ is stable under pullback by all tautological maps, that is
\begin{align}
  &
  \rho^* \bm{1}_{g,n} = \bm{1}_{g-1,n+2} \,,
  &&
  \rho \colon \Mbar_{g-1,n+2} \to \Mbar_{g,n} \,,\\
  &
  \sigma^* \bm{1}_{g,n} = \bm{1}_{g_1,1+|I_1|} \otimes \bm{1}_{g_2,1+|I_2|} \,,
  &&
  \sigma \colon \Mbar_{g_1,1+|I_1|} \times \Mbar_{g_2,1+|I_2|} \to \Mbar_{g,n} \,,\\
  &
  \pi^* \bm{1}_{g,n} = \bm{1}_{g,n+1} \,,
  &&
  \pi \colon \Mbar_{g,n+1} \to \Mbar_{g,n} \,.
\end{align}
The first two equations can be interpreted as a cohomological version of the locality axiom in 2D topological field theories (TQFT for short). Taking inspiration from TQFTs, we now define their cohomological analogue based on the cohomology of $\Mbar_{g,n}$. The original definition, due to M.~Kontsevich and Y.~Manin in the mid '90s \cite{KM94}, was the first attempt at axiomatising topological string theory and has deep connections with the seminal work of B.~Dubrovin on the geometry of 2D TQFTs \cite{Dub96}.

\subsection{Axioms}
Fix once and for all a finite-dimensional $\Q$-vector space $V$, called the \emph{phase space}, equipped with a non-degenerate pairing $\eta \colon V \times V \to \Q$. For convenience, we work in a fixed basis $(e_1,\ldots,e_r)$ of $V$. We denote by $(\eta_{\mu,\nu})$ the matrix elements of the pairing, and by $(\eta^{\mu,\nu})$ the inverse matrix.

\begin{definition}
  A \emph{cohomological field theory} on $(V,\eta)$ consists of a collection $\Omega = (\Omega_{g,n})_{2g-2+n>0}$ of linear maps
  \begin{equation}
    \Omega_{g,n} \colon V^{\otimes n} \longrightarrow H^{2\bullet}(\Mbar_{g,n},\Q) \,,
    \qquad\qquad
    \Omega_{g,n}(e_{\mu_1} \otimes \cdots \otimes e_{\mu_n})
    =
    \Omega_{g;\mu_1,\ldots,\mu_n} \,,
  \end{equation}
  satisfying the following axioms.
  \begin{enumerate}
    \item[i)] \textbf{Symmetry.} Each $\Omega_{g,n}$ is $S_n$-invariant, where the action of the symmetric group $S_n$ permutes simultaneously the marked points of $\Mbar_{g,n}$ and the copies of $V^{\otimes n}$.

    \item[ii)] \textbf{Gluing.} Considering the gluing maps
    \begin{equation}
    \begin{aligned}
      \rho \colon& \Mbar_{g-1,n+2} \longrightarrow \Mbar_{g,n} \,, \\
      \sigma \colon& \Mbar_{g_1,1+|I_1|} \times \Mbar_{g_2,1+|I_2|}  \longrightarrow \Mbar_{g,n} \,,
      \qquad
      g_1 + g_2 = g, \; I_1 \sqcup I_2 = \set{1,\dots,n} \,,
    \end{aligned}
    \end{equation}
    we require
    \begin{equation}
    \begin{aligned}
      \rho^{\ast} \Omega_{g;\mu_1,\ldots,\mu_n}
      & =
      \eta^{\alpha,\beta} \; \Omega_{g-1;\alpha,\beta,\mu_1,\ldots,\mu_n}
      \,, \\
      \sigma^{\ast} \Omega_{g;\mu_1,\ldots,\mu_n}
      & =
      \eta^{\alpha,\beta} \; \Omega_{g_1;\alpha,\mu_{I_1}} \otimes \Omega_{g_2;\beta,\mu_{I_2}} \,.
    \end{aligned}
    \end{equation}
  \end{enumerate}
  If the vector space is equipped with a distinguished non-zero element (which may, without loss of generality, be taken to be $e_1$), we can also ask for a third axiom.
  \begin{enumerate}
    \item[iii)] \textbf{Unit.} Considering the forgetful map
    \begin{equation}
      \pi \colon \Mbar_{g,n+1} \longrightarrow \Mbar_{g,n} \,,
    \end{equation}
    we require
    \begin{equation}
      \pi^{\ast} \Omega_{g;\mu_1,\ldots,\mu_n}
      =
      \Omega_{g;\mu_1,\ldots,\mu_n,1}
      \qquad\text{and}\qquad
      \Omega_{0;\mu,\nu,1}
      =
      \eta_{\mu,\nu} \,.
    \end{equation}
  \end{enumerate}
  In this case, $\Omega$ is called a cohomological field theory \emph{with unit}; the distinguished element is called \emph{unit} or \emph{vacuum}.
\end{definition}

Pictorially, the axioms can be illustrated as follows.
\begin{align}
  &\begin{tikzpicture}[baseline]
    \draw[thick] (0,0) -- ($(0,0) + (150:.85)$);
    \node at ($(0,0) + (150:.85)$) [left] {\tiny$\mu_1$};
    \node[rotate=90] at ($(0,0) + (180:.55)$) {\tiny$\cdots$};
    \draw[thick] (0,0) -- ($(0,0) + (-150:.85)$);
    \node at ($(0,0) + (-150:.85)$) [left] {\tiny$\mu_n$};
    \draw[thick,fill=white] (0,0) circle[radius=.35cm];
    \node at (0,0) {\tiny$\Omega_{g}$};
  \end{tikzpicture}
  \;\;
  \overset{\rho^*}{\longmapsto}
  \;\;
  \begin{tikzpicture}[baseline]
    \draw[thick] (0,0) -- ($(0,0) + (150:.85)$);
    \node at ($(0,0) + (150:.85)$) [left] {\tiny$\mu_1$};
    \node[rotate=90] at ($(0,0) + (180:.55)$) {\tiny$\cdots$};
    \draw[thick] (0,0) -- ($(0,0) + (-150:.85)$);
    \node at ($(0,0) + (-150:.85)$) [left] {\tiny$\mu_n$};
    \draw[thick,rotate=-90] (0,0) to[out=45,in=-90] (.4,.6) to[out=90,in=0] (0,1) to[out=180,in=90] (-.4,.6) to[out=-90,in=135] (0,0);
    \draw[thick,fill=white] (0,0) circle[radius=.35cm];
    \node at (0,0) {\tiny$\Omega_{g-1}$};
    \node at (60:.6) {\tiny$\alpha$};
    \node at (-65:.6) {\tiny$\beta$};
    \node at (1,0) [right] {\footnotesize$\eta$};
  \end{tikzpicture} \\
  &\begin{tikzpicture}[baseline]
    \draw[thick] (0,0) -- ($(0,0) + (150:.85)$);
    \node at ($(0,0) + (150:.85)$) [left] {\tiny$\mu_1$};
    \node[rotate=90] at ($(0,0) + (180:.55)$) {\tiny$\cdots$};
    \draw[thick] (0,0) -- ($(0,0) + (-150:.85)$);
    \node at ($(0,0) + (-150:.85)$) [left] {\tiny$\mu_n$};
    \draw[thick,fill=white] (0,0) circle[radius=.35cm];
    \node at (0,0) {\tiny$\Omega_{g}$};
  \end{tikzpicture}
  \;\;
  \overset{\sigma^*}{\longmapsto}
  \;\;
  \begin{tikzpicture}[baseline]
    \draw[thick] (0,0) -- ($(0,0) + (150:.85)$);
    \node[rotate=90] at ($(0,0) + (180:.55)$) {\tiny$\cdots$};
    \draw[thick] (0,0) -- ($(0,0) + (-150:.85)$);
    \node at ($(0,0) + (-180:.7)$) [left] {\footnotesize$\mu_{I_1}$};
    \draw[thick] (1.6,0) -- ($(1.6,0) + (30:.85)$);
    \node[rotate=90] at ($(1.6,0) + (0:.55)$) {\tiny$\cdots$};
    \draw[thick] (1.6,0) -- ($(1.6,0) + (-30:.85)$);
    \node at ($(1.6,0) + (0:.7)$) [right] {\footnotesize$\mu_{I_2}$};
    \draw[thick] (0,0) -- (1.6,0);
    \draw[thick,fill=white] (0,0) circle[radius=.35cm];
    \draw[thick,fill=white] (1.6,0) circle[radius=.35cm];
    \node at (0,0) {\tiny$\Omega_{g_1}$};
    \node at (1.6,0) {\tiny$\Omega_{g_2}$};
    \node at (.8,0) [above] {\footnotesize$\eta$};
    \node at (.5,0) [below] {\tiny$\alpha$};
    \node at (1.1,0) [below] {\tiny$\beta$};
  \end{tikzpicture} \\
  &\begin{tikzpicture}[baseline]
    \draw[thick] (0,0) -- ($(0,0) + (150:.85)$);
    \node at ($(0,0) + (150:.85)$) [left] {\tiny$\mu_1$};
    \node[rotate=90] at ($(0,0) + (180:.55)$) {\tiny$\cdots$};
    \draw[thick] (0,0) -- ($(0,0) + (-150:.85)$);
    \node at ($(0,0) + (-150:.85)$) [left] {\tiny$\mu_n$};
    \draw[thick,fill=white] (0,0) circle[radius=.35cm];
    \node at (0,0) {\tiny$\Omega_{g}$};
  \end{tikzpicture}
  \;\;
  \overset{\pi^*}{\longmapsto}
  \;\;
  \begin{tikzpicture}[baseline]
    \draw[thick] (0,0) -- ($(0,0) + (150:.85)$);
    \node at ($(0,0) + (150:.85)$) [left] {\tiny$\mu_1$};
    \node[rotate=90] at ($(0,0) + (180:.55)$) {\tiny$\cdots$};
    \draw[thick] (0,0) -- ($(0,0) + (-150:.85)$);
    \node at ($(0,0) + (-150:.85)$) [left] {\tiny$\mu_n$};
    \draw[thick] (0,0) -- ($(0,0) + (0:.85)$);
    \node at ($(0,0) + (0:.85)$) [right] {\tiny$1$};
    \draw[thick,fill=white] (0,0) circle[radius=.35cm];
    \node at (0,0) {\tiny$\Omega_{g}$};
  \end{tikzpicture} \\
  &\begin{tikzpicture}[baseline]
    \draw[thick] (0,0) -- ($(0,0) + (150:.85)$);
    \node at ($(0,0) + (150:.85)$) [left] {\tiny$\mu$};
    \draw[thick] (0,0) -- ($(0,0) + (-150:.85)$);
    \node at ($(0,0) + (-150:.85)$) [left] {\tiny$\nu$};
    \draw[thick] (0,0) -- ($(0,0) + (0:.85)$);
    \node at ($(0,0) + (0:.85)$) [right] {\tiny$1$};
    \draw[thick,fill=white] (0,0) circle[radius=.35cm];
    \node at (0,0) {\tiny$\Omega_{0}$};
  \end{tikzpicture}
  =
  \begin{tikzpicture}[baseline]
    \node at ($(0,0) + (150:.85)$) [left] {\tiny$\mu$};
    \node at ($(0,0) + (-150:.85)$) [left] {\tiny$\nu$};
    \draw[thick] ($(0,0) + (150:.85)$) to[out=0,in=90] (-.2,0) to[out=-90,in=0] ($(0,0) + (-150:.85)$);
  \end{tikzpicture} 
\end{align}
A cohomological field theory (CohFT for short) determines a product $\star$ on $V$, called the \emph{quantum product}:
\begin{equation}
  e_{\mu} \star e_{\nu} = \Omega_{0;\mu,\nu,\alpha} \; \eta^{\alpha,\beta} \; e_{\beta} \,.
\end{equation}
Commutativity and associativity of $\star$ follow from (i) and (ii) respectively. If the CohFT comes with a unit, the quantum product is unital, with $e_1 \in V$ being the identity by (iii).

\begin{exercise}
  Prove that $(V,\eta,\star)$ forms a Frobenius algebra, that is, it satisfies
  \begin{equation}
    \eta(v_1 \star v_2,v_3) = \eta(v_1,v_2 \star v_3) \,.
  \end{equation}
  A Frobenius algebra (with unit $e$) is equivalent to a 2D topological field theory $\mc{Z}$ via the following assignments: $\mathcal{Z}(S^1) = V$ for the Hilbert space of states on the circle and
  \begin{equation}
  \begin{split}
    \mc{Z}\Bigl(\;
    \begin{tikzpicture}[scale=.2,baseline={([yshift=-.5ex]current bounding box.center)}]
      \draw [fill=gray,fill opacity=.1] (0,0) arc (-90:270:.3cm and .5cm);
      \draw [fill=gray,fill opacity=.1] (0,2) arc (-90:270:.3cm and .5cm);
      \shade[bottom color=blue!30!black!10,top color=blue!30!black!60] (0,3) to[out=0,in=90] (1.5,1.5) to[out=-90,in=0] (0,0) arc (-90:90:.3cm and .5cm) to[out=0,in=-90] (.5,1.5) to[out=90,in=0] (0,2) arc (-90:90:.3cm and .5cm);
      \draw (0,3) to[out=0,in=90] (1.5,1.5) to[out=-90,in=0] (0,0);
      \draw (0,2) to[out=0,in=90] (.5,1.5) to[out=-90,in=0] (0,1);
    \end{tikzpicture}
    \,\Bigr)
    &=
    \eta \colon V \otimes V \to \Q \,, \\
    \mc{Z}\Bigl(\;
    \begin{tikzpicture}[scale=.2,baseline={([yshift=-.5ex]current bounding box.center)}]
      \draw [fill=gray,fill opacity=.1] (0,0) arc (-90:270:.3cm and .5cm);
      \draw [fill=gray,fill opacity=.1] (0,2) arc (-90:270:.3cm and .5cm);
      \shade[bottom color=blue!30!black!10,top color=blue!30!black!60] (0,0) to[out=0,in=180] (2.5,1) arc (-90:90:.3cm and .5cm) to[out=180,in=0] (0,3) arc (90:-90:.3cm and .5cm) to[out=0,in=90] (.5,1.5) to[out=-90,in=0] (0,1) arc (90:-90:.3cm and .5cm) -- cycle;
      \draw (0,0) to[out=0,in=180] (2.5,1) arc (-90:90:.3cm and .5cm) to[out=180,in=0] (0,3);
      \draw (0,2) to[out=0,in=90] (.5,1.5) to[out=-90,in=0] (0,1);
      \draw [densely dotted] (2.5,1) arc (270:90:.3cm and .5cm);
    \end{tikzpicture}
    \,\Bigr)
    &=
    \star \colon V\otimes V \to V \,, \\
    \mc{Z}\Bigl(
    \begin{tikzpicture}[scale=.3,baseline={([yshift=-.5ex]current bounding box.center)}]
      \shade[bottom color=blue!30!black!10,top color=blue!30!black!60] (0,1) arc (90:-90:.3cm and .5cm) to[out=180,in=-90] (-1,.5) to[out=90,in=180] (0,1) -- cycle;
      \draw (0,1) arc (90:-90:.3cm and .5cm) to[out=180,in=-90] (-1,.5) to[out=90,in=180] (0,1) -- cycle;
      \draw [densely dotted] (0,1) arc (90:270:.3cm and .5cm);
    \end{tikzpicture}
    \,\Bigr)
    &=
    e \colon \Q \to V \,,
  \end{split}
  \end{equation}
  for the morphisms. The partition function $\mc{Z}(\Sigma_{g,n,m})$ of any genus $g$ surface connecting $n$ initial states to $m$ final states can be reconstructed from the above values using the TFT properties.
\end{exercise}

Associated to any CohFT $\Omega$, we also have a collection of rational numbers called \emph{CohFT correlators} (or ancestor invariants), defined as
\begin{equation}\label{eq:CohFT:corr}
  \Braket{\tau_{\mu_1,d_1} \cdots \tau_{\mu_n,d_n}}_g^{\Omega}
  =
  \int_{\Mbar_{g,n}} \Omega_{g;\mu_1,\ldots,\mu_n} \prod_{i=1}^n \psi_i^{d_i} \,.
\end{equation}
Notice that, for degree reasons, $\sum_{i=1}^n d_i \le 3g-3+n$.

\begin{example}\label{expl:CohFTs}
  Here are some examples of CohFTs in one dimension. Let us take $V = \Q \cdot e_1$ and $\eta(e_1,e_1) = 1$. In this case, we use the simpler notation $\Omega_{g,n}$ for $\Omega_{g,n}(e_1^{\otimes n}) = \Omega_{g;1,\ldots,1}$.
  \begin{itemize}
    \item
    Setting $\Omega_{g,n} = \bm{1}_{g,n}$, the unit element in cohomology, we get a CohFT with unit $e_1$ concentrated in degree zero. It is called the \emph{trivial CohFT}, discussed at the beginning of this section. The associated correlators satisfy the Virasoro constraints \eqref{eq:Vir}, equivalent to topological recursion on the Airy spectral curve $\frac{1}{2}y^2 - x = 0$.

    \item
    The class $\Omega_{g,n} = \exp(2\pi^2 \kappa_1)$ defines a CohFT, sometimes called the \emph{Weil--Petersson CohFT} due to its connection with hyperbolic geometry and JT gravity (cf. \cite{Joh,Tur26}). It is not a CohFT with unit. The associated correlators satisfy the Virasoro constraints (a dilaton-shifted version \eqref{eq:Vir}), equivalent to topological recursion on the sine spectral curve (see \cref{ex:WP:CohFT}).

    \item
    The \emph{Hodge class} $\Omega_{g,n} = \Lambda(u) = \sum_{k=0}^{g} \lambda_k \, u^k$ defines a $1$-parameter family of CohFTs with unit $e_1$. It arises as a vertex contribution in the localisation formula for the topological string amplitudes of $\P^1$. A generalisation is provided by a product of Hodge classes:
    \begin{equation}
      \Omega_{g,n}
      =
      \prod_{m=1}^{D} \Lambda(u_m) \,,
    \end{equation}
    which arises as a vertex term in the localisation formula for the topological string amplitudes of a $D$-dimensional spacetime. A particularly nice case is that of $D = 3$ and the parameters $(u_1,u_2,u_3)$ subjected to the constraint
    \begin{equation}
      \frac{1}{u_1} + \frac{1}{u_2} + \frac{1}{u_3} = 0 \,.
    \end{equation}
    In the context of the localisation formulas, this constraint corresponds to the local Calabi--Yau condition \cite{MV02,LLZ03,OP04} (cf. \cite{Liu}). The connection to Virasoro constraints/topological recursion is known only for $D = 1$ and $D = 3$ with the Calabi--Yau condition (see \cref{ex:CY3}).

    \item
    In \cite{Nor23}, Norbury defines a CohFT, denoted as $\Theta_{g,n} \in H^{2(2g-2+n)}(\Mbar_{g,n},\Q)$ and known as \emph{$\Theta$-class}, that satisfies a different version of the unit axiom, namely
    \begin{equation}
      \psi_{n+1} \cdot \pi^{\ast} \Theta_{g,n} = \Theta_{g,n+1} \,.
    \end{equation}
    It appears in super JT gravity in relation to the fermionic part of the Weil--Petersson volumes \cite{SW20}. Norbury conjectured that the associated partition function coincides with the so-called Brézin--Gross--Witten tau-function of the KdV hierarchy \cite{BG80,GW80}, now proved in \cite{CGG25}. Equivalently, it satisfies Virasoro constraints equivalent to topological recursion on the Bessel spectral curve $\frac{1}{2}y^2 x - 1 = 0$.
  \end{itemize}
  Here are some higher-dimensional CohFTs appearing in the literature.
  \begin{itemize}
    \item
    In \cite{Wit92}, Witten studied a generalisation of his original work on 2D quantum gravity by considering a Wess--Zumino--Witten model at level $k$, conveniently re-parametrised as $k = r-2$. Such a theory defines a CohFT of dimension $r-1$, called the \emph{Witten $r$-spin class}, whose basic components are described as follows. Let $V = \bigoplus_{\mu=1}^{r-1} \Q \cdot e_{\mu}$ with pairing $\eta(e_{\mu},e_{\nu}) = \delta_{\mu+\nu,r}$ and unit $e_1$. The \emph{Witten $r$-spin class} is a CohFT
    \begin{equation}
      W_{g;\mu_1,\ldots,\mu_n}^{r} \in H^{2D_{g;\mu}^{r}}(\Mbar_{g,n},\Q)
    \end{equation}
    of pure complex degree
    \begin{equation}
      D_{g;\mu}^{r}
      =
      \frac{(r-2)(g-1) - n + \sum_{i=1}^n \mu_i}{r} \,.
    \end{equation}
    If \smash{$D_{g;\mu}^{r}$} is not an integer, the corresponding Witten class vanishes. The case $r = 2$ gives the trivial cohomological field theory: \smash{$W_{g;1,\ldots,1}^{2} = \bm{1}_{g,n}$}. In genus zero, the construction was first carried out by Witten \cite{Wit93} using $r$-spin structures. The construction of Witten's class in higher genera was first obtained by Polishchuk and Vaintrob \cite{PV00}. The associated partition function is an $r$-KdV tau function \cite{FSZ10} and it satisfies W-constraints equivalent to topological recursion on the $r$-Airy spectral curve $\frac{1}{r}y^r - x = 0$. \\ In \cite{PPZ15} it was shown that all known relations in the so-called tautological ring of $\Mbar_{g,n}$ (the minimal subalgebra of the cohomology of $\Mbar_{g,n}$ stable under pushforwards and pullbacks by tautological maps) are deduced from the Witten $r$-spin class.

    \item
    In \cite{CGG25}, the authors introduced an $r$-spin version of the $\Theta$-class, denoted $\Theta_{g,n}^{r}$ and satisfying properties analogous to those satisfied by Witten's class. It was proved to satisfy W-constraints, equivalent to the topological recursion on the $r$-Bessel spectral curve $\frac{1}{r}y^r x - 1 = 0$. A further generalisation, $\Theta_{g,n}^{r,s}$ depending on a second integer $s \in \set{1, \ldots, r-1}$, was studied in \cite{BCGS} where W-constraints and a generalised version of topological recursion were proved.

    \item
    In \cite{Chi08}, Chiodo defined a generalisation of the Hodge class, called \emph{$\Omega$-class}, which depends on two integers $r \ge 1$ and $s \in \Z$. Let $V = \bigoplus_{\mu = 1}^r \Q \cdot e_{\mu}$ with pairing given by \smash{$\eta(e_{\mu},e_{\nu}) = \frac{1}{r} \delta_{\mu+\nu \equiv 0 \pmod{r}}$}. The $\Omega$-class is a CohFT of mixed cohomological degree:
    \begin{equation}
      \Omega^{r,s}_{g;\mu_1,\ldots,\mu_n} \in H^{2\bullet}(\Mbar_{g,n},\Q) \,.
    \end{equation}
    It is defined as the total Chern class of a (virtual) vector bundle over the moduli space of Riemann surfaces. The case $r = s = 1$ retrieves the Hodge class. The cases $r \ge 2$ and $s = \pm 1$ are related to the Witten and Theta $r$-spin classes respectively.

    \item
    Let $G$ be a complex, simple, simply-connected Lie group with Lie algebra $\mf{g}$. Fix an integer $\ell > 0$ and define $V$ to be the $\Q$-vector space spanned by irreducible representations $e_{\mu}$ of $\mf{g}$ at level $\ell$. Set $\eta(e_{\mu},e_{\nu}) = \delta_{\mu,\nu^{\star}}$, where $\nu^{\star}$ denotes the dual representation, and let $e_1$ be the vector associated with the trivial representation. The \emph{Verlinde bundle} is the vector bundle
    \begin{equation}
      \mc{V}^{\mf{g},\ell}_{g;\mu_1,\ldots,\mu_n} \longrightarrow \Mbar_{g,n}
    \end{equation}
    whose fibres over a smooth Riemann surface are the spaces of non-abelian theta functions. The Chern characters of the Verlinde bundle
    \begin{equation}
      \msc{V}^{\mf{g},\ell}_{g;\mu_1,\ldots,\mu_n}
      \coloneqq
      \ch(\mc{V}^{\mf{g},\ell}_{g;\mu_1,\ldots,\mu_n}) \in H^{2\bullet}(\Mbar_{g,n},\Q)
    \end{equation}
    form a CohFT with unit \cite{MOPPZ17}. The gluing axiom is a consequence of the fusion rules, while the unit axiom is the propagation of vacua.

    \item
    \emph{Topological string amplitudes} on a fixed target Kähler spacetime $(X,\omega)$ are precisely the CohFT correlators of a CohFT with underlying phase space the graded vector space
    \begin{equation}
      V = \bigoplus_{\beta \in H_2(X,\Z)} H^{\bullet}(X,\Z) \cdot q^{\beta} \,,
      \qquad
      \eta(\gamma_1,\gamma_2) = \int_{X} \gamma_1 \frown \gamma_2 \,.
    \end{equation}
    Here, $q$ is a formal variable, known as the Novikov variable, defined as $q^\beta = e^{- \int_\beta \omega}$ and used to grade contributions by the curve class $\beta \in H_2(X,\Z)$ according to its symplectic area. In this context, the phase space $V$ is infinite-dimensional, but graded by $H_2(X,\Z)$, with finite-dimensional components $H^{\bullet}(X,\Z)$. The unit in cohomology $\bm{1} \in H^0(X,\Z)$ serves as the unit for the associated CohFT. This construction was the motivating example for the axiomatic definition of cohomological field theories by Kontsevich and Manin, see \cite{KM94}. 
  \end{itemize}
\end{example}

\subsection{Givental's action}
\label{subsec:Givental}
We have already seen how $\Mbar_{g,n}$ exhibits a recursive boundary structure. A natural question arises: can we exploit such a recursive structure to define or compute CohFTs? The answer is affirmative, and finds its roots in A.~Givental's work \cite{Giv01a,Giv01b} on localisation computations in topological string theory \cite{GP99}. More precisely, Givental defined two actions on CohFTs, the rotation and translation actions.

\subsubsection{Rotation}
For a fixed $(g,n)$, we have a list of all possible stable graphs parametrising the boundary of $\Mbar_{g,n}$. If we are given a CohFT $\Omega$ on $(V,\eta)$, it is natural to decorate all vertices with cohomology classes provided by $\Omega$ to obtain a cohomology class on $\Mbar_{\Gamma}$. For instance:
\begin{equation}
  \begin{tikzpicture}[baseline,scale=.8]
    \draw[thick] (0,0) -- (150:1.5);
    \draw[thick] (0,0) -- (-150:1.5);
      
    \draw[thick,rotate=-90] (0,0) to[out=45,in=-90] (.45,.75) to[out=90,in=0] (0,1.2) to[out=180,in=90] (-.45,.75) to[out=-90,in=135] (0,0);

    \draw[thick,fill=white] (0,0) circle[radius=.3cm];
    \node at (0,0) {\small$2$};

    \node at (125:.6) {\tiny$\nu_1$};
    \node at (-125:.6) {\tiny$\nu_2$};
    \node at (55:.6) {\tiny$\alpha$};
    \node at (-65:.6) {\tiny$\beta$};

    \node at (2.5,0) {$\rightsquigarrow$};

    \node at (5,0) {$\Omega_{2;\nu_1,\nu_2,\alpha,\beta}$};
  \end{tikzpicture}
\end{equation}
In order to produce a cohomology class on $\Mbar_{g,n}$, we must contract all indices at the edges with a cohomology-valued matrix $E^{\nu_h,\nu_{h'}}$ (a priori arbitrary), the indices at the leaves with a cohomology-valued matrix $L_{\mu_i}^{\nu_i}$ (a priori arbitrary), and pushforward the result via the gluing map $\xi_{\Gamma}$. In the above example, we would get
\begin{equation}
  \Omega_{2;\nu_1,\nu_2,\alpha,\beta} \,
  E^{\alpha,\beta} \,
  L_{\mu_1}^{\nu_1} \,
  L_{\mu_2}^{\nu_2} \,,
\end{equation}
where $\mu_i$ denotes a fixed decoration at the $i$-th leaf (i.e. the $i$-th marked point).

Dividing by the natural automorphism factor and summing over all possible stable graphs, we obtain an expression of the form
\begin{equation}\label{eq:naive:rotation}
  \sum_{\Gamma \text{ of type }(g,n)}
    \frac{1}{|\Aut(\Gamma)|} \,
    \xi_{\Gamma,*}
      \left(
      \prod_{v \in V(\Gamma)}
        \Omega_{g(v);(\nu_h)_{h \rightsquigarrow v}}
      \right) \!
      \left(
      \prod_{e = (h,h') \in E(\Gamma)}
        E^{\nu_{h},\nu_{h'}}
      \right) \!
      \left(
      \prod_{i = 1}^{n}
        L_{\mu_i}^{\nu_i}
      \right) .
\end{equation}
Here $h \rightsquigarrow v$ denotes any half-edge $h$ incident to the vertex $v$.

The natural question is: when is the collection of cohomology classes resulting from \eqref{eq:naive:rotation} forming a CohFT? It turns out that \eqref{eq:naive:rotation} is too naive: the matrices $E^{\mu,\nu}$ and $L_{\mu}^{\nu}$ cannot be arbitrary, but should involve specific combinations $\psi$-classes. This condition is captured by a single element called the rotation matrix.

A \emph{rotation matrix} on $(V,\eta)$ is an $\End(V)$-valued power series that is the identity in degree $0$ and satisfying the symplectic condition with respect to $\eta$:
\begin{equation}
  R_{\mu}^{\nu}(u)
  =
  \delta_{\mu}^{\nu} + \sum_{k \ge 1} (R_k)_{\mu}^{\nu} \, u^k
  \in \Q\bbraket{u} \,,
  \qquad
  R^{\mu}_{\alpha}(u) \, \eta^{\alpha,\beta} \, R^{\nu}_{\beta}(-u) = \eta^{\mu,\nu} \,.
\end{equation}
For a given rotation matrix, define the edge decoration as the following $V^{\otimes 2}$-valued power series in two variables\footnote{Beware that several authors use $R^{-1}$ instead of $R$. Here we follow Givental's convention.}:
\begin{equation}
  E^{\mu,\nu}(u,v)
  =
  \frac{\eta^{\mu,\nu} - R^{\mu}_{\alpha}(u) \, \eta^{\alpha,\beta} \, R^{\nu}_{\beta}(v)}{u + v}
  \in
  \Q\bbraket{u,v}.
\end{equation}
The symplectic condition guarantees that $E^{\mu,\nu}(u,v)$ is regular along $u +v = 0$. Define the scalars $E^{\mu,\nu}_{k,\ell}$ through the expansion $E^{\mu,\nu}(u,v) = \sum_{k,\ell \ge 0} E^{\mu,\nu}_{k,\ell} u^k v^{\ell}$.

\begin{definition}
  Consider a CohFT $\Omega$ on $(V,\eta)$ together with a rotation matrix $R$. We define a new collection of cohomology-valued linear maps
  \begin{equation}
    R\Omega_{g,n} \colon V^{\otimes n} \longrightarrow H^{2\bullet}(\Mbar_{g,n},\Q)
  \end{equation}
  as follows. For each stable graph $\Gamma$ of type $(g,n)$, define a contribution through the following construction:
  \begin{itemize}
    \item place $\Omega_{g(v);(\nu_h)_{h \rightsquigarrow v}}$ at each vertex $v$ of $\Gamma$, with arbitrary decorations $\nu_h$ at the half-edges connected to $v$;

    \item place $R_{\mu_i}^{\nu_i}(\psi_{i})$ at the $i$-th leaf of $\Gamma$,

    \item place $E^{\nu_{h},\nu_{h'}}(\psi_{h},\psi_{h'})$ at every edge $e = (h,h')$ of $\Gamma$,

    \item contract all the indices.
  \end{itemize}
  In other words, we get a cohomology class:
  \begin{equation}
    \mr{Cont}_{\Gamma;\mu_1\ldots,\mu_n}
    =
    \left(
      \prod_{v \in V(\Gamma)}
        \Omega_{g(v);(\nu_h)_{h \rightsquigarrow v}}
    \right) \!
    \left(
      \prod_{e = (h,h') \in E(\Gamma)}
        E^{\nu_{h},\nu_{h'}}(\psi_{h},\psi_{h'})
    \right) \!
    \left(
      \prod_{i = 1}^{n}
        R_{\mu_i}^{\nu_i}(\psi_{i})
    \right) .
  \end{equation}
  Although the expressions $E^{\nu_{h},\nu_{h'}}(\psi_{h},\psi_{h'})$ and $R_{\mu_i}^{\nu_i}(\psi_{i})$ have, a priori, infinitely many terms, they terminate due to cohomological degree reasons.

  Define $R\Omega_{g;\mu_1,\ldots,\mu_n}$ to be the sum of contributions of all stable graphs, after pushforward to the moduli space weighted by automorphism factors:
  \begin{equation}
    R\Omega_{g;\mu_1,\ldots,\mu_n}
    =
    \sum_{\Gamma \text{ of type }(g,n)}
      \frac{1}{|\Aut(\Gamma)|} \,
      \xi_{\Gamma,*}
      \mr{Cont}_{\Gamma;\mu_1\ldots,\mu_n} \,.
  \end{equation}
\end{definition}

Let us analyse some examples in low topologies.
\begin{itemize}
  \item
  $R\Omega_{0,3}$. There is a single stable graph of type $(0,3)$, and for dimensional reasons, the decoration $R(\psi_{i})$ at the leaves is simply the identity. Thus, we find
  \begin{equation}
    R\Omega_{0,3}
    =
    \Omega_{0,3}.
  \end{equation}

  \item
  $R\Omega_{0,4}$. The stable graphs of type $(0,4)$ are the following:
  \begin{equation}
    \begin{tikzpicture}[baseline,scale=.5]
      \node at (-2.25,0) {$\Gamma_0 =$};
      \draw[thick] (0,0) -- (135:1);
      \node at (135:1.3) {\tiny$1$};
      \draw[thick] (0,0) -- (-135:1);
      \node at (-135:1.3) {\tiny$2$};
      \draw[thick] (0,0) -- (-45:1);
      \node at (-45:1.3) {\tiny$3$};
      \draw[thick] (0,0) -- (45:1);
      \node at (45:1.3) {\tiny$4$};

      \draw[thick,fill=white] (0,0) circle[radius=.3cm];
      \node at (0,0) {\tiny$0$};

      \begin{scope}[xshift = 10cm]
        \node at (-3.5,0) {$\Gamma_{ij|k\ell} =$};
        \draw[thick] (-1,0) -- (1,0);

        \draw[thick] (-1,0) -- ($(-1,0) + (135:1)$);
        \node at ($(-1,0) + (135:1.3)$) {\tiny$i$};
        \draw[thick] (-1,0) -- ($(-1,0) + (-135:1)$);
        \node at ($(-1,0) + (-135:1.3)$) {\tiny$j$};
        
        \draw[thick] (1,0) -- ($(1,0) + (45:1)$);
        \node at ($(1,0) + (45:1.3)$) {\tiny$k$};
        \draw[thick] (1,0) -- ($(1,0) + (-45:1)$);
        \node at ($(1,0) + (-45:1.3)$) {\tiny$\ell$};

        \draw[thick,fill=white] (-1,0) circle[radius=.3cm];
        \node at (-1,0) {\tiny$0$};
        \draw[thick,fill=white] (1,0) circle[radius=.3cm];
        \node at (1,0) {\tiny$0$};
      \end{scope}
    \end{tikzpicture}
  \end{equation}
  for $ij|k\ell \in \{12|34, 13|24, 14|23\}$.
  The contribution of the stable graph $\Gamma_0$ is given by
  \begin{equation}
  \begin{split}
    \mr{Cont}_{\Gamma_0;\mu_1,\mu_2,\mu_3,\mu_4}
    =
    \Omega_{0;\mu_1,\mu_2,\mu_3,\mu_4}
    &+
    \Omega_{0;\alpha,\mu_2,\mu_3,\mu_4} (R_1)_{\mu_1}^{\alpha} \psi_1
    +
    \Omega_{0;\alpha,\mu_1,\mu_3,\mu_4} (R_1)_{\mu_2}^{\alpha} \psi_2 \\
    &+
    \Omega_{0;\alpha,\mu_1,\mu_2,\mu_4} (R_1)_{\mu_3}^{\alpha} \psi_3
    +
    \Omega_{0;\alpha,\mu_1,\mu_2,\mu_3} (R_1)_{\mu_4}^{\alpha} \psi_4 \,.
  \end{split}
  \end{equation}
  The contribution of the stable graph $\Gamma_{ij|k\ell}$ is given by
  \begin{equation}
    \mr{Cont}_{\Gamma_{ij|kl};\mu_1,\mu_2,\mu_3,\mu_4}
    =
    \Omega_{0;\mu_i,\mu_j,\alpha} \, E^{\alpha,\beta}_{0,0} \, \Omega_{0;\beta,\mu_k,\mu_{\ell}} \,.
  \end{equation}
  It can be shown that $\xi_{\Gamma_{ij|k\ell},*}\bm{1} = [\Gamma_{ij|k\ell}] = \kappa_1$, so that we find
  \begin{equation}
  \begin{split}
    R\Omega_{0;\mu_1,\mu_2,\mu_3,\mu_4}
    =
    \Omega_{0;\mu_1,\mu_2,\mu_3,\mu_4}
    &+
    \Omega_{0;\alpha,\mu_2,\mu_3,\mu_4} (R_1)_{\mu_1}^{\alpha} \psi_1
    +
    \Omega_{0;\alpha,\mu_1,\mu_3,\mu_4} (R_1)_{\mu_2}^{\alpha} \psi_2 \\
    &+
    \Omega_{0;\alpha,\mu_1,\mu_2,\mu_4} (R_1)_{\mu_3}^{\alpha} \psi_3
    +
    \Omega_{0;\alpha,\mu_1,\mu_2,\mu_3} (R_1)_{\mu_4}^{\alpha} \psi_4 \\
    &+
    \left(
      \sum_{ij|k\ell}
        \Omega_{0;\mu_i,\mu_j,\alpha} \, E^{\alpha,\beta}_{0,0} \, \Omega_{0;\beta,\mu_k,\mu_{\ell}}
    \right) \kappa_1
  \end{split}
  \end{equation}

  \item
  $R\Omega_{1,1}$. There are two stable graphs of type $(1,1)$:
  \begin{equation}
    \begin{tikzpicture}[baseline,scale=.5]
      \node at (-3,0) {$\Gamma=$};
      \draw[thick] (0,0) -- (180:1.5);
      \node at (180:1.7) {\tiny$1$};

      \draw[thick,fill=white] (0,0) circle[radius=.3cm];
      \node at (0,0) {\tiny$1$};

      \begin{scope}[xshift = 10cm]
        \node at (-3,0) {$\Gamma'=$};
        \draw[thick] (0,0) -- (180:1.5);
        \node at (180:1.7) {\tiny$1$};
      
        \draw[thick,rotate=-90] (0,0) to[out=45,in=-90] (.45,.75) to[out=90,in=0] (0,1.2) to[out=180,in=90] (-.45,.75) to[out=-90,in=135] (0,0);

        \draw[thick,fill=white] (0,0) circle[radius=.3cm];
        \node at (0,0) {\tiny$0$};

        \node at (1.8,0) {$\vphantom{\Gamma'=}$.};
      \end{scope}
    \end{tikzpicture}
  \end{equation}
  The contribution of $\Gamma$ is
  \begin{equation}
    \mr{Cont}_{\Gamma;\mu} = \Omega_{1;\mu} + \Omega_{1;\nu} \, (R_1)_{\mu}^{\nu} \, \psi_1 \,.
  \end{equation}
  For the one-loop diagram $\Gamma'$, we find
  \begin{equation}
    \mr{Cont}_{\Gamma';\mu} = \Omega_{0;\mu,\alpha,\beta} \, E^{\alpha,\beta}_{0,0} \,.
  \end{equation}
  It can be shown that $\tfrac{1}{2} \xi_{\Gamma',\ast} 1 = [\Gamma'] = 12 \psi_1$, so that
  \begin{equation}
    R\Omega_{1;\mu}
    =
    \Omega_{1;\mu} + \left( \Omega_{1;\nu} \, (R_1)_{\mu}^{\nu} + 12 \, \Omega_{0;\mu,\alpha,\beta} \, E^{\alpha,\beta}_{0,0} \right) \psi_1.
  \end{equation}
\end{itemize}

The main point of this construction is that the resulting collection of cohomology-valued maps $R\Omega$ forms a CohFT.

\begin{proposition}
  The collection of cohomology-valued linear maps $R\Omega = (R\Omega_{g,n})_{2g-2+n > 0}$ forms a CohFT on $(V,\eta)$. Moreover, rotations form a right group action.
\end{proposition}

\subsubsection{Translation}
The rotation action exploits the gluing map by attaching CohFTs through a sort of 2-point correlator, the rotation matrix. There is one more tautological map we can take into account: the forgetful map. Diagrammatically, the forgetful map prunes a leaf of the diagram, which can be decorated (before forgetting it) with a sort of 1-point correlator. As in the case of rotations, the correct approach is to decorate the forgotten leaf with a specific combination of $\psi$-classes. This is taken into account by the translation.

A \emph{translation} is a $V$-valued power series vanishing in degrees $0$ and $1$:
\begin{equation}
  T^{\mu}(u) = \sum_{d \ge 1} (T_d)^{\mu} u^{d+1} \in u^2 \Q\bbraket{u} \,.
\end{equation}

\begin{definition}
  Consider a CohFT $\Omega$ on $(V,\eta)$, together with a translation $T$. We define a collection of cohomology-valued linear maps
  \begin{equation}
    T\Omega_{g,n} \colon V^{\otimes n} \rightarrow H^{2\bullet}(\Mbar_{g,n},\Q)
  \end{equation}
  by setting
  \begin{equation}
    T\Omega_{g;\mu_1,\ldots,\mu_n}
    =
    \sum_{m \ge 0}
      \frac{1}{m!} \pi_{m,\ast}
      \Omega_{g;\mu_1, \ldots, \mu_n,\nu_1,\ldots,\nu_m} \,
      T^{\nu_1}(\psi_{n+1}) \cdots T^{\nu_m}(\psi_{n+m}) \,.
  \end{equation}
  Here $\pi_m \colon \Mbar_{g,n+m} \to \Mbar_{g,n}$ is the map forgetting the last $m$ marked points. Notice that the vanishing of $T$ in degree $0$ and $1$ ensures that the above sum is actually finite.
\end{definition}

\begin{proposition}
  The collection of cohomology-valued linear maps $T\Omega = (T\Omega_{g,n})_{2g-2+n > 0}$ forms a CohFT on $(V,\eta)$. Moreover, translations form an abelian group action.
\end{proposition}

One can also check the composition law for a combination of rotation and translation. The result is parallel to the action of rotation and translation on the plane, hence the name.

\subsubsection{Examples and Teleman's theorem}
Several CohFTs are expressed through Givental's action. We present here the cases of the Weil--Petersson class and the Hodge class.

\begin{exercise}\label{ex:WP:CohFT}
  Prove that $\exp(2\pi^2 \kappa_1)$ is the CohFT obtained from the trivial one under the action of the following translation:
  \begin{equation}
    T(u)
    =
    - \sum_{k \ge 1} \frac{(-2\pi^2)^k}{k!} u^{k+1}
    =
    u\bigl( 1 - e^{-2\pi^2 u} \bigr) \,.
  \end{equation}
\end{exercise}

\begin{theorem}[Mumford's formula]
  The Hodge class $\Lambda(t)$ is the CohFT obtained from the trivial one under the action of the following translation and rotation (in this order) \cite{Mum83}:
  \begin{equation}
  \begin{split}
    R(u)
    &=
    \exp\left(
      - \sum_{m \ge 1} \frac{B_{m+1}}{m(m+1)} (tu)^m
    \right), \\
    T(u)
    &=
    u \bigl( 1 - R(u) \bigr) \,,
  \end{split}
  \end{equation}
  where $B_m$ is the $m$-th Bernoulli number. After re-summing the stable graphs sum, we deduce that
  \begin{equation}\label{eq:Mumford:formula}
    \Lambda(t)
    =
    \exp\Biggl( \sum_{m \ge 1} \frac{B_{m+1}}{m(m+1)} t^m \biggl(
      \kappa_m
      -
      \sum_{i=1}^n \psi_i^m
      +
      \delta_m
    \biggr)\Biggr),
  \end{equation}
  where $\delta_m = \frac{1}{2} \, j_{\ast} \bigl( \sum_{k + \ell = m-1} \psi^{k} (\psi')^{\ell} \bigr)$, and $j$ is the inclusion of all codimension-one boundary strata (i.e. stable graphs with a single edge). The classes $\psi$ and $\psi'$ are the two $\psi$-classes at the nodes.
\end{theorem}

Givental's action is extremely powerful for two reasons. First, as we will see shortly, it gives a recursive way of computing CohFT correlators. Secondly, it might produce relations in cohomology! Take for instance Mumford's formula. One knows from geometric reasons that the Hodge class $\Lambda(t)$ vanishes in degree $d > g$ (it is the Chern polynomial of a rank $g$ bundle). On the other hand, Mumford's formula for $\Lambda(t)$ gives a certain class in any degree. Denoting by $\mc{H}_{g,n}^{d}$ the component of Mumford's formula in complex degree $d$ (i.e. the coefficient of $t^d$ in the right-hand side of \cref{eq:Mumford:formula}), we obtain the following tautological relations: for every $d > g$, $\mc{H}_{g,n}^{d} = 0$ in $H^{2d}(\Mbar_{g,n},\Q)$. The first non-trivial example of such tautological relations is the degree $1$ relation in genus $0$:
\begin{equation}
  \mc{H}_{0,n}^{1}
  =
  \kappa_1 - \sum_{i=1}^n \psi_i + \delta_1 = 0
  \qquad
  \text{in }
  H^{2}(\Mbar_{0,n},\Q) \,.
\end{equation}
Pixton--Pandharipande--Zvonkine \cite{PPZ15} exploited this argument in the case of Witten $3$-spin class to prove all known relations in cohomology.

\begin{exercise}
  Prove, using Mumford's formula, that $\Lambda(t) \Lambda(-t) = 1$. This is sometimes referred to as \emph{Mumford's relation}. Deduce the relations $\lambda_g^2 = 0$.
\end{exercise}

Another reason why Givental's action is extremely valuable is its range of applicability, a result proved by Teleman \cite{Tel12}. Teleman proved that all CohFTs whose underlying quantum product is semisimple are contained in the orbit of the trivial CohFT under the Givental action. Under an additional homogeneity condition, he provided an algorithm to explicitly compute the rotation and the translation matrix.

\begin{theorem}[{Teleman's classification}]
  Let $\Omega$ be a CohFT on $(V,\eta)$. If $\Omega$ is semisimple and homogeneous, then there exist explicit $R$ and $T$ such that
  \begin{equation}
    \Omega_{g,n} = RTw_{g,n} \,,
  \end{equation}
  where $w_{g,n} = \Omega_{g,n}|_{H^0(\Mbar_{g,n},\Q)}$ is the associated 2D TFT. If $\Omega$ is semisimple (but not homogeneous), then there exist $R$ and $T$ such that the above equation holds, but they are defined up to a diagonal ambiguity.
\end{theorem}

In other words, Teleman's theorem classifies all semisimple, homogeneous CohFTs as the orbit under the Givental action of semisimple 2D TFT. Pretty neat!

\subsection{Connection to topological recursion}
Givental's action provides a recursive construction of CohFTs. As the correlators of the trivial CohFT are computed recursively via topological recursion, a natural question arises: is it possible to recursively compute all correlators obtained in the Givental orbit of the trivial CohFT? The answer is affirmative and it beautifully connects to the theory of topological recursion.

Consider a spectral curve $\mc{S} = (\Sigma,x,y,B)$ with $r$ simple ramification points. Choose local coordinates $\zeta_\mu$ around a ramification point $\mu$ such that $x = \zeta_{\mu}^2 + x(\mu)$. Consider the auxiliary functions $\xi^{\mu}$ and the associated meromorphic differentials $\dd\xi^{\mu,k}$, defined as
\begin{equation}
  \xi^\mu(z)
  =
  \int^z \left.\frac{B(w,\cdot)}{\dd\zeta_{\mu}(w)}\right|_{w = \mu},
  \qquad
  \dd\xi^{\mu,k}(z)
  =
  \dd\biggl( \left( - \frac{1}{\zeta_{\mu}} \frac{\dd}{\dd\zeta_{\mu}} \right)^{k} \xi^{\mu}(z) \biggr).
\end{equation}
Set $t^{\mu} = -2 \frac{\dd y(z)}{\dd\zeta_{\mu}(z)} \big|_{z = \mu}$. Define the ($r$-copies of the trivial) CohFT\footnote{
  In the remaining part of this section, we work over $\C$ rather than $\Q$.
} on $V = \bigoplus_{\mu = 1}^r \C \cdot e_{\mu}$ by setting $\eta(e_{\mu},e_{\nu}) = \delta_{\mu,\nu}$ and
\begin{equation}
  w_{g;\mu_1,\ldots,\mu_n} = \frac{\delta_{\mu_1,\ldots,\mu_n}}{(t^{\mu_i})^{2g-2+n}} \,.
\end{equation}
Define the rotation matrix $R$ and the translation $T$ by setting
\begin{align}
  \label{eq:rotation:TR}
  R^{\nu}_{\mu}(u)
  & =
  - \sqrt{\frac{u}{2\pi}} \int_{\gamma_{\nu}}
    e^{-\frac{x - x(\nu)}{2u}} \,
    \dd\xi^\mu
    \,, \\
  \label{eq:translation:TR}
  T^{\mu}(u)
  & =
  \Biggl(
    u \, t^{\mu}
    +
    \frac{1}{\sqrt{2\pi u}}
    \int_{\gamma_{\mu}}
      e^{-\frac{x - x(\mu)}{2u}} \,
      \omega_{0,1}
  \Biggr).
\end{align}
Here $\gamma_{\mu}$ is the formal steepest descent path for $x(z)$ emanating from the ramification point $\mu$; locally it can be taken along the real axis in the $\zeta_{\mu}$-plane. Moreover, the equations are intended as equalities between formal power series in $u$, where on the right-hand side we take an asymptotic expansion as $u \to 0$.

Through the Givental action, we can then define a CohFT
\begin{equation}\label{eq:CohFT:TR}
  \Omega_{g,n} = RTw_{g,n} \colon V^{\otimes n} \longrightarrow H^{2\bullet}(\Mbar_{g,n},\C)
\end{equation}
from the data $(w,R,T)$ through a sum over stable graphs as explained in \cref{subsec:Givental}. The connection with the topological recursion correlators is given by the following theorem \cite{Eyn14b,DOSS14}.

\begin{theorem}[{CohFT/TR correspondence}] \label{thm:Eyn:DOSS}
  Fix a compact spectral curve $\mc{S} = (\Sigma,x,y,\omega_{0,2})$ and define the CohFT $\Omega$ as in \labelcref{eq:CohFT:TR}. Then the topological recursion correlators compute the CohFT correlators:
  \begin{equation}\label{eq:Eyn:DOSS}
    \omega_{g,n}(z_1,\dots,z_n)
    =
    \Braket{\tau_{\mu_1,d_1} \cdots \tau_{\mu_n,d_n}}^{\Omega}_g \;
    \dd\xi^{\mu_1,d_1}(z_1) \cdots \dd\xi^{\mu_n,d_n}(z_n) \,.
  \end{equation}
  Conversely, if we are given a CohFT in the Givental orbit of a semisimple 2D TFT, we can define a (local) spectral curve via \cref{eq:rotation:TR,eq:translation:TR} that computes the correlators as in \cref{eq:Eyn:DOSS}.
\end{theorem}

In a nutshell, the correspondence between CohFTs and topological recursion can be summarised as in \cref{tbl:CohFT:TR}.

\begin{table}
  \centering
  \begin{tabular}{c|c}
    \toprule
      CohFT & Topological recursion \\
    \midrule
      $\dim(V)$ & \# ramification points \\[.5ex]
      trivial CohFT & $\frac{\dd y}{\dd\zeta}$ \\[.5ex]
      translation & $\omega_{0,1}$ \\[.5ex]
      rotation & $\dd \xi$ \\[.5ex]
      edge contribution & $\omega_{0,2}$ \\
    \bottomrule
  \end{tabular}
  \caption{The correspondence between CohFT and topological recursion data.}
  \label{tbl:CohFT:TR}
\end{table}

\begin{exercise}\label{ex:TR:CohFT:WP}
  Show that the CohFT associated with the following spectral curve
  \begin{equation}
    \left(
      \P^1, \;
      x(z) = \frac{z^2}{2},\;
      y(z) = \frac{\sin(2\pi z)}{2\pi}, \;
      \omega_{0,2}(z_1,z_2) = \frac{\dd z_1 \dd z_2}{(z_1 - z_2)^2}
    \right)
  \end{equation}
  is the Weil--Petersson CohFT $\exp(2\pi^2 \kappa_1)$. See also \cref{ex:WP:CohFT}; cf. \cite[section 8.1]{Eyn} for a proof.
\end{exercise}

\begin{exercise}\label{ex:CY3}
  Show that the CohFT associated with the following spectral curve
  \begin{equation}
    \left(
      \P^1, \;
      x(z) = -f \log(z) - \log(1-z),\;
      y(z) = -\log(z), \;
      \omega_{0,2}(z_1,z_2) = \frac{\dd z_1 \dd z_2}{(z_1 - z_2)^2}
    \right) \,.
  \end{equation}
  is the \emph{triple Hodge} \smash{$\Lambda(1)\Lambda(f)\Lambda(-f-1)$}. See \cite[section~7]{Eyn} for a proof. This is the CohFT underlying the (framed) topological vertex \cite{MV02,LLZ03,OP04}, and the topological recursion formula for the triple Hodge class is none other than the BKMP remodelling conjecture for the vertex. The large framing limit recovers the so-called Lambert curve from \cite{BM08} that computes Hurwitz numbers (cf. \cite[section~8.2]{Eyn}).

  \smallskip

  {\small \emph{\faLightbulb \ Hint.} Recall the integral representation of the Euler Beta function
  \begin{equation}
    \mathrm{B}(p,q) = \frac{\Gamma(p)\Gamma(q)}{\Gamma(p+q)} = \int_{0}^1 t^{p-1} (1-t)^{q-1} \, dt
  \end{equation}
  and the asymptotic expansion of the Euler Gamma function
  \begin{equation}
    e^{\frac{1}{v}} \sqrt{2\pi} \,
    \frac{ (-v)^{\frac{1}{v}+a+\frac{1}{2}} }{\Gamma(a - v^{-1})}
    \sim
    \exp\Biggl(
      \sum_{m=1}^{\infty} \frac{B_{m+1}(a)}{m(m+1)} v^{m}
    \Biggr).
  \end{equation}
  Here $B_{m+1}(a)$ are Bernoulli polynomials, and specialise to Bernoulli numbers at both $a = 0$ and $a = 1$: $B_{m+1}(0) = B_{m+1}(1) = B_{m+1}$.
  }
\end{exercise}

We conclude by noting that the correspondence between topological recursion and cohomological field theories has found recent applications across various areas of high-energy physics, including gravity \cite{Ebe24,CEMR24,CEMR25,AC25} and gauge theory \cite{GM}.

\newpage
\section{Further directions}
\label{sec:next}

\subsection{Moduli of hyperbolic surfaces}
In JT gravity, the path integral of the theory is over the space of hyperbolic metrics (rather than the space of complex structures). In other words, the `correct' moduli space is that of \emph{hyperbolic structures}:
\begin{equation}\label{eq:MgnL}
  \Mhyp_{g,n}(L_1,\ldots,L_n)
  =
  \left.
  \Set{ X | \substack{
    \displaystyle\text{$X$ is a hyperbolic surface of genus $g$ } \\[.5ex]
    \displaystyle\text{with $n$ labelled geodesics boundaries} \\[.5ex]
    \displaystyle\text{of lengths $L_1,\ldots,L_n$}
  }}
\middle/ \sim \right.
\end{equation}
where $X \sim X'$ if and only if there exists an isometry from $X$ to $X'$ preserving the labelling of the boundary components.

How is that related to the moduli space of Riemann surfaces? A non-trivial result, which is a consequence of the Riemann uniformisation theorem, is that \smash{$\Mhyp_{g,n}(L)$} is homeomorphic to the moduli space of Riemann surfaces discussed in \cref{sec:moduli}.

\begin{theorem}
  The space \smash{$\Mhyp_{g,n}(L)$} is a smooth real orbifold of dimension $2(3g-3+n)$. Moreover, for all $L \in \R_+^n$, it is homeomorphic (as a smooth real orbifold) to the moduli space of smooth Riemann surfaces:
  \begin{equation}
    \Mhyp_{g,n}(L) \cong \M_{g,n} \,.
  \end{equation}
\end{theorem}

For any fixed $L \in \R_+^n$, the moduli space \smash{$\Mhyp_{g,n}(L)$} is naturally equipped with a symplectic form, called the Weil--Petersson form and denoted $\WP$. In particular, we can define the volumes
\begin{equation}
  \VWP_{g,n}(L)
  =
  \int_{\Mhyp_{g,n}(L)} \frac{\WP^{3g-3+n}}{(3g-3+n)!} \,.
\end{equation}
A toy example of such a structure is the fibration over $\R_+ \ni L$ by spheres $S^2(L)$ of radius $L$. Although all fibres are homeomorphic to $\P^1$, each fibre carries a specific symplectic geometry that depends on the point $L$ on the base. For instance, the area of $S^2(L)$ is $4\pi L^2$. However, we can transfer the particular geometry to $\P^1$ and get an $L$-dependent expression on $\P^1$. For instance, under the isomorphism $S^2(L) \cong \P^1$ provided by stereographic projection, we find that the symplectic form on $S^2(L)$ is mapped to the following polynomial in $L$:
\begin{equation}
  4L^2 \, \Re \frac{\dd z \, \dd\bar{z}}{(1 + |z|^2)^2} \,.
\end{equation}
The analogous result for the Weil--Petersson form and the isomorphism \smash{$\Mhyp_{g,n}(L) \cong \M_{g,n}$} is a result due to Wolpert (for the case $L_i = 0$) and Mirzakhani (for the general case) \cite{Wol85,Mir07+}.

\begin{theorem}
  Under the homeomorphism \smash{$\Mhyp_{g,n}(L) \cong \M_{g,n}$}, the Weil--Petersson form extends as a closed form to $\Mbar_{g,n}$ and defines the cohomology class
  \begin{equation}
    2\pi^2 \kappa_1 + \frac{1}{2} \sum_{i=1}^n L_i^2 \, \psi_i \,.
  \end{equation}
\end{theorem}

An immediate consequence of the above result is that the Weil--Petersson volumes are finite (this was not obvious because \smash{$\Mhyp_{g,n}(L)$} is not compact) and is a symmetric polynomial in boundary lengths squared whose coefficients are intersection numbers involving $\psi$-classes and $\exp(2\pi^2\kappa_1)$:
\begin{equation}
  \VWP_{g,n}(L)
  =
  \sum_{\substack{d_1,\dots,d_n \ge 0 \\ d_1 + \cdots + d_n \le 3g-3+n}}
  \int_{\Mbar_{g,n}}
    e^{2\pi^2\kappa_1}
    \prod_{i=1}^n \psi_i^{d_i} \frac{L_i^{2d_i}}{2^{d_i} d_i!} \,.
\end{equation}
These intersection numbers are precisely in the form of CohFT correlators, and as such can be computed by topological recursion!

\begin{exercise}
  Consider the spectral curve
  \begin{equation}\label{eq:SC:Mirz}
    \left(
      \P^1, \;
      x(z) = \frac{z^2}{2},\;
      y(z) = \frac{\sin(2\pi z)}{2\pi}, \;
      \omega_{0,2}(z_1,z_2) = \frac{\dd z_1 \dd z_2}{(z_1 - z_2)^2}
    \right) \,,
  \end{equation}
  see also \cref{ex:TR:CohFT:WP}. Show that the topological recursion correlators associated with the above spectral curve compute the differential of the Laplace transform of the Weil--Petersson volumes:
  \begin{equation}
    \omega_{g,n}(z_1,\ldots,z_n)
    =
    \dd_{z_1} \cdots \dd_{z_n}
    \left(
      \prod_{i=1}^n \int_{0}^{\infty} \dd L_i \, e^{-z_i L_i}
    \right)
      \VWP_{g,n}(L_1,\ldots,L_n) \,.
  \end{equation}
\end{exercise}

A statement equivalent to the topological recursion (for the volumes rather than their Laplace transform) was proved by M.~Mirzakhani in a remarkable series of papers \cite{Mir07,Mir07+}. Her approach is entirely geometric (rather than algebraic) and is based on the following simple idea due to McShane \cite{McS98}.

\begin{figure}[t]
  \centering
  {
  \tikzset{every picture/.style=thick}
  \begin{tikzpicture}[x=1pt,y=1pt,scale=.4]
    \draw[thin,opacity=.5] (183.7593, 729.9248) -- (189.9928, 733.5251) -- (193.3305, 727.5382);
    \draw[BrickRed, opacity=.5](211.7511, 793.1645) .. controls (216.627, 798.7579) and (226.5504, 797.5348) .. (238.0775, 793.4391) .. controls (249.6046, 789.3434) and (262.7355, 782.375) .. (263.8057, 774.7405);
    \draw[BrickRed](187.257, 723.9633) .. controls (211.8857, 737.2462) and (215.3571, 749.7916) .. (215.348, 758.2338) .. controls (215.3389, 766.676) and (211.8493, 771.0151) .. (210.1373, 776.5968) .. controls (208.4252, 782.1785) and (208.4908, 789.0028) .. (211.7511, 793.1645);
    \draw[BrickRed](263.8057, 774.7405) .. controls (262.57, 770.2811) and (256.0757, 774.9102) .. (251.7821, 780.9872) .. controls (247.4885, 787.0642) and (245.3954, 794.589) .. (248.5761, 808.2784);
    \draw[ForestGreen](220.133, 642.413) .. controls (230.051, 633.704) and (247.0255, 646.852) .. (261.5128, 662.0927) .. controls (276, 677.3333) and (288, 694.6667) .. (290.9638, 711.7898) .. controls (293.9275, 728.913) and (287.855, 745.826) .. (277.369, 748.285);
    \draw[ForestGreen, opacity=.5](277.369, 748.285) .. controls (267.726, 753.964) and (255.863, 744.982) .. (245.2648, 733.1577) .. controls (234.6667, 721.3333) and (225.3333, 706.6667) .. (218.0524, 690.0133) .. controls (210.7715, 673.36) and (205.543, 654.72) .. (220.133, 642.413);
    \draw(176, 672) .. controls (188, 672) and (216, 646) .. (240.6667, 624.3333) .. controls (265.3333, 602.6667) and (286.6667, 585.3333) .. (309.3333, 573.6667) .. controls (332, 562) and (356, 556) .. (464.699, 587.884);
    \draw(460, 612) circle[radius=1, cm={16,-4,0,24,(0,0)}];
    \draw(176, 704) circle[radius=1, cm={16,0,-2,32,(0,0)}];
    \draw(428, 752) ellipse[x radius=12, y radius=20];
    \draw(328, 728) .. controls (349.3333, 714.6667) and (373.3333, 709.3333) .. (400, 712);
    \draw(341.458, 720.808) .. controls (359.1527, 725.6027) and (374.1273, 722.4447) .. (386.382, 711.334);
    \draw(344, 616) .. controls (372, 596) and (408, 616) .. (412, 624);
    \draw(357.035, 609.515) .. controls (368, 620) and (392, 624) .. (405.362, 617.681);
    \draw(260, 808) circle[radius=1, cm={11.3137,0,-8,20,(0,0)}];
    \draw(172.845, 735.916) .. controls (188, 740) and (192, 764) .. (205, 784) .. controls (218, 804) and (240, 820) .. (247.936, 826.087);
    \draw(271.457, 789.268) .. controls (260, 772) and (260, 762) .. (266.6667, 755) .. controls (273.3333, 748) and (286.6667, 744) .. (300.9975, 746.049) .. controls (315.3283, 748.098) and (330.6567, 756.196) .. (350.3601, 761.9667) .. controls (370.0635, 767.7375) and (394.142, 771.181) .. (428.052, 772);
    \draw(456.951, 636.322) .. controls (436, 636) and (424, 654) .. (418, 668.3333) .. controls (412, 682.6667) and (412, 693.3333) .. (414, 704.6667) .. controls (416, 716) and (420, 728) .. (428.934, 732.061);
    \draw(296, 672) .. controls (316, 656) and (326, 654) .. (337, 655) .. controls (348, 656) and (360, 660) .. (368, 672);
    \draw(309.972, 662.04) .. controls (320, 672) and (344, 676) .. (361.221, 664.33);
    \node at (305, 816.5) {$\de_m \Sigma$};
    \node[ForestGreen] at (208.708, 622.152) {$\gamma$};

    \begin{scope}[xshift=13cm]
      \draw[thin,opacity=.5](202.4667, 677.5777) -- (209.1964, 679.7214) -- (211.3373, 673.2011);
      \draw[BrickRed, opacity=.5](256.336, 663.317) .. controls (256, 656) and (252, 650) .. (246, 646.3333) .. controls (240, 642.6667) and (232, 641.3333) .. (220, 641.6667) .. controls (208, 642) and (192, 644) .. (176.61, 650.314);
      \draw[BrickRed](205.017, 671.535) .. controls (220, 676) and (232, 676) .. (241, 674) .. controls (250, 672) and (256, 668) .. (256.336, 663.317);
      \draw[ForestGreen](279.374, 667.219) .. controls (290.525, 667.722) and (291.2625, 685.861) .. (285.9773, 701.866) .. controls (280.692, 717.871) and (269.384, 731.742) .. (260.42, 725.438);
      \draw[ForestGreen, opacity=.5](260.42, 725.438) .. controls (252.539, 719.814) and (254.2695, 707.907) .. (258.108, 694.265) .. controls (261.9465, 680.623) and (267.893, 665.246) .. (279.374, 667.219);
      \draw[ForestGreen](294.233, 652.127) .. controls (302.983, 650.866) and (301.4915, 631.433) .. (293.465, 616.0387) .. controls (285.4385, 600.6445) and (270.877, 589.289) .. (259.472, 595.569);
      \draw[ForestGreen, opacity=.5](259.472, 595.569) .. controls (250.614, 600.51) and (255.307, 616.255) .. (262.4798, 629.9933) .. controls (269.6525, 643.7315) and (279.305, 655.463) .. (294.233, 652.127);
      \draw(189.688, 687.994) .. controls (208, 688) and (236, 708) .. (254.6667, 721.3333) .. controls (273.3333, 734.6667) and (282.6667, 741.3333) .. (297.3333, 748) .. controls (312, 754.6667) and (332, 761.3333) .. (352, 765.6667) .. controls (372, 770) and (392, 772) .. (428.052, 772);
      \draw(194.593, 624.022) .. controls (208, 624) and (230, 612) .. (250.3333, 600.6667) .. controls (270.6667, 589.3333) and (289.3333, 578.6667) .. (310.6667, 570.3333) .. controls (332, 562) and (356, 556) .. (464.699, 587.884);
      \draw(460, 612) circle[radius=1, cm={16,-4,0,24,(0,0)}];
      \draw(192, 656) circle[radius=1, cm={16,0,-2,32,(0,0)}];
      \draw(448, 688) ellipse[x radius=12, y radius=20];
      \draw(428, 752) ellipse[x radius=12, y radius=20];
      \draw(428.493, 732.017) .. controls (412, 728) and (424, 708) .. (447.774, 707.996);
      \draw(447.984, 668) .. controls (428, 664) and (426, 654) .. (430, 647) .. controls (434, 640) and (444, 636) .. (456.329, 636.278);
      \draw(248, 668) .. controls (260, 660) and (276, 656) .. (290, 653) .. controls (304, 650) and (316, 648) .. (336, 660);
      \draw(258.448, 662.348) .. controls (284, 672) and (316, 668) .. (328.024, 655.739);
      \draw(328, 728) .. controls (349.3333, 714.6667) and (373.3333, 709.3333) .. (400, 712);
      \draw(341.458, 720.808) .. controls (359.1527, 725.6027) and (374.1273, 722.4447) .. (386.382, 711.334);
      \draw(344, 616) .. controls (372, 596) and (408, 616) .. (412, 624);
      \draw(357.035, 609.515) .. controls (368, 620) and (392, 624) .. (405.362, 617.681);
      \node[ForestGreen] at (256, 740) {$\gamma$};
      \node[ForestGreen] at (252, 576) {$\gamma'$};
    \end{scope}

    \begin{scope}[xshift=26cm]
      \draw[thin,opacity=.5](207.7156, 659.2374) -- (216.1697, 659.2375) -- (216.2145, 651.1117);;
      \draw[BrickRed, opacity=.5](234.863, 609.1302) .. controls (239.9432, 606.4022) and (244.9454, 611.3676) .. (249.6765, 617.4074) .. controls (254.4076, 623.4472) and (258.8675, 630.5614) .. (261.4032, 639.7206) .. controls (263.9389, 648.8798) and (264.5504, 660.0841) .. (258.448, 662.348);
      \draw[BrickRed](208.1211, 651.3745) .. controls (222.1133, 651.443) and (227.1733, 649.4002) .. (229.9307, 644.5661) .. controls (232.688, 639.7321) and (233.1427, 632.1069) .. (232.6948, 624.9749) .. controls (232.247, 617.8428) and (230.8965, 611.2038) .. (234.863, 609.1302);
      \draw[BrickRed](258.448, 662.348) .. controls (253.6513, 664.6681) and (249.0122, 661.1917) .. (244.4677, 656.9486) .. controls (239.9231, 652.7056) and (235.4732, 647.6958) .. (232.16, 637.9908);
      \draw[ForestGreen](279.374, 667.219) .. controls (290.525, 667.722) and (291.2625, 685.861) .. (285.9773, 701.866) .. controls (280.692, 717.871) and (269.384, 731.742) .. (260.42, 725.438);
      \draw[ForestGreen, opacity=.5](260.42, 725.438) .. controls (252.539, 719.814) and (254.2695, 707.907) .. (258.108, 694.265) .. controls (261.9465, 680.623) and (267.893, 665.246) .. (279.374, 667.219);
      \draw[ForestGreen](294.245, 655.236) .. controls (306.214, 655.877) and (305.107, 633.9385) .. (296.2727, 617.2915) .. controls (287.4385, 600.6445) and (270.877, 589.289) .. (259.472, 595.569);
      \draw[ForestGreen, opacity=.5](259.472, 595.569) .. controls (250.614, 600.51) and (255.307, 616.255) .. (262.9025, 630.1885) .. controls (270.498, 644.122) and (280.996, 656.244) .. (294.245, 655.236);
      \draw(189.688, 687.994) .. controls (208, 688) and (236, 708) .. (254.6667, 721.3333) .. controls (273.3333, 734.6667) and (282.6667, 741.3333) .. (297.3333, 748) .. controls (312, 754.6667) and (332, 761.3333) .. (352, 765.6667) .. controls (372, 770) and (392, 772) .. (428.052, 772);
      \draw(460, 612) circle[radius=1, cm={16,-4,0,24,(0,0)}];
      \draw(448, 688) ellipse[x radius=12, y radius=20];
      \draw(428, 752) ellipse[x radius=12, y radius=20];
      \draw(447.984, 668) .. controls (428, 664) and (426, 654) .. (430, 647) .. controls (434, 640) and (444, 636) .. (456.329, 636.278);
      \draw(258.448, 662.348) .. controls (284, 668) and (302, 672) .. (317.6667, 676) .. controls (333.3333, 680) and (346.6667, 684) .. (359.3333, 690) .. controls (372, 696) and (384, 704) .. (394, 712) .. controls (404, 720) and (412, 728) .. (428, 732);
      \draw(344, 616) .. controls (372, 596) and (408, 616) .. (412, 624);
      \draw(357.035, 609.515) .. controls (368, 620) and (392, 624) .. (405.362, 617.681);
      \draw(312, 716) .. controls (344, 712) and (372, 732) .. (368, 744);
      \draw(324.731, 715.617) .. controls (328, 724) and (352, 736) .. (366.568, 735.717);
      \draw(248, 668) .. controls (260, 660) and (276, 656) .. (290.6667, 655.3333) .. controls (305.3333, 654.6667) and (318.6667, 657.3333) .. (330, 662) .. controls (341.3333, 666.6667) and (350.6667, 673.3333) .. (371.3333, 684.6667) .. controls (392, 696) and (424, 712) .. (448, 708);
      \draw(194.593, 624.022) .. controls (208, 624) and (230, 612) .. (250.3333, 600.6667) .. controls (270.6667, 589.3333) and (289.3333, 578.6667) .. (310.6667, 570.3333) .. controls (332, 562) and (356, 556) .. (464.699, 587.884);
      \draw(192, 656) circle[radius=1, cm={16,0,-2,32,(0,0)}];
      \node[ForestGreen][ForestGreen] at (256, 740) {$\gamma$};
      \node[ForestGreen] at (252, 576) {$\gamma'$};
    \end{scope}
  \end{tikzpicture}
  }
  \caption{
    The geodesic $\alpha_p$ (in red) and some of its possible behaviour, together with the simple closed curve(s) it determines (in green). On the left, the arc $\alpha_p$ intersects the boundary component $\de_m \Sigma$ (\textsc{b}${}_m$-type), and it determines a single simple closed curve $\gamma$. In the two other cases, $\alpha_p$ intersect $\de_1\Sigma$ and itself respectively (\textsc{c}-type), determining two simple closed curves $(\gamma,\gamma')$.}
  \label{fig:orthogeodesic:behaviour}
\end{figure}

Consider a fixed hyperbolic surface $(\Sigma,h)$ with geodesic boundaries. Pick a random (with respect to the hyperbolic measure) point $p \in \de_1\Sigma$ and consider the geodesic $\alpha_p$ starting at $p$ orthogonally to $\de_1\Sigma$. Then one of the following mutually exclusive situations must arise (cf. \cref{fig:orthogeodesic:behaviour}).
\begin{enumerate}
  \item[\textsc{i})]\label{enum:orthogeodesic:i}
  The geodesic $\alpha_p$ intersects $\de_m \Sigma$ for some $m \in \set{2,\dots,n}$, without intersecting itself.

  \item[\textsc{ii})]\label{enum:orthogeodesic:ii}
  The geodesic $\alpha_p$ intersects $\de_1 \Sigma$, without intersecting itself.

  \item[\textsc{iii})]\label{enum:orthogeodesic:iii}
  The geodesic $\alpha_p$ intersects itself.

  \item[\textsc{iv})]\label{enum:orthogeodesic:iiv}
  The geodesic $\alpha_p$ never intersects itself or a boundary component (it spirals indefinitely).
\end{enumerate}
On the one hand, the probability of finding (\textsc{iv}) is zero by a result of Birman--Series. Thus, we simply have that $1 = \P_{\textup{\textsc{i}}} + \P_{\textup{\textsc{ii}}} + \P_{\textup{\textsc{iii}}}$, i.e. the total probability is expressed as the sum of finding (\textsc{i}), (\textsc{ii}), or (\textsc{iii}).
In order to compute such probabilities, Mirzakhani proceeded as follows. Consider the union of $\de_1\Sigma$, the geodesic $\alpha_p$ and, only for case (\textsc{i}), $\de_m \Sigma$.
A sufficiently small neighbourhood of this union of curves is a topological pair of pants. By taking geodesic representatives of the boundary components, we obtain an embedded hyperbolic pair of pants whose geodesic boundary is $(\de_1 \Sigma,\de_m\Sigma,\gamma)$, which we call a \textsc{b}${}_m$-case, or $(\de_1 \Sigma,\gamma,\gamma')$, which we call a \textsc{c}-case (see again \cref{fig:orthogeodesic:behaviour}). Thus, we can write
\begin{equation}
  1 = \sum_{m=2}^n \P_{\textup{\textsc{b}}_m} + \P_{\textup{\textsc{c}}} \,.
\end{equation}
Mirzakhani computed the probabilities $\P_{\textup{\textsc{b}}_m}$ and $\P_{\textup{\textsc{c}}}$ as functions of the hyperbolic lengths of the corresponding pairs of pants, deriving her celebrated identity:
\begin{equation}\label{eqn:Mirzakhani:idnty}
  1
  =
  \sum_{m=2}^n \sum_{\gamma}
    B\bigl( L_1,L_m,\ell(\gamma) \bigr)
  +
  \tfrac{1}{2} \sum_{\gamma,\gamma'}
    C\bigl( L_1,\ell(\gamma),\ell(\gamma') \bigr) \,,
\end{equation}
where $B$ and $C$ are the explicit hyperbolic functions
\begin{equation}\label{eqn:Mirz:BC}
\begin{aligned}
  B(L,L',\ell)
  & =
  1 - \frac{1}{L} \log\left(
    \frac{\cosh(\frac{L'}{2}) + \cosh(\frac{L+\ell}{2})}{\cosh(\frac{L'}{2}) + \cosh(\frac{L-\ell}{2})}
  \right), \\
  C(L,\ell,\ell')
  & =
  \frac{2}{L} \log\left(
    \frac{e^{\frac{L}{2}} + e^{\frac{\ell+\ell'}{2}}}{e^{-\frac{L}{2}} + e^{\frac{\ell+\ell'}{2}}}
  \right).
\end{aligned}
\end{equation}
Integration of the constant function $1$ over the moduli space gives the Weil--Petersson volumes on the left-hand side, while the right-hand side can be expressed as a specific integration formula involving volumes of lower topological complexity thanks to the removal of pairs of pants.

\begin{theorem}[{Mirzakhani's recursion}]
  The Weil--Petersson volumes are uniquely determined by the following recursion on $2g - 2 + n > 1$
  \begin{multline}\label{eqn:TR:WP}
    \VWP_{g,n}(L_1,\dots,L_n)
    =
    \sum_{m=2}^n \int_{0}^{\infty}
      \dd\ell \, \ell \, B(L_1,L_m,\ell) \,
      \VWP_{g,n-1}(\ell,L_2\dots,\widehat{L_m},\dots,L_n) \,
      \\
    +
    \frac{1}{2} \int_{0}^{\infty} \int_{0}^{\infty}
      \dd\ell \, \dd\ell' \, \ell \, \ell' \,
      C(L_1,\ell,\ell') \Biggl(
      \VWP_{g-1,n+2}(\ell,\ell',L_2,\dots,L_n) \\
      +
      \sum_{\substack{g_1 + g_2 = g \\ I_1 \sqcup I_2 = \set{2,\dots,n}}}
      \VWP_{g_1,1+|I_1|}(\ell,L_{I_1}) \,
      \VWP_{g_2,1+|I_2|}(\ell',L_{I_2})
    \Biggr)
  \end{multline}
  with the conventions $\VWP_{0,1} = \VWP_{0,2} = 0$ and the base cases $\VWP_{0,3}(L_1,L_2,L_3) = 1$, $\VWP_{1,1}(L) = \frac{L^2}{48} + \frac{\pi^2}{12}$.
\end{theorem}

Topological recursion on the spectral curve given in \labelcref{eq:SC:Mirz} produces precisely the Laplace transform of Mirzakhani's recursion \cite{EO}.

\begin{table}
  \footnotesize
  \begin{center}
  \begin{tabular}{c | l}
    \toprule
    $(g,n)$ &
      $V_{g,n}^{\textup{WP}}(L_1,\ldots,L_n)$ \\
    \midrule
    $(0,3)$ &
      $1$
    \\[.5ex]
    $(0,4)$ &
      $\tfrac{1}{2} m_{(1)} + 2\pi^2$
    \\[.5ex]
    $(0,5)$ &
      $\tfrac{1}{8} m_{(2)} +
      \tfrac{1}{2} m_{(1^2)} +
      3\pi^2 m_{(1)} +
      10\pi^4$
    \\[.5ex]
    $(0,6)$ &
      $\tfrac{1}{48} m_{(3)} +
      \tfrac{3}{16} m_{(2,1)} +
      \tfrac{3}{4} m_{(1^3)} +
      \tfrac{3\pi^2}{2} m_{(2)} +
      6\pi^2 m_{(1^2)} +
      26\pi^4 m_{(1)} +
      \tfrac{244\pi^6}{3}$
    \\[.5ex]
    $(0,7)$ &
      $\begin{aligned}[t]
        &
        \tfrac{1}{384} m_{(4)} +
        \tfrac{1}{24} m_{(3,1)} +
        \tfrac{3}{32} m_{(2^2)} +
        \tfrac{3}{8} m_{(2,1^2)} +
        \tfrac{3}{2} m_{(1^4)} +
        \tfrac{5\pi^2}{12} m_{(3)} +
        \tfrac{15\pi^2}{12} m_{(2,1)} \\
        & \quad
        +
        15\pi^2 m_{(1^3)} +
        20\pi^4 m_{(2)} +
        80\pi^4 m_{(1^2)} +
        \tfrac{910\pi^6}{3} m_{(1)} +
        \tfrac{2758\pi^8}{3}
      \end{aligned}$
    \\[.5ex]
    \midrule
    $(1,1)$ &
      $\tfrac{1}{48} m_{(1)} +
      \tfrac{\pi^2}{12}$
    \\[.5ex]
    $(1,2)$ &
      $\tfrac{1}{192} m_{(2)} +
      \tfrac{1}{96} m_{(1^2)} +
      \tfrac{\pi^2}{12} m_{(1)} +
      \tfrac{\pi^4}{4}$
    \\[.5ex]
    $(1,3)$ &
      $\tfrac{1}{1152} m_{(3)} +
      \tfrac{1}{192} m_{(2,1)} +
      \tfrac{1}{96} m_{(1^3)} +
      \tfrac{\pi^2}{24} m_{(2)} +
      \tfrac{\pi^2}{8} m_{(1^2)} +
      \tfrac{13\pi^4}{24} m_{(1)} +
      \tfrac{14\pi^6}{9}$
    \\[.5ex]
    $(1,4)$ &
      $\begin{aligned}[t]
        &
        \tfrac{1}{9216} m_{(4)} +
        \tfrac{1}{768} m_{(3,1)} +
        \tfrac{1}{384} m_{(2^2)} +
        \tfrac{1}{128} m_{(2,1^2)} +
        \tfrac{1}{64} m_{(1^4)} +
        \tfrac{7\pi^2}{576} m_{(3)} \\
        & \quad
        +
        \tfrac{\pi^2}{12} m_{(2,1)} +
        \tfrac{\pi^2}{4} m_{(1^3)} + 
        \tfrac{41\pi^4}{96} m_{(2)} +
        \tfrac{17\pi^4}{12} m_{(1^2)} +
        \tfrac{187\pi^6}{36} m_{(1)} +
        \tfrac{529\pi^8}{36}
      \end{aligned}$
    \\[.5ex]
    \midrule
    $(2,1)$ &
      $\tfrac{1}{442368} m_{(4)} +
        \tfrac{29\pi^2}{138240} m_{(3)} +
        \tfrac{139\pi^4}{23040} m_{(2)} +
        \tfrac{169\pi^6}{2880} m_{(1)} +
        \tfrac{29\pi^8}{192}$
    \\[.5ex]
    $(2,2)$ &
      $\begin{aligned}[t]
        &
        \tfrac{1}{4423680} m_{(5)} +
        \tfrac{1}{294912} m_{(4,1)} +
        \tfrac{29}{2211840} m_{(3,2)} +
        \tfrac{11\pi^2}{276480} m_{(4)} +
        \tfrac{29\pi^2}{69120} m_{(3,1)} +
        \tfrac{7\pi^2}{7680} m_{(2^2)} \\
        & \quad
        +
        \tfrac{19\pi^4}{7680} m_{(3)} +
        \tfrac{181\pi^4}{11520} m_{(2,1)} +
        \tfrac{551\pi^6}{8640} m_{(2)} +
        \tfrac{7\pi^6}{36} m_{(1^2)} +
        \tfrac{1085\pi^8}{1728} m_{(1)} +
        \tfrac{787\pi^{10}}{480}
      \end{aligned}$
    \\[.5ex]
    \midrule
    $(3,1)$ &
      $\begin{aligned}[t]
        &
        \tfrac{1}{53508833280} m_{(7)} +
        \tfrac{77\pi^2}{9555148800} m_{(6)} +
        \tfrac{3781\pi^4}{2786918400} m_{(5)} +
        \tfrac{47209\pi^6}{418037760} m_{(4)} +
        \tfrac{127189\pi^8}{26127360} m_{(3)} \\
        & \quad
        +
        \tfrac{8983379\pi^{10}}{87091200} m_{(2)} +
        \tfrac{8497697\pi^{12}}{9331200} m_{(1)} +
        \tfrac{9292841\pi^{14}}{4082400}
      \end{aligned}$
    \\[.5ex]
    \bottomrule
  \end{tabular}
  \caption{
    A list of Weil--Petersson polynomials $\VWP_{g,n}(L)$ computed via topological recursion. Here $m_{\lambda}$ is the monomial symmetric polynomial associated with the partition $\lambda$, evaluated at $L_1^2,\dots,L_n^2$.}
  \label{table:WP:poly}
  \end{center}
\end{table}

\subsection{String theory and moduli of maps}
As mentioned in the text, topological string theory is also intimately connected to the moduli space of Riemann surfaces. Topological string theory (or, in mathematical terms, Gromov--Witten theory) aims to compute worldsheets of the strings in a fixed target spacetime $X$ as parametrised Riemann surfaces, that is, maps
\begin{equation}
  f \colon (\Sigma,p_1,\ldots,p_n) \longrightarrow X \,.
\end{equation}
Here $p_1,\ldots,p_n$ are marked points on $\Sigma$ and can be thought of as the initial/final states of the worldsheet $\Sigma$. The path integral of the theory is then an integral over the moduli space of such maps:
\begin{equation}
  \M_{g,n}(X,\beta) = \left. \Set{
    (\Sigma,p_1,\dots,p_n,f) | \substack{
      \displaystyle f \colon (\Sigma,p_1,\ldots,p_n) \to X  \\[.5ex]
      \displaystyle \text{$f_{\ast}[\Sigma] = \beta$}
      }
    }
    \middle/ \sim \right. ,\\[1ex]
\end{equation}
where $\beta \in H_2(X,\Z)$ is a fixed class (called the degree). The proper definition of $\M_{g,n}(X,\beta)$ and its compactification is a highly delicate mathematical problem (much more intricate than the moduli space of Riemann surfaces). The computation of the associated correlators is even more subtle.

Witten's conjecture can be seen as the tip of the iceberg of this theory: it corresponds to the case of $X = \set{*}$, a zero-dimensional target. Eguchi, Hori, and Xiong \cite{EHX97} extended the Virasoro constraints for the point and conjectured that the partition function of every target obeys the Virasoro conditions. In a remarkable series of papers \cite{OP06a,OP06b,OP06c}, Okounkov and Pandharipande gave a complete solution in the one-dimensional case, proving the conjecture of Eguchi--Hori--Xiong.
Beyond the case of a point and of complex curves, Virasoro constraints have also been shown to hold for special classes of targets (of arbitrary dimension), namely:
\begin{itemize}
  \item for toric Fano manifolds and manifolds satisfying a semisimplicity assumption, as shown by Givental--Teleman \cite{Giv01a,Tel12},

  \item even more explicitly for toric Calabi--Yau threefolds following the Bouchard--Klemm--Mariño--Pasquetti ``remodelling conjecture'' \cite{BKMP09}, now a theorem \cite{EO15,FLZ20}.
\end{itemize}

\newpage
\printbibliography

\end{document}